\documentclass[12pt, a4paper]{amsart}
\usepackage{amsmath}
\usepackage{geometry,amsthm,graphics,tabularx,amssymb,shapepar}
\usepackage{latexsym,amsfonts,amssymb,amsmath,amsxtra,bbm,color,mathrsfs}
\usepackage{amscd}

\usepackage{hyperref}
\usepackage{pdfsync}
\usepackage[all]{xypic}

\usepackage{verbatim}

\newcommand{\BA}{{\mathbb {A}}}

\newcommand{\BC}{{\mathbb {C}}}

\newcommand{\BN}{{\mathbb {N}}}

\newcommand{\BQ}{{\mathbb {Q}}}
\newcommand{\BR}{{\mathbb {R}}}

\newcommand{\BZ}{{\mathbb {Z}}}

\newcommand{\CC}{{\mathcal {C}}}

\newcommand{\CE}{{\mathcal {E}}}

\newcommand{\CG}{{\mathcal {G}}}
\newcommand{\CH}{{\mathcal {H}}}

\newcommand{\CM}{{\mathcal {M}}}

\newcommand{\CO}{{\mathcal {O}}}
\newcommand{\CP}{{\mathcal {P}}}

\newcommand{\CS}{{\mathcal {S}}}

\newcommand{\CX}{{\mathcal {X}}}

\newcommand{\RB}{{\mathrm {B}}}

\newcommand{\RH}{{\mathrm {H}}}

\newcommand{\RN}{{\mathrm {N}}}
\newcommand{\RO}{{\mathrm {O}}}

\newcommand{\RU}{{\mathrm {U}}}

\newcommand{\Aut}{{\mathrm{Aut}}}

\newcommand{\GL}{{\mathrm{GL}}}

\newcommand{\Hom}{{\mathrm{Hom}}}

\newcommand{\Ind}{{\mathrm{Ind}}}

\newcommand{\Ker}{{\mathrm{Ker}}}

\renewcommand{\Re}{{\mathrm{Re}}}

\newcommand{\rk}{{\mathrm{k}}}

\newcommand{\sgn}{\operatorname{sgn}}

\newcommand{\od}{\operatorname{d}}

\newcommand{\oL}{\operatorname{L}}

\newcommand{\oH}{\operatorname{H}}

\newcommand{\oZ}{\operatorname{Z}}

\newcommand{\g}{\mathfrak g}

\renewcommand{\k}{\mathfrak k}

\renewcommand{\l}{\mathfrak l}

\renewcommand{\rk}{\mathrm k}

\newcommand{\Z}{\mathbb{Z}}
\newcommand{\C}{\mathbb{C}}
\newcommand{\R}{\mathbb R}

\newcommand{\K}{\mathbb{K}}

\newcommand{\A}{\mathbb{A}}

\newcommand{\ve}{{\vee}}

\newcommand{\abs}[1]{\lvert#1\rvert}

\newcommand{\la}{\langle}
\newcommand{\ra}{\rangle}

\newcommand{\be}{\begin {equation}}
\newcommand{\ee}{\end {equation}}
\newcommand{\bee}{\begin {equation*}}
\newcommand{\eee}{\end {equation*}}

\newcommand{\cf}{\emph{cf.}~}

\theoremstyle{Theorem}

\theoremstyle{Theorem}

\theoremstyle{Theorem}

\theoremstyle{Theorem}
\newtheorem{prp}{Proposition}[section]
\newtheorem{corp}[prp]{Corollary}
\newtheorem{lemp}[prp]{Lemma}
\newtheorem{thmp}[prp]{Theorem}

\theoremstyle{Plain}

\newtheorem{remarkp}[prp]{Remark}

\theoremstyle{Definition}

\newtheorem{dfnp}[prp]{Definition}

\begin{document}

\title[Period relations]{Period relations for Rankin-Selberg convolutions for $\GL(n)\times \GL(n-1)$}

\author[J.-S. Li]{Jian-Shu Li}
\address{Institute for Advanced Study in Mathematics, Zhejiang University\\
  Hangzhou, 310058, China}\email{jianshu@zju.edu.cn}

\author[D. Liu]{Dongwen Liu}
\address{School of Mathematical Sciences,  Zhejiang University\\
  Hangzhou, 310058, China}\email{maliu@zju.edu.cn}

\author[B. Sun]{Binyong Sun}
\address{Institute for Advanced Study in Mathematics, Zhejiang University\\
  Hangzhou, 310058, China}\email{sunbinyong@zju.edu.cn}




\begin{abstract}
We formulate and prove the  archimedean period relations for Rankin-Selberg convolutions for $\GL(n)\times \GL(n-1)$. As a consequence, we prove  the period relations for critical values of the Rankin-Selberg L-functions for $\GL(n)\times \GL(n-1)$ over arbitrary number fields.
\end{abstract}

 \maketitle

\tableofcontents

\section{Introduction}\label{s1}


The cases of $\GL(n)\times \GL(n-1)$ and $\GL(n)\times \GL(n)$ are fundamental in  the general Rankin-Selberg theory, and many problems for general Rankin-Selberg convolutions are reduced to these two cases. The goal of this article is to give an unconditional proof of the period relations for critical values of Rankin-Selberg L-functions for $\GL(n)\times \GL(n-1)$ over arbitrary number fields, which is a long-standing problem and has been studied by many authors. 
In the framework of Langlands program, it is compatible with the celebrated conjecture of  Deligne \cite{D} on the rationality of critical values 
of L-functions attached to pure motives. More general conjectures concerning  period relations for  critical values of  Rankin-Selberg L-functions are formulated by Blasius in \cite{B}. 

 \subsection{Whittaker periods} 
 Let $\rk$ be a number field, and write $\A$ for the adele ring of $\rk$.  Denote by $\rk_v$ the completion of $\rk$ at a place $v$. Write
\[
\rk_\infty:=\rk\otimes_{\BQ}\BR= \prod_{v|\infty} \rk_v \hookrightarrow \rk \otimes_{\BQ}\BC = \prod_{\iota\in \CE_\rk}\BC,
\]
where $\CE_\rk$ is the set of field embeddings $\iota: \rk\hookrightarrow\BC$. 

Let $\Pi$ be an irreducible smooth automorphic representation of $\GL_n(\A)$ ($n\geq 1$) that is cuspidal or (more generally) tamely isobaric as defined in \eqref{uni-iso00}. 
Suppose that $\Pi$ is regular algebraic in the sense of Clozel (see \cite{Clo}).  
By \cite[Section 3]{Clo}, up to isomorphism there is a unique  irreducible algebraic representation $F_\mu$ of $\GL_n(\rk\otimes_\BQ\BC)$, say of  highest weight $\mu = \{\mu^\iota\}_{\iota\in \CE_\rk}  \in (\BZ^n)^{\CE_\rk}$, such that the total continuous cohomology 
\be\label{cohospace}
\RH^*_{\rm ct}(\GL_n(\rk_\infty)^0; F_\mu^\vee\otimes \Pi_\infty)\neq \{0\}.
\ee
Here $\Pi_\infty:=\widehat \otimes_{v|\infty} \Pi_v$ is the infinite part of $\Pi$, a superscript ``$\,^\vee$" over a representation indicates the contragredient representation, and a superscript ``$0$" over a Lie group indicates the identity connected component of the Lie group. Moreover $\mu$ is pure in the sense that there exists $w_\mu\in \BZ$ such that 
\[
\mu^{\iota}_1+\mu^{\bar{\iota}}_n=\mu^{\iota}_2+\mu^{\bar{\iota}}_{n-1}=\dots = \mu^{\iota}_n+\mu^{\bar{\iota}}_1=w_\mu
\]
for all $\iota\in \CE_\rk$. Here we write $\mu^\iota= (\mu^\iota_1,\ldots, \mu^\iota_n)$, and $\bar \iota$ is the composition of 
\[
\rk\xrightarrow{\iota}\C\xrightarrow{\textrm{complex conjugation}}\C.
\]
The representation $F_\mu$ is called the coefficient system of $\Pi$.

Let $\Pi_f:=\otimes'_{v\nmid\infty}\Pi_v$ be the finite part of $\Pi$.  The rationality field $\BQ(\Pi)$ of $\Pi$ is the fixed field of the group of field automorphisms $\sigma\in {\rm Aut}(\BC)$ such that
${}^\sigma (\Pi_f) = \Pi_f$. This is a number field contained in $\C$. 
By \cite[Theorem 3.13]{Clo} and \cite[Lemma 1.2]{G}, for every $\sigma\in \Aut(\BC)$, there exists a unique  irreducible smooth  automorphic representation ${}^\sigma \Pi$ of $\GL_n(\A)$ that is tamely isobaric and regular algebraic, and whose finite part $({}^\sigma\Pi)_f$ is isomorphic to ${}^\sigma(\Pi_f) $.   See Section \ref{sec4.3} for more details.  

The cohomology space in \eqref{cohospace} is naturally a representation of the group
\[
\pi_0(\GL_n(\rk_\infty)):=\GL_n(\rk_\infty)/\GL_n(\rk_\infty)^0.
\]
By using the determinant homomorphism, the latter group is identified with 
\[
\pi_0(\rk_\infty^\times):=\rk_\infty^\times/(\rk_\infty^\times)^0=\{\pm 1\}^{\CE_\rk^\R},
\]
where $\CE_\rk^\R$ denotes the set of real places of $\rk$, which is  identified with a subset of $\CE_\rk$. 
The group of characters of $\pi_0(\rk_\infty^\times)$ is denoted by $\widehat{\pi_0(\rk^\times_\infty)}$, which is obviously identified with the group of quadratic characters of $\rk_\infty^\times$.

For every archimedean local field $\K$, put
\[
  b_{n,\K}:=\left\{
         \begin{array}{ll}
           {\lfloor}\frac{n^2}{4}{\rfloor}, & \hbox{if $ \ \ \K\cong \R$;} \\
           \frac{n(n-1)}{2}, & \hbox{if $ \ \ \K\cong \C$.}
         \end{array}
       \right.
 \]
 Write 
 \[
  b_{n,\infty}:=\sum_{v|  \infty} b_{n,\rk_v}.
 \]
Let $\varepsilon_{\Pi_\infty}$ denote the central character of $F_\mu^\vee\otimes \Pi_\infty$. Note that $\varepsilon_{\Pi_\infty}$ is a quadratic character of $\rk_\infty^\times$, and is trivial when $n$ is even.
By \cite[Lemma 3.14]{Clo},
\[
  \oH_{\mathrm{ct}}^i(\GL_n(\rk_\infty)^0;
                           F_\mu^\vee\otimes \Pi_\infty)=\{0\},\qquad \textrm{if $\ \ i<b_{n,\infty}$,}
\]
and as a representation of $\pi_0(\GL_n(\rk_\infty))$, 
\be\label{botcoh}
  \oH_{\mathrm{ct}}^{b_{n,\infty}}(\GL_n(\rk)^0;      F_\mu^\vee\otimes \Pi_\infty)\cong 
  \left\{
                \begin{array}{ll}
                  {\bigoplus}_{\varepsilon\in \widehat{\pi_0(\rk^\times_\infty)}}\, \varepsilon, & \hbox{if  $n$ is even;}\\
                  \varepsilon_{\Pi_\infty}, & \hbox{if $n$ is odd.} 
                \end{array}
              \right.
\ee
We are particularly interested in the bottom degree cohomology space  \eqref{botcoh}. 

For every $\varepsilon\in \widehat{\pi_0(\rk^\times_\infty)}$ that occurs in the bottom degree cohomology space  \eqref{botcoh},  by comparing the Betti and de Rham cohomologies of the (tower of) locally symmetric spaces attached to $\GL_n(\BA)$, 
Raghuram and Shahidi define a nonzero complex
number, to be called the Whittaker period for $\Pi$ and $\varepsilon$ (see \cite[Definition/Proposition 3.3]{RS2}).
These Whittaker periods play an important role in the arithmetic study of special values of Rankin-Selberg L-functions. However, their definition of Whittaker period is not canonical since it depends on an  arbitrarily fixed  generator of  the $\varepsilon$-eigenspace of \eqref{botcoh}. In Section \ref{sec4}, based on the nonvanishing hypothesis that is proved in \cite{Sun}, we will canonically define Raghuram-Shahidi's  Whittaker period by fixing a canonical generator of the concerning  $\varepsilon$-eigenspace.

With a slight variation, we define the Whittaker period $\Omega_{\varepsilon}(\Pi)$ for every $\varepsilon\in \widehat{\pi_0(\rk^\times_\infty)}$ that occurs in
\[
\CH(\Pi_\infty):=  \oH_{\mathrm{ct}}^{b_{n,\infty}}(\GL_n(\rk)^0;      F_\mu^\vee\otimes \Pi_\infty)\otimes \widetilde{\mathfrak{O}}_{n,\infty},  
\]
where $\widetilde{\mathfrak{O}}_{n,\infty}$ is a certain one-dimensional complex vector space defined by orientations (see \eqref{wtildefo}), which is naturally a representation of $\pi_0(\rk^\times_\infty)$ that is isomorphic to 
$ \sgn_\infty^{\frac{(n-1)(n-2)}{2}}
$. Here $\sgn_\infty$ is the quadratic character of $\rk_\infty^\times$ that is nontrival on $\rk_v^\times$ for every real place $v$ of $\rk$. Note that the isomorphism class of the representation 
$\CH((\,^\sigma \Pi)_\infty)$ of $\pi_0(\rk^\times_\infty)$ is independent of $\sigma\in \Aut(\C)$ (see Remark \ref{remks54}). 

In fact, by fixing a generator of a certain one-dimensional $\mathbb Q(\Pi)$-vector space, we will simultaneously define a family 
\[
\{\Omega_{\varepsilon}(\Pi')\}_{\Pi'\in \{\,^\sigma\Pi\, : \,  \sigma\in \Aut(\C)\}}
\]
of Whittaker periods, which are nonzero complex numbers. Moreover, the family is unique up to scalar multiplication by  $\mathbb Q(\Pi)^\times$ in the following sense  (see Lemma \ref{uniquefamily}): Suppose that another generator yields another family  $\{\Omega'_{\varepsilon}(\Pi')\}_{\Pi'\in \{\,^\sigma\Pi\, : \,  \sigma\in \Aut(\C)\}}$ of Whittaker periods. Then for all $\Pi_1, \Pi_2\in \{\,^\sigma\Pi\, : \,  \sigma\in \Aut(\C)\}$ and all $\sigma\in \Aut(\C)$ such that $\,^\sigma\Pi_1=\Pi_2$,
\[
  \sigma\left(\frac{\Omega'_{\varepsilon}(\Pi_1)}{\Omega_{\varepsilon}(\Pi_1)}\right)=\frac{\Omega'_{\varepsilon}(\Pi_2)}{\Omega_{\varepsilon}(\Pi_2)}.
\]
In particular, like Deligne's periods for pure motives, the Whittaker period $\Omega_\varepsilon(\Pi)$ is uniquely defined up to scalar multiplication by  $\mathbb Q(\Pi)^\times$. See Section \ref{sec4.4} for details. When $n=1$, the Whittaker period $\Omega_\varepsilon(\Pi)\in \BQ(\Pi)^\times$. 


\begin{remarkp}
By comparing Deligne's conjecture and the global period relation (Theorem \ref{thm: global period relation} of this article), Hara and Namikawa  (\cite[Theorem 1.1]{HN2}) supply a conjectural description of our Whittaker period $\Omega_\varepsilon(\Pi)$ in terms of Deligne's periods and Yoshida's fundamental periods (see \cite{Y}). They also partially prove Theorem \ref{thm: global period relation} in \cite[Theorem 6.11]{HN2} under the assumption that  \cite[Conjecture 6.8]{HN2}  holds true. 
We remark that the archimedean period relation (Theorem \ref{thmap}) that is proved in this article  implies their Conjecture 6.8.

\end{remarkp}

\subsection{Period relations}
Suppose that $n\geq 2$ and  $\Pi$  is cuspidal. Let $\Sigma$ be an irreducible smooth automorphic representation of $\GL_{n-1}(\BA)$ that is tamely isobaric and regular algebraic.  Assume that the coefficient systems $F_\mu$ and $F_\nu$ of $\Pi$ and $\Sigma$ respectively are balanced, that is, there is an integer $j$ such that
\[
\Hom_{\GL_{n-1}(\rk\otimes_\BQ\BC)}(F_\mu^\vee\otimes F_\nu^\vee, \otimes_{\iota\in\CE_\rk} {\det}^{j})\neq \{0\}.
\]
We call such an integers $j$ a balanced place (for $F_\mu$ and $F_\nu$). These balanced places $j$  are in bijection with the critical places
          $\frac{1}{2}+j$ of $\Pi\times \Sigma$ (see Section \ref{sec5.2}). As before, $F_\nu$ has highest weight $\nu = \{\nu^\iota\}_{\iota\in \CE_\rk}  \in (\BZ^{n-1})^{\CE_\rk}$, and $\nu^\iota= (\nu^\iota_1,\ldots, \nu^\iota_{n-1})$.
          
          Let $\frac{1}{2}+j$ be a critical place of $\Pi\times \Sigma$. Put
         \[
         \Omega_{\mu,\nu, j}:= \mathrm i^{\,j \frac{n(n-1)}{2}   [\rk\, :\, \BQ]+\sum_{\iota\in \CE_\rk}\sum_{i=1}^{n-1} (n-i)(\mu^\iota_i+\nu^\iota_i)}\qquad (\mathrm i:=\sqrt{-1}).
         \]
Let 
  $\chi: \rk^\times\backslash \BA^\times \rightarrow \BC^\times$ be a finite order Hecke character. We are concerned with the rationality of the critical value $\oL(\frac{1}{2}+j,\Pi\times \Sigma\times \chi)$, when both the critical place $\frac{1}{2}+j$ and the finite order Hecke character $\chi$ vary. Here $\oL(s,\Pi\times \Sigma\times \chi)$ denotes the completed Rankin-Selberg L-function.  Define the composition field
\[
\BQ(\Pi,\Sigma,\chi):=\BQ(\Pi)\BQ(\Sigma)\BQ(\chi)\subset \BC.
\]
Similar to $\Pi_\infty$, we have the archimedean parts $\Sigma_\infty$ and $\chi_\infty$ of $\Sigma$ and $\chi$ respectively. 


The main result of this article is the following global period relation.

    \begin{thmp}\label{thm: global period relation}
                   Let the notations and assumptions be as above. Then
 \begin{equation}\label{eq:mainthm0}   \frac{\oL(\frac{1}{2}+j,\Pi\times \Sigma\times \chi)}{ \Omega_{\mu,\nu,j} \cdot \CG(\chi_\Sigma)\cdot \mathcal G(\chi)^{\frac{n(n-1)}{2}} \cdot \Omega_{\varepsilon_n} (\Pi)\cdot \Omega_{\varepsilon_{n-1}}(\Sigma) } \in \BQ(\Pi, \Sigma, \chi),\end{equation}
where $\chi_\Sigma$ is the central character of $\Sigma$,  ``$\CG$" indicates   the Gauss sum (see  \eqref{gauss-sum000} and \eqref{gauss-sum002}), 
 and $\varepsilon_{n}, \varepsilon_{n-1}$ are the quadratic characters of $\rk_\infty^\times$ given by
 \[
 (\varepsilon_{n}, \varepsilon_{n-1}):=
 \left\{
   \begin{array}{ll}
   \left( \varepsilon_{\Sigma_\infty}\cdot \sgn_\infty^{\frac{(n-2)(n-3)}{2}+j} \cdot \chi_\infty, \, \varepsilon_{\Sigma_\infty}\cdot \sgn_\infty^{\frac{(n-2)(n-3)}{2}} \right),&\quad \textrm{if $n$ is even;}\smallskip \\
   \left( \varepsilon_{\Pi_\infty}\cdot \sgn_\infty^{\frac{(n-1)(n-2)}{2}}, \, \varepsilon_{\Pi_\infty}\cdot \sgn_\infty^{\frac{(n-1)(n-2)}{2}+j} \cdot \chi_\infty \right),&\quad \textrm{if $n$ is odd.}
   \end{array}
 \right.
 \]
Moreover, the quotient \eqref{eq:mainthm0} is $\Aut(\BC)$-equivariant in the sense that
 \be \label{eq:galois}
 \begin{aligned}
  & \sigma\left( \frac{\oL(\frac{1}{2}+j,\Pi\times \Sigma\times \chi)}{\Omega_{\mu, \nu, j}\cdot \CG(\chi_\Sigma)\cdot  \mathcal G(\chi)^{\frac{n(n-1)}{2}} \cdot \Omega_{\varepsilon_{n}}(\Pi)\cdot \Omega_{\varepsilon_{n-1}}(\Sigma) }  \right)\\
  = \ &  \frac{\oL(\frac{1}{2}+j, {}^\sigma\Pi \times {}^\sigma \Sigma\times {}^\sigma\chi)}{  \Omega_{\mu,\nu,j}\cdot  \CG(\chi_{{}^\sigma\Sigma})\cdot \mathcal G({}^\sigma\chi)^{\frac{n(n-1)}{2}} \cdot \Omega_{\varepsilon_{n}}({}^\sigma\Pi)\cdot \Omega_{\varepsilon_{n-1}}({}^\sigma\Sigma)} 
  \end{aligned}
 \ee
 for every $\sigma\in  \Aut(\BC)$.
    \end{thmp}


The proof of Theorem \ref{thm: global period relation} crucially depends on three local results that are responsible for the occurrence of the denominator in \eqref{eq:mainthm0}. More precisely:
\begin{itemize}
    \item The definition of the canonical Whittaker periods $\Omega_{\varepsilon_n} (\Pi)$ and $ \Omega_{\varepsilon_{n-1}}(\Sigma)$  relies on the non-vanishing hypothesis that was proposed by Kazhdan-Mazur in 1970's and proved by Sun in 2017 \cite{Sun}. 
    \item The appearance of the term $\CG(\chi_\Sigma)\cdot  \mathcal G(\chi)^{\frac{n(n-1)}{2}}$ is a consequence of the non-archimedean period relation (Proposition \ref{nonarchrel}), which is essentially due to Harder \cite[Section III]{H} for $n=2$ and Mahnkopf \cite[Section 3.4]{Mah} and Raghuram \cite[Section 3.3]{Rag1} in general. 
    \item The explicit  calculation of $\Omega_{\mu, \nu, j}$ is a consequence of the archimedean period relation (Theorem \ref{thmap}). The key contribution of this article is a proof of the  archimedean period relation, based on the preparatory work in \cite{LLSS}. The  proof  is much more involved than that of the non-archimedean period relation. 
\end{itemize}


In what follows we comment on some previous works  concerning Theorem \ref{thm: global period relation}. The first result was obtained by Shimura in 1959 (\cite[Section 9]{Sh1}). 
He proved that for certain nonzero complex numbers $\{\Omega_\epsilon\}_{\epsilon\in \{\pm 1\}}$,
\be\label{shram}
  \frac{\oL(k,\Delta)}{(2\pi\mathrm i)^k\cdot \Omega_{(-1)^k}}\in \BQ\quad \textrm{for  all }k=1,2, \cdots, 11. 
\ee
Here $\Delta$ is Ramanujan's cusp form of weight 12 and level 1 given by   
\[
\Delta(z) = q\prod_{n=1}^\infty 
(1-q^n
)^{24}=\sum_{n=1}^\infty \tau(n) q^n\qquad (q:=e^{2\pi\mathrm i\cdot z}),
\]
and the (incompleted) L-function $\oL(s,\Delta)$ is given by 
\[
 \oL(s,\Delta)=\sum_{n=1}^\infty \frac{\tau(n)} {n^{s}}\qquad (\textrm{when the real part of $s$ is sufficiently large}).
\]
When $n=2$, $\rk=\BQ$, $\chi$ and $\Sigma$ are trivial, and $\Pi$ is the  automorphic representation associated with $\Delta$, Theorem \ref{thm: global period relation} is a reformulation of  the relation \eqref{shram}.

After the aforementioned  pioneering work of Shimura,   a series of results towards Theorem \ref{thm: global period relation} for $n=2$ were obtained by Manin  \cite{Man1, Man2, Man3}, Shimura  \cite{Sh2, Sh3, Sh4}, and  Harder \cite{H}. Theorem \ref{thm: global period relation} for $n=2$ was finally proved in full generality by  Hida in 1994 (\cite[Theorem I]{Hi}).



For general $n$, the representation-theoretic problems behind Theorem \ref{thm: global period relation} are much more difficult than the case of $n=2$. The non-archimedean period relation is responsible for the rationality of 
$\oL(\frac{1}{2}+j,\Pi\times \Sigma\times \chi)$ when the finite order Hecke character $\chi$ varies, and the archimedean period relation is responsible for the rationality of 
$\oL(\frac{1}{2}+j,\Pi\times \Sigma\times \chi)$ when the critical place  $\frac{1}{2}+j$ varies. The  non-archimedean period relation is much easier to prove than the  archimedean period relation.
 Partly because of this reason,  more complete results on the rationality of  $\oL(\frac{1}{2}+j,\Pi\times \Sigma\times \chi)$ have been obtained for fixed $j$ and varying $\chi$, in a series of works including \cite{Sch2, KMS, Mah, KS, Rag1, Rag2, GH, G}. See also the survey paper \cite{HL} for more relevant works.
 
 However, it is also crucial to understand  the rationality of  $\oL(\frac{1}{2}+j,\Pi\times \Sigma\times \chi)$ when $\chi$ is fixed and $j$ varies, as in Shimura's result \eqref{shram}. For example, as explained in the Introduction of \cite{HN1}, this is essentially important for the Kummer congruence (also called Manin congruence) in the construction  of $p$-adic Rankin-Selberg L-functions (see \cite{Jan4}). Only some partial or conditional results (for varying $j$)  have been obtained in this direction (see \cite{Jan5, HR,  GL, HN1, Rag3}).

We have some more specific comments that compare Theorem \ref{thm: global period relation} with the existing results in the literature: 
\begin{itemize}
    \item The number field $\rk$ is assumed to be $\BQ$ in \cite{KMS, KS, Mah, Rag1}. It is assumed to be imaginary quadratic or CM in \cite{GH, G, GL}, with extra assumptions on $\Pi$ and $\Sigma$. In \cite[Theorem A]{GL} the rationality for varying $j$ is obtained under the hypotheses that certain central $\oL$-values are non-vanishing, which themselves remain a variety of difficult open problems.
    \item Under a hypothesis that is more or less equivalent to the archimedean period relation, a less precise version of Theorem 
    \ref{thm: global period relation} is proved in \cite{Jan5} for general $\rk$. For $n=3$, $\rk=\BQ$ and $\chi=1$, based on the explicit calculation of certain Rankin-Selberg zeta integrals in \cite{HIM}, Theorem \ref{thm: global period relation} is proved in \cite{HN1}.
    \item Roughly speaking, Theorem \ref{thm: global period relation} asserts that the transcendency of the critical L-values is  captured by the Whittaker periods. Harder and Raghuram prove in \cite[Theorem 7.21]{HR} that the transcendency of the ratio of  two successive critical L-values is 
    captured by the ``relative period" (which is in fact  the ratio of two Whittaker periods). They use Langlands-Shahidi method, and their result is proved for more general Rankin-Selberg $\oL$-functions and for $\rk$ totally real. This is extended to the case that $\rk$ is totally imaginary in \cite{Rag3}.  The results in the case of $\GL(n)\times \GL(n-1)$ are  immediate consequences of Theorem \ref{thm: global period relation}.
     \item For standard automorphic L-functions of symplectic type, a rationality  result similar to \eqref{eq:mainthm0}  is proved by Jiang, Sun and Tian in \cite{JST}. The reciprocity law, namely the equation \eqref{eq:galois}, is not proved in \cite{JST} for those L-functions. 
\end{itemize}

 In this article, we complete the story by giving an unconditional proof of Theorem \ref{thm: global period relation}, which is over arbitrary number fields. As we mentioned earlier, the key ingredient is  the archimedean period relation whose proof is very much involved. 

Last but not least, it is clear that the period relations (Theorem \ref{thm: global period relation}) have further applications towards the arithmetic study of other L-functions and Deligne's conjecture (see \cite{Mah, RS1, Rag1, Rag2, C1, C2, C3, HN2}), and they are also indispensable for the study of $p$-adic  L-functions (see \cite{Man2, Man3, Sch1, Sch2, Sch3, KMS} and \cite{Jan1}--\cite{Jan4}).

The article is organized as follows. In Section \ref{sec21} we translate 
 general cohomological representations to the cohomological representation with trivial coefficient system. This is the main idea used in the proof of the archimedean period relations (Theorem \ref{thmap}), which is formulated in Section   \ref{secap} and proved in Section \ref{sec3}. To this end, we recall the main result of \cite{LLSS} and use it to compare the Rankin-Selberg integrals with the integrals over a certain open orbit. 
In Section \ref{secnap} we reformulate the  non-archimedean period relations (Proposition \ref{nonarchrel}) and provide a proof of it for completeness.  
In Section \ref{sec4} we define the Whittaker periods of irreducible automorphic representations that are tamely isobaric and regular algebraic,  and study their properties under Galois twist.  We formulate the  global modular symbols and modular symbols at infinity, and explain their relationship in Section \ref{sec5}, which amounts to the unfolding of global Rankin-Selberg integrals as in \cite{JS2}. Finally, the global period relation Theorem \ref{thm: global period relation} is proved 
 in Section \ref{sec5} based on the results established in earlier sections. 

\section{Cohomological representations and their translations} \label{sec21}

In this section we introduce some generalities for cohomological representations, and give an explicit construction of the translation from the cohomological representations with trivial coefficient system to general ones.

\subsection{Cohomological representations} Let $\K$ be an archimedean local field. Thus it is a topological field that is topologically isomorphic to $\R$ or $\C$. Its complexification
\[
  \K\otimes_{\R}\C=\prod_{\iota\in \mathcal E_\K} \C,
\]
where $\mathcal E_\K$ denotes the set of all continuous field embeddings $\iota: \K\rightarrow \C$.

 Fix an integer $n\geq 1$, and fix a weight
\[
  \mu^\iota=(\mu_1^\iota\geq \mu_2^\iota\geq \cdots \geq \mu_n^\iota)\in \Z^n
\]
for every $\iota\in \mathcal E_\K$. Write $\mu:=\{\mu^\iota\}_{\iota\in \mathcal E_\K}$, and denote by $F_{\mu}$ the irreducible algebraic representation of $\GL_n(\K\otimes_\R \C)=\prod_{\iota\in \mathcal E_\K} \GL_n(\C)$ of highest weight $\mu$.

We say that $\mu$ is pure if 
\[
  \mu_{1}^\iota+\mu_{n}^{\bar \iota}=\mu_{2}^{\iota}+\mu_{n-1}^{\bar \iota}=\cdots=\mu_{n}^{\iota}+\mu_{1}^{\bar \iota},
\]
for every $\iota\in \mathcal E_\K$, where $\bar \iota$ denotes the composition of $\iota$ with the complex conjugation.
We suppose that $\mu$ is pure.  Denote by $\Omega(\mu)$ the set of isomorphism classes of irreducible Casselman-Wallach representations $\pi_\mu$ 
of $\GL_{n}(\K)$ such that
\begin{itemize}
                           \item
                            $\pi_\mu$ is generic, and  essentially unitarizable in the sense that $\pi_\mu\otimes \chi'$ is unitarizable for some character $\chi'$ of $\GL_n(\K)$; and
                           \item the total continuous cohomology
                           \[
                             \oH_{\mathrm{ct}}^*(\GL_n(\K)^0;
                           F_\mu^\vee\otimes \pi_\mu)\neq \{0\}.
                         \]
                             \end{itemize}
We remark that no such $\pi_\mu$ exists when $\mu$ is not pure. 

By \cite[Section 3]{Clo},
\begin{equation}\label{pimu}
  \#(\Omega(\mu))=\left\{
                \begin{array}{ll}
                  2, & \hbox{if $\K\cong \BR$ and and $n$ is odd;}\\
                  1, & \hbox{otherwise.} 
                \end{array}
              \right.
\end{equation}
Write $\sgn_{\K^\times}:\K^\times \rightarrow \C^\times$ for the quadratic character that is nontrivial if and only if $\K\cong \R$, and define the sign character
\[
\sgn:=\sgn_{\K^\times}\circ\det 
\]
of a general linear group $\GL_n(\K)$. 
Then in the first case of \eqref{pimu} the two members of $\Omega(\mu)$ are twists of each other by the sign character, and in the second case of \eqref{pimu} the only representation in $\Omega(\mu)$ is isomorphic to its own twist by the sign character. 
Recall that by \cite[Lemma 3.14]{Clo},
\[
  \oH_{\mathrm{ct}}^i(\GL_n(\K)^0;
                           F_\mu^\vee\otimes \pi_\mu)=\{0\},\qquad \textrm{if $i<b_{n,\K}$,}
\]
and
\[
   \oH_{\mathrm{ct}}^{b_{n,\K}}(\GL_n(\K)^0;      F_\mu^\vee\otimes \pi_\mu) \cong
  \left\{
                \begin{array}{ll}
                  1_{\K^\times}\oplus \sgn_{\K^\times}, & \hbox{if $\K\cong \BR$ and and $n$ is even;}\\
                  \varepsilon_{\pi_\mu}, & \hbox{otherwise,} 
                \end{array}
              \right.
\]
as representations of $\pi_0(\K^\times)$, where $1_{\K^\times}$ denotes the trivial character of $\K^\times$, and $\varepsilon_{\pi_\mu}$ denotes the central character 
of $F_\mu^\vee\otimes \pi_\mu$. Here and henceforth we make the identification 
 \[
 \pi_0(\GL_n(\K)) = \pi_0(\K^\times) \quad (\pi_0 \textrm{ indicates the set of connected components})
 \]
  through the determinant map $\GL_n(\K)\to \K^\times$. Note that $\varepsilon_{\pi_\mu}$ is equal to either  $1_{\K^\times}$ or   $\sgn_{\K^\times}$ for $\K\cong \BR$ and $n$ odd, and is trivial otherwise.

For every commutative ring $R$, let $\mathrm B_n(R)$  be the subgroup of $\GL_n(R)$ consisting of all the upper triangular matrices, and let $\mathrm N_n(R)$ be the subgroup of matrices in $\mathrm B_n(R)$  whose diagonal entries are  $1$. Likewise  let $\bar{\mathrm B}_n(R)$  be the subgroup of $\GL_n(R)$ consisting of all the lower triangular matrices, and let $\bar{\mathrm N}_n(R)$ be  the subgroup of matrices in $\bar{\mathrm B}_n(R)$  whose diagonal entries are  $1$.
Let $\mathrm{T}_n(R)$ be the subgroup of diagonal matrices in $\GL_n(R)$. 

Note that the invariant spaces $(F_\mu)^{\mathrm N_n(\K\otimes_\R \C)}$ and $(F_\mu^\vee)^{\bar{\mathrm N}_n(\K\otimes_\R \C)}$ are one-dimensional. We shall fix a generator $v_\mu\in (F_\mu)^{\mathrm N_n(\K\otimes_\R \C)}$ and a generator $v_\mu^\vee \in (F_\mu^\vee)^{\bar{\mathrm N}_n(\K\otimes_\R \C)}$ such that their pairing 
\be\label{pairvv}
\la v_\mu, v_\mu^\ve\ra =1.
\ee
To be more concrete, we
 realize  $F_\mu$  as the algebraic induction 
\be \label{algind1}
F_\mu= {}^{\rm alg} \Ind^{\GL_n(\K\otimes_\R \C)}_{\bar{\mathrm B}_n(\K\otimes_\R \C)}\chi_\mu,
\ee
and realize $v_\mu$ as the $ {\mathrm N}_n(\K\otimes_\R \C)$-invariant  
algebraic function $f$ 
in $F_\mu$  such that 
\[
f(1_n)=1\qquad (1_n\textrm{ denotes the identity element of $\GL_n(\K\otimes_\R \C)$}).
\]Here $\chi_\mu = \otimes_{\iota \in {\CE_\K}}\chi_{\mu^{\iota}}$ denotes the algebraic character of  $\mathrm T_n(\K\otimes_\R \C)$ corresponding to the weight $\mu\in (\Z^n)^{\CE_\K}$, to be viewed as an algebraic  character of $\bar{\mathrm B}_n(\K\otimes_\R \C)$ as usual. Similarly, we
 realize  $F_\mu^\vee$  as the algebraic induction 
\be \label{algind1}
F_\mu^\vee= {}^{\rm alg} \Ind^{\GL_n(\K\otimes_\R \C)}_{{\mathrm B}_n(\K\otimes_\R \C)}\chi_{-\mu},
\ee
and realize $v_\mu^\vee$ as the ${\bar{\mathrm N}}_n(\K\otimes_\R \C)$-invariant  
algebraic function $f^\vee$ 
in $F_\mu^\vee$  such that $f^\vee(1_n)=1$. The invariant pairing $\la\,,\,\ra: F_\mu\times F_\mu^\vee\rightarrow \C$ is determined by the equality \eqref{pairvv}. Note that as a linear functional on $F_\mu^\vee$, $v_\mu$ equals the evaluation map at $1_n$. Similarly, $v_\mu^\vee$ equals the evaluation map at $1_n$ as a linear functional on $F_\mu$.

Fix a unitary character 
\be\label{psir}
\psi_\R: \R\rightarrow \C^\times,\quad x\mapsto e^{2\pi \mathrm{i}x},
\ee
which induces 
 a unitary character
\be \label{psiK}
  \psi_\K: \K\rightarrow \C^\times, \quad x\mapsto \psi_\R\left(\sum_{\iota\in \mathcal E_\K} \iota(x)\right).
\ee
This further induces a unitary character
\be \label{psi_n}
  \psi_{n,\K}: \mathrm{N}_n(\K)\rightarrow \C^\times, \quad [x_{i,j}]_{1\leq i,j\leq n} \mapsto \psi_\K\left((-1)^n\cdot \sum_{i=1}^{n-1} x_{i,i+1}\right).
\ee
By abuse of notation, we will still use $\psi_{n,\K}$ to denote the space $\BC$ carrying the representation of $\mathrm{N}_n(\K)$ corresponding to the character $\psi_{n,\K}$. Similar notation will be freely used for other characters. 
Let $\pi_\mu\in \Omega(\mu)$. Recall that the space $\Hom_{ \mathrm{N}_n(\K)}(\pi_\mu, \psi_{n,\K})$ is one-dimensional. Fix a generator \be\label{genw}
\lambda_\mu\in \Hom_{ \mathrm{N}_n(\K)}(\pi_\mu, \psi_{n,\K}),
\ee
to be called the Whittaker functional on $\pi_\mu$. 

Write $0_{n,\K}$ for the zero element of $(\Z^n)^{\mathcal E_\K}$. Then $F_{0_{n,\K}}$ is the trivial representation. Specifying the above argument to the case when $\mu=0_{n,\K}$, we take a representation $\pi_{0_{n,\K}}\in \Omega(0_{n,\K})$, together with the Whittaker  functional $\lambda_{0_{n,\K}}\in \Hom_{ \mathrm{N}_n(\K)}(\pi_{0_{n,\K}}, \psi_{n,\K})\setminus \{0\}$.

Throughout this article, we assume that the representation $\pi_{0_{n,\K}} \in \Omega(0_{n,\K})$ is chosen such that $\pi_{0_{n,\K}}$ and $F_\mu^\vee \otimes \pi_\mu$ have the same central character, to be denoted by $\varepsilon_{n,\K}$.

 \subsection{Explicit translations} \label{sec2.2}
 We will prove the following result in this subsection.
 
\begin{prp} \label{prpjmath}
 There is a unique element $\jmath_\mu\in \Hom_{\GL_n(\K)}(\pi_{0_{n,\K}},  F_\mu^\vee\otimes \pi_\mu)$ such that the diagram
\[
 \begin{CD}
          \pi_{0_{n,\K}}
                  @>\jmath_\mu >> F_\mu^\vee\otimes \pi_\mu \\
            @V \lambda_{0_{n,\K}}VV            @VV v_\mu\otimes \lambda_{\mu} V\\
         \C@=\C \\
  \end{CD}
\]
commutes. Moreover, $\jmath_\mu$ induces a linear isomorphism
\[
 \jmath_\mu :  \oH_{\mathrm{ct}}^i(\GL_n(\K)^0; \pi_{0_{n,\K}})\xrightarrow{\sim}  \oH_{\mathrm{ct}}^i(\GL_n(\K)^0;  F_\mu^\vee\otimes \pi_\mu)
\]
of representations of $\pi_0(\K^\times)$ for each  $i\in \Z$.
\end{prp}


It is known from the Vogan-Zuckerman theory of cohomological representations (see Proposition 1.2 and Section 5 of \cite{VZ}) that 
\be \label{mult1}
\dim \Hom_{\GL_n(\K)}(\pi_{0_{n,\K}},  F_\mu^\vee\otimes \pi_\mu)=1.
\ee 
We first recall the realization of $\pi_\mu$ and introduce a certain principal series representation $I_\mu$ of $\GL_n(\K)$.  Define a character 
\[
\rho_n:=\otimes^n_{i=1}\abs{\, \cdot \,}_\K^{\frac{n+1}{2}-i}
\qquad (\abs{\,\cdot\,}_\K \textrm{ denotes the normalized absolute value})
\]
of $\mathrm T_n(\K)$. 
For $\iota\in \CE_\K$,  define  the half-integers
\be \label{lmu}
\tilde{\mu}_i^\iota : = \mu_i^\iota + \frac{n+1-2i}{2}, \quad i =1, \ldots, n.
\ee

For $\K \cong \R$, $a, b\in\BC$ with $a-b\in \BZ\setminus \{0\}$, denote by $D_{a, b}$ the relative discrete series representation of $\GL_2(\K)$ with infinitesimal character $(a,b)$. 
 If $n$ is even, then 
\[
\pi_\mu\cong \Ind^{\GL_n(\K)}_{\bar{{\mathrm P}}_n(\K)}\left(D_{\tilde{\mu}^\iota_1, \tilde{\mu}^\iota_n}\otimes \cdots \otimes D_{\tilde{\mu}^\iota_{\frac{n}{2}}, \tilde{\mu}^\iota_{\frac{n}{2}+1}}\right)
\quad (\textrm{normalized smooth induction}),
\]
where $\bar{{\mathrm P}}_n$ is the lower triangular parabolic subgroup of type $(2,\ldots, 2)$. If $n$ is odd, then 
\[
\pi_\mu\cong \Ind^{\GL_n(\K)}_{\bar{{\mathrm P}}_n(\K)}\left(D_{\tilde{\mu}^\iota_1, \tilde{\mu}^\iota_n}\otimes \cdots \otimes D_{\tilde{\mu}^\iota_{\frac{n-1}{2}}, \tilde{\mu}^\iota_{\frac{n+3}{2}}}\otimes  (\cdot)^{\tilde{\mu}^\iota_{\frac{n+1}{2}}}\varepsilon_{n,\K} \right),
\]
where $\bar{{\mathrm P}}_n$ is the lower triangular parabolic subgroup of type $(2,\ldots, 2, 1)$, and we recall that $\varepsilon_{n,\K}=1_{\K^\times}$ or $\sgn_{\K^\times}$ is the common central character of 
$F_\mu^\vee\otimes \pi_\mu$ and $\pi_{0_{n,\K}}$.

For $\K\cong \C$,  we have that
\[
\pi_\mu\cong \Ind^{\GL_n(\K)}_{\bar{\mathrm B}_n(\K)} \left(\iota^{\tilde{\mu}^\iota_1} \bar{\iota}^{\tilde{\mu}^{\bar{\iota}}_n}\otimes\cdots\otimes \iota^{\tilde{\mu}^\iota_n} \bar{\iota}^{\tilde{\mu}^{\bar{\iota}}_1}\right),
\]
where for $a, b\in \BC$ with $a-b\in\BZ$, $\iota^a \bar{\iota}^b$ denotes the character
\[
\iota^a \bar{\iota}^b: \K^\times \to \BC^\times, \quad z\mapsto \iota(z)^{a-b} (\iota(z)\bar{\iota}(z))^b.
\]

In both the real and complex cases, we define the principal series representation 
\[
I_\mu := \Ind^{\GL_n(\K)}_{\bar{\mathrm B}_n(\K)} (\chi_\mu \cdot \rho_n\cdot ( \varepsilon_{n,\K}\circ \det))= {}^{\rm u}\Ind^{\GL_n(\K)}_{\bar{\mathrm B}_n(\K)} (\chi_\mu \cdot  ( \varepsilon_{n,\K}\circ \det)),
\]
so that $I_\mu$ and $\pi_\mu$ have the same central character. Here ${}^{\rm u}\Ind$ stands for the unnormalized smooth parabolic induction. 
 
\begin{lemp} \label{lem2.1}
The principal series representation  $I_\mu$ has a unique irreducible generic quotient representation, which is isomorphic to $\pi_\mu$. 
\end{lemp}

\begin{proof}
By \cite[Lemma 2.5]{J2}, $I_\mu^\vee$ has a unique irreducible generic subrepresentation, hence $I_\mu$ has a unique irreducible generic quotient representation. 

For $\K\cong\BR$, by the well-known realization of relative discrete series representations of $\GL_2(\BR)$ as quotients of principal series representations,  $\pi_\mu$ is a quotient of 
\[
\Ind^{\GL_n(\K)}_{\bar{\mathrm B}_n(\K)} (w(\chi_\mu\cdot \rho_n)\cdot( \varepsilon_{n,\K}\circ \det))
\]
for a certain $w\in W_n$. Here $W_n$ is the subgroup of permutation matrices in 
$\textrm{GL}_n(\BZ)$ which acts on $\mathrm T_n(\K)$ by conjugation, and thus acts on the set of characters of $\mathrm T_n(\K)$. 
The above  representation and $I_\mu$ have the same irreducible  constituents. Hence $\pi_\mu$ is isomorphic to the unique irreducible generic subquotient of $I_\mu$, which is in fact a quotient as we mentioned at the beginning of the proof. 

For $\K\cong\BC$, the lemma follows easily from \cite[Theorem 6.2]{JL} (for the case of $\GL_2(\BC)$) and parabolic induction in stages. 
\end{proof}

The group $\RN_n(\K)$ is equipped with the Haar measure
\be\label{mea1}
\od\! u:=\prod_{1\leq i< j\leq n} \od\!u_{i,j}, \qquad u=[u_{i,j}]_{1\leq i,j\leq n}\in \RN_n(\K),
 \ee
where $\K$ is equipped with the self-dual Haar measure with respect to $\psi_\K$.  By \cite[Theorem 15.4.1]{Wa2}, 
\be \label{I-whit-mult1}
\dim \Hom_{\mathrm{N}_n(\K)}(I_\mu, \psi_{n,\K}) =1
\ee
and there is a unique $\lambda_\mu'\in \Hom_{\mathrm{N}_n(\K)}(I_\mu, \psi_{n,\K})$ such that 
\be \label{la'} 
\lambda_\mu'(f)=\int_{\mathrm{N}_n(\K)} f(u)\overline{\psi_{n,\K}}(u)\od\! u
\ee
for all $f\in I_\mu$ such that $f|_{\mathrm{N}_n(\K)}\in \CS(\mathrm{N}_n(\K))$. Here and henceforth, for a Nash manifold $X$, denote by $\CS(X)$ the space of Schwartz functions on $X$
(see \cite{Du, AG1}).

As usual, an element $u\otimes f\in F_\mu^\vee\otimes I_\mu$ is identified with the function 
\[
\GL_n(\K)\rightarrow F_\mu^\vee, \quad g\mapsto f(g) \cdot u.
\]
Then  $F_\mu^\vee\otimes I_\mu$ is identified with the space of  $F_\mu^\vee$-valued  smooth functions $\varphi$ on $\GL_n(\K)$ satisfying  that 
\[
\varphi(bx)=\left(((\varepsilon_{n,\K}\circ \det)\cdot \chi_\mu)(b)\right)\cdot \varphi(x),\quad \textrm{for all }b\in \bar{\mathrm B}_n(\K),\  x\in \GL_n(\K),
\]
on which $\GL_n(\K)$ acts by
\[
(g . \varphi)(x):=g. (\varphi(xg)),\quad \textrm{where } \, g, \,  x \in \GL_n(\K).
\]
Define a $\GL_n(\K)$-homomorphism
\be \label{imu}
\begin{array}{rcl}
\imath_\mu:I_{0_{n,\K}}&\rightarrow &  F_\mu^\vee \otimes I_\mu,\\
f&\mapsto & \left(g\mapsto  f(g)\cdot (g^{-1}. v_\mu^\vee)\right).
\end{array}
\ee

\begin{lemp} \label{lem2.2}
The map $\imath_\mu$ satisfies that 
\be \label{eq2.2}
(v_\mu \otimes \lambda_\mu' )  \circ \imath_\mu =  \lambda_{0_{n,\K}}'.
\ee
\end{lemp}

\begin{proof}
Recall that $v_\mu\in F_\mu$ is ${\mathrm N}_n(\K\otimes_\BR \C)$-invariant, and $\langle v_\mu, v_\mu^\vee\rangle =1$.
For $f\in I_{0_n,\K}$ with $f|_{\mathrm{N}_n(\K)}\in \CS(\mathrm{N}_n(\K))$, we have that
\[
\begin{aligned}
& \left((v_\mu \otimes \lambda_\mu' )  \circ \imath_\mu\right)(f) \\
=  \ & \int_{\mathrm{N}_n(\K)} \langle v_\mu, \imath_\mu(f)(u)\rangle \overline{\psi_{n,\K}}(u)\od\! u \\
= \ & \int_{\mathrm{N}_n(\K)} f(u)\overline{\psi_{n,\K}}(u)\od\! u \\
= \ & \lambda_{0_{n,\K}}'(f).
\end{aligned}
\]
This proves \eqref{eq2.2}, in view of \cite[Theorem 15.4.1]{Wa2}.
\end{proof}

By Lemma \ref{lem2.1} and \eqref{I-whit-mult1}, there is a unique $p_\mu\in \Hom_{\GL_n(\K)}(I_\mu, \pi_\mu)$ such that 
\be \label{eq-p}
\lambda_\mu\circ p_\mu=\lambda_\mu'.
\ee
Let $J_\mu:=\Ker(p_\mu)$, which is the largest  subrepresentation of $I_\mu$ such that
\[
\textrm{Dim} \, J_\mu < \textrm{Dim}\,  I_\mu.
\]
Here and below, $\textrm{Dim}$ indicates the Gelfand-Kirillov dimension of a Casselman-Wallach representation  of $\GL_n(\K)$. Likewise we have $J_{0_{n,\K}}:=\Ker(p_{0_{n,\K}})\subset I_{0_{n,\K}}$. 

\begin{lemp} \label{lem2.3}
It holds that
\[
\imath_\mu(J_{0_{n,\K}})\subset F_\mu^\vee \otimes J_\mu.
\]
\end{lemp}

\begin{proof} It suffices to show that 
\[
\tilde{\imath}_\mu(F_\mu \otimes J_{0_{n, \K}} )\subset J_\mu,
\]
where $\tilde{\imath}_\mu\in \Hom_{\GL_n(\K)}(F_\mu\otimes I_{0_{n, \K}}, I_\mu)$ is  the linear map induced by $\iota_\mu$.
This follows from the fact that (see \cite[Lemma 2.2]{V})
\[
\textrm{Dim} \, ( F_\mu \otimes J_{0_{n, \K}} ) = \textrm{Dim} \, J_{0_{n, \K}}.
\]
\end{proof}

By Lemma \ref{lem2.3}, there is a unique $\jmath_\mu\in \Hom_{\GL_n(\K)}(\pi_{0_{n,\K}},  F_\mu^\vee\otimes \pi_\mu)$ such that the  diagram 
\be \label{I-pi}
 \begin{CD}
           F_\mu^\vee\otimes I_\mu 
                  @> {\rm id}\otimes p_\mu  >>   F_\mu^\vee\otimes \pi_\mu \\
            @A  \imath_\mu A  A         @AA  \jmath_\mu  A\\
           I_{0_{n, \K}}  @> p_{0_{n, \K}} >> \pi_{0_{n,\K}}\\
  \end{CD}\ \ \qquad(\textrm{$\rm id$ indicates the identity map})
\ee
commutes.
By \eqref{eq2.2} and \eqref{eq-p}, 
\[
\begin{aligned}
& (v_\mu\otimes \lambda_\mu)\circ \jmath_\mu\circ p_{0_{n,\K}} \\
= \ & (v_\mu\otimes \lambda_\mu)\circ ( {\rm id}_{F_\mu^\vee}\otimes p_\mu) \circ \imath_\mu \\
= \ & (v_\mu \otimes \lambda_\mu')\circ \imath_\mu \\
= \ &  \lambda'_{0_{n,\K}} \\
= \ & \lambda_{0_{n, \K}}\circ p_{0_{n, \K}},
\end{aligned}
\]
which implies that 
\[
(v_\mu\otimes \lambda_\mu)\circ \jmath_\mu = \lambda_{0_{n, \K}}.
\]
This proves the existence part of Proposition \ref{prpjmath}. The uniqueness follows from \eqref{mult1}. The last statement of the proposition follows from \cite[Section 5]{VZ}.

\section{Archimedean period relations}\label{secap}

In this section we explain the statement of the archimedean period relation (Theorem \ref{thmap}), 
whose proof will be given in the next section.

\subsection{Some cohomology spaces}
For simplicity, write
\[
 \oH_\mu:= \oH_{\mathrm{ct}}^{b_{n,\K}}(\GL_n(\K)^0;  F_\mu^\vee\otimes \pi_\mu),
\]
which is of dimension 1 or 2. As in Proposition \ref{prpjmath}, we have a linear isomorphism
\[
\jmath_\mu :\oH_{0_{n,\K}}   \xrightarrow{\sim} \oH_\mu
\]
of representations of $\pi_0(\K^\times)$.


Fix a maximal compact subgroup 
\be \label{cpt-inf}
K_{n,\K} :=\left\{ \begin{array}{ll} \RO(n), & {\rm if \ }\K\cong\BR; \\
\RU(n), & {\rm if \ }\K\cong\BC
\end{array}\right.
\ee
 of $\GL_n(\K)$.  The determinant homomorphism yields  identifications
 \[
   \pi_0(\GL_n(\K))=\pi_0(K_{n,\K})=\pi_0(\K^\times). 
 \]

 We use the corresponding  lower case gothic letter to denote the Lie algebra of a Lie group. For example the Lie algebra of $K_{n,\K}$ will be denoted by $\frak{k}_{n,\K}$. Put
\[
  d_{n,\K}:=b_{n+1,\K}+b_{n,\K}=\dim_\R (\g\l_n(\K)/\k_{n,\K})=\left\{
         \begin{array}{ll}
         \frac{n(n+1)}{2}, & \hbox{if $\K\cong \R$;} \smallskip \\
           n^2, & \hbox{if $\K\cong \C$.}
         \end{array}
       \right.
\]
Define a one-dimensional real vector space
\[
  \omega_{n,\K}(\R):=\wedge^{d_{n,\K}} (\g\l_n(\K)/\k_{n,\K}).
\]
Put
\[
  \omega_{n,\K}:=\omega_{n,\K}(\R)\otimes_\R \C,
\]
which is naturally  a representation of $\pi_0(\K^\times)$ that is isomorphic  to $\sgn_{\K^\times}^{n-1}$. Here and henceforth, we also view $1_{\K^\times}$ and $\sgn_{\K^\times}$ as representations of $\pi_0(\K^\times)$. Then there is an identification
\be \label{omega-coh}
   \oH_{\mathrm{ct}}^{d_{n,\K}}(\GL_n(\K)^0;  \omega_{n,\K})=1_{\K^\times}
\ee
of representations of $\pi_0(\K^\times)$.

Write $ \omega_{n,\K}^+$ and $\omega_{n,\K}^-$ for the two connected components of  $ \omega_{n,\K}(\R)\setminus \{0\}$, which are viewed as left invariant orientations on 
$\GL_n(\K)/K_{n,\K}^0$. The complex orientation space of  $ \omega_{n,\K}(\R)$ is defined to be the one-dimensional space
\be \label{orient}
  \mathfrak{O}_{n,\K}:=\frac{\C\cdot \omega_{n,\K}^+ \oplus \C\cdot  \omega_{n,\K}^- }{\{a( \omega_{n,\K}^+ + \omega_{n,\K}^-)\, : \,  a\in \C\}}.
\ee
Then $\pi_0(\K^\times)=\pi_0(K_{n,\K})$ acts on   $\mathfrak{O}_{n,\K}$ by $\sgn_{\K^\times}^{n-1}$, through the right translation on $\GL_n(\K)/K_{n,\K}^0$. We  identify
$
\omega_{n,\K}^*\otimes  \mathfrak{O}_{n,\K}
$
with the space of invariant measures on $\GL_n(\K)/ K_{n,\K}^0$ in the obvious way. Here and as usual, a superscript $*$ over a vector space indicates the dual space. Denote by $\frak{M}_{n,\K}$ the one-dimensional space of invariant measures on $\GL_n(\K)$. 
By push-forward of measures through the map $\GL_n(\K)\to \GL(\K)/ K_{n,\K}^0$, we have an identification 
\[
\frak{M}_{n,\K}  = \omega_{n,\K}^*\otimes  \mathfrak{O}_{n,\K}.
\]
In view of this and \eqref{omega-coh}, we have that
\[
   \oH_{\mathrm{ct}}^{d_{n,\K}}(\GL_n(\K)^0;    \mathfrak{M}_{n,\K}^*) \otimes \frak{O}_{n,\K}= \BC.
\]

\subsection{Archimedean modular symbols and archimedean period relations}\label{secams}
Now we assume that $n\geq 2$. Let $\nu\in (\Z^{n-1})^{\mathcal E_\K}$ be a highest weight and assume that it is pure. Then as before we have representations $F_\nu$  and $\pi_\nu$ of $\GL_{n-1}(\K\otimes_\BR\BC)$ and $\GL_{n-1}(\K)$ respectively, an element  $v_\nu\in F_\nu$, an element $v_\nu^\vee \in F_\nu^\vee$, and a Whittaker  functional  $\lambda_\nu$ on $\pi_\nu$. The representation $\pi_{0_{n-1,\K}}$ is determined by $\pi_\nu$ as before, and 
we have a linear isomorphism
 \[
\jmath_\nu :\oH_{0_{n-1,\K}}   \xrightarrow{\sim} \oH_\nu
\]
of  representations of $\pi_0(\K^\times) $.

As usual, we have an embedding 
\begin{equation}\label{emh}
 \imath: \GL_{n-1}(R)\hookrightarrow \GL_n(R),\quad g\mapsto \left[
                      \begin{array}{cc}
                        g & 0 \\
                        0 & 1 \\
                      \end{array}
                    \right],
\end{equation}
where $R$ is an arbitrary commutative ring.  We also view $\GL_{n-1}(R)$ as a subgroup of $\GL_n(R)\times \GL_{n-1}(R)$ via the diagonal embedding
\[
\GL_{n-1}(R) \hookrightarrow \GL_n(R)\times \GL_{n-1}(R), \quad g\mapsto \left(\begin{bmatrix} g & 0 \\ 0 & 1 \end{bmatrix}, g \right).
\]
For all $k,l\in \BN$,  denote by  $R^{k\times l}$ the set of $k\times l$ matrices with entries in $R$.

Put $\xi:=(\mu,\nu)$. Write $F_\xi:=F_\mu \otimes F_\nu$. Assume that $\xi$ is {\it balanced}  in the sense that there is an integer $j$ such that
 \begin{equation}\label{homf002}
  \Hom_{\GL_{n-1}(\K\otimes_\R \C)}(F_\xi^\vee, \otimes_{\iota\in \mathcal E_\K}{\det}^j) \neq 0.
\end{equation}

For each $k\in \BN$, write
\[
 w_k:=\left[
        \begin{array}{cccc}
          0 & \cdots & 0 & 1 \\
          0 & \cdots & 1 & 0 \\
          &\cdots & \cdots &  \\
          1 & 0 & \cdots & 0 \\
        \end{array}
      \right]\in  \GL_k(\BZ).
\]
Following \cite{LLSS}, define a family  $\{z_k\in \GL_k(\BZ)\}_{k\in \BN}$ of matrices
inductively by
\[
 z_0:=\varnothing\ \  (\textrm{the unique element of $\GL_0(\BZ)$}), \quad  z_1:=[1],
 \]
and
\be\label{zk}
       z_k :=
 \left[
            \begin{array}{cc}
                            w_{k-1}& 0 \\
                       0 & 1 \\
                     \end{array}
                   \right]
                    \left[
            \begin{array}{cc}
                       {}^t z_{k-2}^{-1}& 0 \\
                       0 & 1_2 \\
                     \end{array}
                   \right]
                    \left[
            \begin{array}{cc}
                       {}^tz_{k-1}  w_{k-1} z_{k-1}& {}^t e_{k-1}\\
                       0 & 1 \\
                     \end{array}
                   \right], \quad \textrm{for all  $k\geq 2$.}
                   \ee
Here and as usual, a left superscript $t$ over a matrix  indicates the transpose, 
$1_2$ stands for the $2\times 2$ identity matrix, and $e_{k-1}:=[0,\cdots, 0,1]\in \BZ^{1\times(k-1)}$. 

Let $j\in \Z$ be as in \eqref{homf002}. The following proposition follows from the fact that 
\[
\left(\bar{\mathrm B}_n(\C) \times \bar{\mathrm B}_{n-1}(\C)\right) \cdot (z_n, z_{n-1}) \cdot \GL_{n-1}(\BC)
\]
is Zariski open in $\GL_n(\C) \times \GL_{n-1}(\C)$ (see \cite[Lemma 1.1]{LLSS}).  

\begin{prp}\label{phixij}
There is a unique element
\[
  \phi_{\xi,j}\in\Hom_{\GL_{n-1}(\K\otimes_\R \C)}(F_\xi^\vee , \otimes_{\iota\in \mathcal E_\K}{\det}^j)
\]
such that
\[
 \phi_{\xi,j}((z_n^{-1} . v_\mu^\vee)\otimes (z_{n-1}^{-1}. v_\nu^\vee))=1.
\]
\end{prp}

Fix a quadratic character $\chi_\K$ of $\K^\times$. Define  characters 
\be\label{chikt}
\chi_{\K,t}:=\chi_\K\cdot \abs{\,\cdot\,}_\K^t,\quad\quad \chi_{\K}^{(j)}:=\chi_{\K}\cdot \sgn_{\K^\times}^j\qquad (t\in \C),
\ee
and more generally,
\[
 \chi_{\K,t}^{(j)}:=\chi_{\K}\cdot\abs{\,\cdot\,}_\K^t\cdot \sgn_{\K^\times}^j,
\]
of the group $\K^\times$.
When no confusion arises, for every commutative ring $R$, every character of $R^\times$ is identified with a character of $\GL_{n-1}(R)$ via the pullback through the determinant homomorphism. In particular, $\chi_{\K,t}^{(j)}$ is also viewed as a character of $\GL_{n-1}(\K)$. 

Put 
\[
 \pi_\xi:= \pi_\mu\widehat \otimes \pi_\nu \qquad (\textrm{the completed projective tensor product}).
\]
We have the normalized Rankin-Selberg integral (see \cite{J2})
\begin{eqnarray*}
 \oZ^\circ_\xi (\cdot, s, \chi_\K)&\in& \Hom_{\GL_{n-1}(\K)}\left(\pi_\xi \otimes \frak{M}_{n-1,\K}, \, \chi_{\K, -s+\frac{1}{2}}\right)\\
 &=&\Hom_{\GL_{n-1}(\K)}\left( \pi_\xi, \chi_{\K, -s+\frac{1}{2}}\otimes  \mathfrak{M}^*_{n-1,\K}\right),
\end{eqnarray*}
such that  
\be\label{nrs}
 \begin{aligned}
 &   \quad \oZ^\circ_\xi( f\otimes f' \otimes m, s, \chi_\K) := \frac{1}{\oL(s, \pi_\mu\times\pi_\nu)} \\
  &   \cdot \int_{\mathrm N_{n-1}(\K)\backslash \GL_{n-1}(\K)} \lambda_\mu\left(\begin{bmatrix} g & 0 \\ 0 & 1 \end{bmatrix}.f\right) \lambda_{\nu} (g. f') \cdot \chi_{\K}(\det g) \cdot \abs{\det g}_\K^{s-\frac{1}{2}} \od\!\bar{m} (g)
\end{aligned}
\ee
for all $f\in \pi_\mu$, $f' \in \pi_\nu$, $m \in \frak{M}_{n-1,\K}$, and $s\in \BC$ with the real part ${\rm Re}(s)$ sufficiently large (it extends to all $s\in \BC$ by holomorphic continuation). Here and henceforth,  $\bar{m}$ is the quotient measure 
on $\RN_{n-1}(\K) \backslash \GL_{n-1}(\K)$ induced by $m$. Recall that a Haar measure on $\mathrm{N}_{n-1}(\K)$ has been fixed as in \eqref{mea1}.

Put
\[
\RH_{\chi_\K^{(j)}} := \RH_{\rm ct}^0( \GL_{n-1}(\K)^0; \chi_\K^{(j)}).
\]
Let $\phi_{\xi,j}$ be as in Proposition \ref{phixij}. Then we have a $\GL_{n-1}(\K)$-equivariant continuous linear map
\be \label{phioZ}
\begin{aligned}
 \phi_{\xi,j}\otimes  \oZ_{\xi}^\circ(\cdot, \frac{1}{2}+j, \chi_\K): F_\xi^\vee\otimes \pi_\xi \rightarrow & \  (\otimes_{\iota \in \mathcal E_\K}{\det}^{j})\otimes  (\chi_{\K, -j}\otimes   \mathfrak{M}^*_{n-1,\K}) \\
= & \  \chi_\K^{(j)} \otimes   \mathfrak{M}^*_{n-1,\K}. 
 \end{aligned}
\ee
By restriction of cohomology, this induces a linear map
\begin{eqnarray*}
\wp_{\xi, \chi_\K, j}&:&  \oH_\mu\otimes \oH_\nu\otimes \oH_{\chi_\K^{(j)}} \otimes \mathfrak{O}_{n-1,\K} \\
 & = & \oH_{\mathrm{ct}}^{d_{n-1,\K}}(\GL_n(\K)^0\times \GL_{n-1}(\K)^0;  F_\xi^\vee\otimes \pi_\xi\otimes \chi_\K^{(j)}) \otimes \frak{O}_{n-1,\K}\\
&\rightarrow &\oH_{\mathrm{ct}}^{d_{n-1,\K}}(\GL_{n-1}(\K)^0;   \mathfrak{M}^*_{n-1,\K})\otimes \mathfrak{O}_{n-1,\K} = \BC.
\end{eqnarray*}
We call this map the {\it archimedean modular symbol}, which is nonzero by the non-vanishing hypothesis that is proved in \cite{Sun}.

Specifying the above argument to the case when $\xi=\xi_0:=(0_{n, \K}, 0_{n-1, \K})$ and $j=0$, we get a linear map (with $\chi_\K$ replaced by $\chi_\K^{(j)}$)
\be \label{wp00}
\wp_{\xi_0,\chi_\K^{(j)}, 0} : \oH_{0_{n,\K}}\otimes \oH_{0_{n-1,\K}}\otimes \oH_{\chi_\K^{(j)}}\otimes  \mathfrak{O}_{n-1,\K} \to \BC.
\ee


The archimedean period relation  is the following theorem.

\begin{thmp}\label{thmap}
Let the notations and assumptions be as above. Let
\be \label{OmegaKj}
\Omega_{\mu,\nu, j}':= 
 {\mathrm i}^{ j  \frac{n(n-1)}{2} [\K\, :\,  \BR]}\cdot   c'_\mu \cdot c_\nu \cdot  \varepsilon_{\mu,\nu},
\ee
where 
\[
\begin{aligned}
& c'_\mu := \prod^{n-1}_{i=1}  ((-1)^n   \mathrm{i})^{(n-i)\sum_{\iota\in \CE_\K}\mu^{\iota}_i}, \\
&  c_\nu: = \prod^{n-1}_{i=1}  ((-1)^n   \mathrm{i})^{(n-i)\sum_{\iota \in \CE_\K}\nu^{\iota}_i}, \quad and\\
 & \varepsilon_{\mu,\nu}:=\prod_{i>k, \, i+k\leq n}(-1)^{\sum_{\iota\in \CE_\K}(\mu_i^\iota+\nu_k^\iota)}.
\end{aligned}
\]
Then the diagram
\[
 \begin{CD}
           \oH_\mu\otimes \oH_\nu \otimes \oH_{\chi_\K^{(j)}} \otimes \mathfrak{O}_{n-1,\K}
                  @>\Omega'_{\mu,\nu, j }\cdot \wp_{\xi,\chi_\K, j} >>  \BC  \\ 
            @ A \jmath_\mu\otimes \jmath_\nu\otimes {\rm id} \otimes {\rm id} A A          @  | \\
           \oH_{0_{n,\K}}\otimes \oH_{0_{n-1,\K}} \otimes \oH_{\chi_\K^{(j)}}\otimes \mathfrak{O}_{n-1,\K} @>\wp_{\xi_0,\chi_\K^{(j)}, 0}>>  \BC \\
  \end{CD}
\]
commutes.
\end{thmp}

\begin{remarkp}
The Rankin-Selberg integrals for minimal $K$-type vectors of principal series representations of $\GL_n(\K)\times\GL_{n-1}(\K)$ have been explicitly calculated by Ishii and Miayzaki in \cite{IM}.
The Rankin-Selberg integrals for minimal $K$-type vectors of irreducible generalized principal series representations  of $\GL_3(\K)\times \GL_2(\K)$ have been explicitly calculated by  Hirano, Ishii and Miyazaki in \cite{HIM}.  It should be also possible to prove Theorem \ref{thmap} when $n\leq 3$ or  $\K\cong\BC$, by using these results and the method in \cite{Sun}.
\end{remarkp}



\section{Proof of archimedean period relations} \label{sec3}

In this section we prove the archimedean period relations Theorem \ref{thmap}. Retain the notation of the last section.

Put
\be \label{jmath}
\jmath_\xi:=\jmath_\mu \otimes \jmath_\nu \in \Hom_{\GL_n(\K)\times\GL_{n-1}(\K)}(\pi_{\xi_0},  F_\xi^\vee\otimes \pi_\xi).
\ee
We will prove the following result, which implies  Theorem \ref{thmap} by specifying $s$ to $\frac{1}{2}$.

\begin{thmp}\label{thmap2}
The diagram 
\[
 \begin{CD}
            F_\xi^\vee  \otimes \pi_\xi=(F_\mu^\vee \otimes \pi_\mu)\widehat \otimes (F_\nu^\vee\otimes \pi_\nu)
                  @>  \Omega'_{\mu,\nu, j}\cdot \phi_{\xi,j}\otimes \oZ_{\xi}^\circ(\cdot, s+j, \chi_\K) >> \chi_{\K, \frac{1}{2}-s}^{(j)}\otimes  \mathfrak{M}^*_{n-1,\K} \\
            @A  \jmath_\xi  A   A        @ | \\
           \pi_{\xi_0} = \pi_{0_{n,\K}}\widehat \otimes \pi_{0_{n-1,\K}}@>\ \oZ_{\xi_0}^\circ(\cdot, s, \chi_\K^{(j)})>> \chi_{\K, \frac{1}{2}-s}^{(j)}\otimes  \mathfrak{M}^*_{n-1,\K}\\
  \end{CD}
\]
commutes for all $s\in \BC$, where $\Omega'_{\mu,\nu, j}$ is  given by \eqref{OmegaKj} as in Theorem \ref{thmap}. 
\end{thmp}

Here  $\phi_{\xi, j}\otimes \oZ_{\xi}^\circ(\cdot, s+j, \chi_\K)$ is defined in the way similar to \eqref{phioZ}.

\subsection{Reduction to principal series representations}

Recall that in Section \ref{sec2.2} we have defined a principal series representation  $I_\mu$ with a Whittaker functional $\lambda_\mu' \in \Hom_{\mathrm N_n(\K)}(I_\mu, \psi_{n, \K})$, and a unique 
$
p_\mu\in \Hom_{\GL_{n}(\K)}(I_\mu, \pi_\mu)
$
such that 
\[
\lambda_\mu \circ p_\mu =\lambda_\mu'.
\] 
We have also defined 
$
\imath_\mu\in \Hom_{\GL_n(\K)}(I_{0_{n,\K}},  F_\mu^\vee \otimes I_\mu)
$
such that
\[
(v_\mu \otimes \lambda_\mu' )  \circ \imath_\mu =  \lambda_{0_{n,\K}}'
\]
and that the diagram \eqref{I-pi} commutes. We have similar data for $\nu$. 
Put
\[
I_\xi := I_\mu\widehat \otimes I_\nu, \quad  p_\xi:=p_\mu \otimes p_\nu, \quad
 \imath_\xi:=\imath_\mu \otimes \imath_\nu.
\]
Define the normalized Rankin-Selberg integral
\[
\begin{aligned}
 \oZ_\xi^{\diamond}(\cdot, s, \chi_\K) & = \frac{1}{\oL(s, \pi_\mu\times\pi_\nu)} \oZ_\xi(\cdot, s, \chi_\K) \\
 &  \in \Hom_{\GL_{n-1}(\K)}\left(I_\xi \otimes \frak{M}_{n-1,\K}, \,  \chi_{\K, -s+\frac{1}{2}}\right)
 \end{aligned}
 \]
 as the composition 
 \[
I_\xi \otimes \frak{M}_{n-1,\K} \xrightarrow{p_\xi\otimes{\rm id}} \pi_\xi \otimes \frak{M}_{n-1,\K}  \xrightarrow{ \oZ_\xi^{\circ}(\cdot, s, \chi_\K)} \chi_{\K, -s+\frac{1}{2}}. 
  \]
Then 
\[
 \begin{aligned}
 &   \quad \oZ^\diamond_{\xi}( f\otimes f' \otimes m, s, \chi_\K)   =    \frac{1}{\oL(s, \pi_\mu\times\pi_\nu)}    \\
&  \cdot \int_{\mathrm N_{n-1}(\K)\backslash \GL_{n-1}(\K)} \lambda'_\mu\left(\begin{bmatrix} g & 0 \\ 0 & 1 \end{bmatrix}.f\right) \lambda'_{\nu} (g. f') \cdot \chi_{\K}(\det g) \cdot \abs{\det g}_\K^{s-\frac{1}{2}}  \od\! \bar{m} (g), 
\end{aligned}
\]
for $f\in I_\mu$, $f' \in I_\nu$, $ m \in \frak{M}_{n-1,\K}$,  and $s\in \C$ with  ${\rm Re}(s)$ sufficiently large.

In view of all the above, by multiplicity one theorem \cite{AG2, SZ}, there exists a unique entire function $\Xi_{\mu, \nu, j}(s)$  such that the diagram  
\be \label{cube}
\xymatrix{
  F_\xi^\vee \otimes I_\xi \ar[rd]^(0.6){{\rm id}\otimes p_\xi} \ar[rrr]^{ \phi_{\xi, j}\otimes \oZ_\xi^\diamond(\cdot, s+j,  \chi_\K)}  
  &&&  \chi_{\K, \frac{1}{2}-s}^{(j)}\otimes \frak{M}^*_{n-1,\K}  \ar@{=}[rd]\\
    & F_\xi^\vee \otimes \pi_\xi \ar[rrr]^(0.3){  \phi_{\xi, j}\otimes \oZ_\xi^\circ(\cdot, s+j,  \chi_\K)} 
    &&&   \chi_{\K, \frac{1}{2}-s}^{(j)}\otimes  \frak{M}^*_{n-1,\K} \\
    I_{\xi_0} \ar[rd]_{p_{\xi_0}}\ar[uu]^{\imath_\xi} \ar[rrr]^(0.6){\oZ_{\xi_0}^\diamond(\cdot, s,  \chi_\K^{(j)})} |!{[ur];[dr]}\hole
    && &  \chi_{\K, \frac{1}{2}-s}^{(j)}\otimes  \frak{M}^*_{n-1,\K} \ar@{=}[rd] \ar^(0.4){\Xi_{\mu,\nu, j}(s) }[uu]|!{[ull];[ur]}\hole  \\
    & \pi_{\xi_0} \ar[uu]^(0.6){\jmath_\xi} \ar[rrr]^{ \oZ_{\xi_0}^\circ(\cdot, s,  \chi_\K^{(j)})} 
    &&&  \chi_{\K, \frac{1}{2}-s}^{(j)}\otimes \frak{M}^*_{n-1,\K} \ar^{\Xi_{\mu,\nu, j}(s) }[uu]\\
}
\ee
commutes for all $s\in \BC$.

In the rest of this section we compute the function $\Xi_{\mu, \nu, j}(s)$ and show that it is a nonzero constant whose inverse is equal to $\Omega'_{\mu, \nu, j}$ given by \eqref{OmegaKj}. The main ingredient of the computation is \cite{LLSS}.

\subsection{Integral over the open orbit} \label{sec:lss} For the convenience of the reader, we describe the main result of \cite{LLSS}. Write $\widehat{\K^\times}$ for the set of all (unitary or not) characters of $\K^\times$.
Let $\varrho=(\varrho_1,\ldots, \varrho_n) \in (\widehat{\K^\times})^n$, viewed as a character of $\bar{\mathrm B}_n(\K)$ as usual, and let
 \[
 I(\varrho) :=\Ind^{\GL_n(\K)}_{\bar{\rm B}_n(\K)} \varrho
 \]
 be the corresponding principal series representation of $\GL_n(\K)$. Similarly let $\varrho' = (\varrho_1',\ldots, \varrho_{n-1}')\in (\widehat{\K^\times})^{n-1}$ and let $I(\varrho')$ be the corresponding principal series representation of $\GL_{n-1}(\K)$.  We have a meromorphic family of unnormalized Rankin-Selberg integrals
\[
 \oZ(\cdot, s, \chi_\K)  
  \in \Hom_{\GL_{n-1}(\K)}\left(I(\varrho)\widehat{\otimes} I(\varrho') \otimes \frak{M}_{n-1,\K}, \,  \chi_{\K, -s+\frac{1}{2}}\right)
 \]
 such that
 \[
 \begin{aligned}
&  \oZ( f\otimes f' \otimes  m, s, \chi_\K)   \\
 = \ & \int_{\mathrm N_{n-1}(\K)\backslash \GL_{n-1}(\K)} \lambda'_\varrho\left(\begin{bmatrix} g & 0 \\ 0 & 1 \end{bmatrix}.f\right) \lambda'_{\varrho'} (g. f') \cdot \chi_{\K}(\det g) \cdot \abs{\det g}_\K^{s-\frac{1}{2}}  \od\! \bar{m} (g)
 \end{aligned}
\]
for all $f\in I(\varrho)$, $f' \in I(\varrho')$, $m \in \frak{M}_{n-1,\K}$, and $s\in \C$ with  ${\rm Re}(s)$ sufficiently large, where $\lambda'_\varrho\in \Hom_{\RN_n(\K)}(I(\varrho), \psi_{n,\K})$ and 
$\lambda'_{\varrho'}\in \Hom_{\RN_{n-1}(\K)}(I(\varrho'), \psi_{n-1,\K})$
 are defined in the way similar to \eqref{la'}.
 
 Let
\be\label{z}
z:=(z_n, z_{n-1})\in \GL_n(\mathbb{Z})\times \GL_{n-1}(\BZ),
\ee
where $z_n\in \GL_n(\Z)$ is defined inductively in \eqref{zk}. 
 The right action of $\GL_{n-1}(\K)$ on the flag variety $(\bar{\mathrm B}_n(\K) \times \bar{\mathrm B}_{n-1}(\K) )\backslash (\GL_n(\K) \times \GL_{n-1}(\K)) $ has a unique  open orbit 
 \be \label{openorbit}
\left(\big(\bar{\mathrm B}_n(\K)\times \bar{\mathrm B}_{n-1}(\K)\big)   z\right) \cdot \GL_{n-1}(\K).
 \ee
Note that  
\[
I(\varrho)\widehat{\otimes} I(\varrho') = \Ind^{\GL_n(\K)\times \GL_{n-1}(\K)}_{\bar{\rm B}_n(\K)\times \bar{\rm B}_{n-1}(\K)}\varrho\otimes \varrho'.
\]
Following \cite{LLSS}, we first formally define
\[
\Lambda(\cdot, s, \chi_\K )\in \Hom_{\GL_{n-1}(\K)}\left(I(\varrho) \widehat \otimes I(\varrho') \otimes \frak{M}_{n-1,\K}, \,  \chi_{\K, -s+\frac{1}{2}}\right)
\]
as the integral over the above open orbit, that is, 
\be \label{lambda}
\begin{aligned}
 & \quad  \Lambda( \phi \otimes m, s, \chi_\K)   \\
:= &     \int_{  \GL_{n-1}(\K)}\phi \left(z_n \left[
                                                                \begin{array}{cc}
                                                                  g & 0 \\
                                                                  0 & 1 \\
                                                                \end{array}
                                                              \right], z_{n-1} g\right)\cdot \chi_\K(\det g) \cdot \abs{\det g}^{s-\frac{1}{2}}_{\K}\od\! m(g)
\end{aligned}
\ee
for $ \phi\in  I(\varrho) \widehat\otimes I(\varrho')$ and $m\in \frak{M}_{n-1,\K}$.

Define
\[
\sgn(\varrho, \varrho', \chi_\K):=\prod_{i>k,  \, i+k\leq n} (\varrho_i \cdot \varrho'_k \cdot\chi_\K)(-1),
\]
and a meromorphic function 
\[
\gamma_{\psi_\K^{(n)}}(s, \varrho, \varrho', \chi_\K):=\prod_{i+k\leq n} \gamma(s, \varrho_i \cdot \varrho'_k \cdot \chi_\K, \psi_\K^{(n)}),
\]
where $\psi_\K^{(n)}$ is the additive character $\K\to \BC^\times$, $x\mapsto \psi_\K((-1)^n x)$, 
\[
\gamma(s, \omega, \psi_\K^{(n)}) = \varepsilon(s,\omega, \psi_\K^{(n)}) \cdot \frac{\oL(1-s, \omega^{-1})}{\oL(s, \omega)}
\]
is the local gamma factor of a character $\omega \in \widehat{\K^\times}$, and $\varepsilon(s, \omega, \psi_\K^{(n)})$ is the local epsilon factor, defined following \cite{T, J1, Ku}.
For convenience also define
\[
\varepsilon_{\psi_\K^{(n)}}(s, \varrho, \varrho', \chi_\K):=\prod_{i+k\leq n} \varepsilon(s, \varrho_i \cdot \varrho'_k\cdot  \chi_\K, \psi_\K^{(n)}).
\]
Finally define a meromorphic function 
\[
\Gamma_{\psi_\K^{(n)}}(s, \varrho, \varrho', \chi_\K):=\sgn(\varrho, \varrho', \chi_\K)  \cdot \gamma_{\psi_\K^{(n)}}(s, \varrho, \varrho', \chi_\K).
\]

For a character $\omega\in \widehat{\K^\times}$, denote by ${\rm ex}(\omega)$ the real number such that 
\[
\abs{\omega} = \abs{\,\cdot\,}_\K^{{\rm ex}(\omega)}.
\] 
Consider the complex manifold 
\[
\mathcal{M}:=\BC\times (\widehat{\K^\times})^n \times (\widehat{\K^\times})^{n-1} 
\]
and its nonempty open subset 
\[
\Omega:=\left\{ (s, \varrho, \varrho') \in \mathcal{M} \left|\, \begin{aligned} & \mathrm{ex}(\varrho_i)+\mathrm{ex}(\varrho'_k)+\Re(s)<1  \textrm{ whenever } i+k\leq n,\\
& \mathrm{ex}(\varrho_i)+\mathrm{ex}(\varrho'_k)+\Re(s)>0 \textrm{ whenever } i+k> n
\end{aligned}\right.\right\}.
\]

\begin{thmp} \cite[Theorem 1.6 (b)]{LLSS} \label{thm:LSS}
Assume that $(s,\varrho,\varrho')\in\Omega$. Then the integral \eqref{lambda} converges absolutely, and 
\be \label{eq:LSS}
 \Lambda ( \phi \otimes m, s, \chi_\K)   = \Gamma_{\psi_\K^{(n)}}(s, \varrho, \varrho', \chi_\K)  \cdot \oZ( \phi \otimes m, s, \chi_\K). 
\ee
\end{thmp}
We remark that the right hand side of \eqref{eq:LSS} is  holomorphic as a  function of the variable $s\in \Omega_{\varrho, \varrho'}:=\{s\in \C\, : \,  (s, \varrho, \varrho')\in \Omega\}$ (see \cite[Remark 1.7]{LLSS}).

  Let $(I(\varrho) \widehat\otimes I(\varrho'))^\sharp \subset I(\varrho) \widehat\otimes I(\varrho')$ be the subspace of $\phi \in  I(\varrho) \widehat\otimes I(\varrho')$ such that 
 \[
\phi |_{z\cdot \GL_{n-1}(\K)}\in
 \mathcal{S}(z\cdot \GL_{n-1}(\K)).
 \]
 Then for every $\varrho\in (\widehat{\K^\times})^n$, $\varrho'\in (\widehat{\K^\times})^{n-1}$ and $\phi \in (I(\varrho)\widehat\otimes I(\varrho'))^\sharp$, the integral \eqref{lambda} converges absolutely and is an entire function of $s\in\BC$. We deduce the following consequence of Theorem \ref{thm:LSS}.
 
 \begin{corp} \label{cor:LSS}
 For every $\varrho\in (\widehat{\K^\times})^n$, $\varrho'\in (\widehat{\K^\times})^{n-1}$ and $\phi \in (I(\varrho)\widehat\otimes I(\varrho'))^\sharp$, the equality \eqref{eq:LSS} holds as entire functions of $s\in\BC$.
 \end{corp}
 
\begin{proof}
Let $\CC$ be the connected component of $(\widehat{\K^\times})^n$ containing $\varrho$, and let $\CC'$ be the  connected component of $(\widehat{\K^\times})^{n-1}$ containing $\varrho'$.
Write $K_\K:=K_{n,\K}\times K_{n-1,\K}$.  Define
\[
 C^\infty_{\CC, \CC'}(K_\K):=\left\{ f \in C^\infty(K_\K) \left| \begin{array}{l} f(b\cdot k) = ( \varrho \otimes \varrho')(b) \cdot f(k), \\
  \textrm{for all } b\in K_\K\cap (\bar{\mathrm B}_n(\K) \times \bar{\mathrm B}_{n-1}(\K)), \ k\in K_\K
  \end{array}\right.\right\},
\] 
which only depends on $\CC$ and $\CC'$, not on the particular choices of $\varrho$ and $\varrho'$.

Consider the natural  map 
\[
K_\K \to (\bar{\mathrm B}_n(\K) \times \bar{\mathrm B}_{n-1}(\K) )\backslash (\GL_n(\K) \times \GL_{n-1}(\K)),
\]
which is surjective by the Iwasawa decomposition. Let $K_\K^\sharp \subset K_\K$ be the preimage of the open orbit  \eqref{openorbit} under the above map. Fix $f \in C^\infty_{\CC, \CC'}(K_\K)$ such that 
$f|_{K_\K^\sharp}\in \CS(K_\K^\sharp)$. Then there is a unique  
\[
\phi_{\varrho,\varrho'}:=\phi_{f,\varrho,\varrho'}\in (I(\varrho)\widehat\otimes I(\varrho'))^\sharp
\]
such that 
\[
\phi_{\varrho, \varrho'}|_{K_\K} = f.
\]

Let $\CM^\circ:=\BC\times \CC\times \CC'$, which is a connected component of $\CM$. When $(\varrho, \varrho')$ varies in $\CC\times \CC' $, the integral  
$\Lambda(\phi_{\varrho,\varrho'}\otimes m, s, \chi_\K)$ is clearly holomorphic on $\CM^\circ$. 
By \cite[Section 8.1]{J2}, we also have that
\[
 \Gamma_{\psi_\K^{(n)}}(s, \varrho, \varrho', \chi_\K)  \cdot \oZ( \phi_{\varrho,\varrho'} \otimes m, s, \chi_\K)  
 \]
is meromorphic on $\CM^{\circ}$.  Since the equality 
\[
\Lambda(\phi_{\varrho,\varrho'}\otimes m, s, \chi_\K)  =  \Gamma_{\psi_\K^{(n)}} (s, \varrho, \varrho', \chi_\K)  \cdot \oZ( \phi_{\varrho,\varrho'} \otimes m, s, \chi_\K) 
\]
holds on $\Omega\cap \CM^\circ$, which is nonempty and open, it holds over all $\CM^\circ$ by the uniqueness of meromorphic  continuation. The corollary then follows by noting that every $\phi\in  (I(\varrho)\widehat\otimes I(\varrho'))^\sharp$ equals $\phi_{f,\varrho, \varrho'}$ for some $f \in C^\infty_{\CC, \CC'}(K_\K)$ such that 
$f|_{K_\K^\sharp}\in \CS(K_\K^\sharp)$.
\end{proof}

\subsection{A commutative diagram} We now specify the above discussion  to the principal series representations $I_\mu$ and $I_\nu$.  
Define $\varrho^\mu = (\varrho^\mu_1, \ldots, \varrho^\mu_n)\in (\widehat{\K^\times})^n$, where
\[
\varrho^\mu_i  := \varepsilon_{n,\K} \abs{\,\cdot\,}_\K^{\frac{n+1}{2}-i} \prod_{\iota \in \CE_\K} \iota^{\mu_i^{\iota}}\in \widehat{\K^\times}, \quad i =1,\ldots, n,
\]
so that
$
I_\mu = I(\varrho^\mu)
$
in the above notation. Likewise we define $\varrho^\nu = (\varrho^\nu_1,\ldots, \varrho^\nu_{n-1})\in (\widehat{\K^\times})^{n-1}$ so that $I_\nu = I(\varrho^\nu)$, and
put
$I_\xi^\sharp: = (I_\mu \widehat\otimes I_\nu)^\sharp$.  Similar to $\varepsilon_{n,\K}$,  $\varepsilon_{n-1,\K}$ denotes the common central character of $F_\nu^\vee\otimes \pi_\nu$ and $\pi_{0_{n-1,\K}}$.

The integral \eqref{lambda} defines a nonzero linear functional 
\[
\Lambda_\xi (\cdot, s, \chi_\K)\in \Hom_{\GL_{n-1}(\K)}\left(I_\xi^\sharp, \,  \chi_{\K, -s+\frac{1}{2}} \otimes \frak{M}_{n-1,\K}^* \right).
\]
By Corollary \ref{cor:LSS},
\be \label{LLSS}
\begin{aligned}
\Lambda_\xi (\cdot, s, \chi_\K) & = \Gamma_{\psi_\K^{(n)}} (s, \varrho^\mu, \varrho^\nu, \chi_\K) \cdot \oZ_\xi(\cdot, s, \chi_\K) \\
& = \Gamma_{\psi_\K^{(n)}} (s, \varrho^\mu, \varrho^\nu, \chi_\K) \cdot \oL(s, \pi_\mu\times \pi_\nu) \cdot \oZ^\diamond_\xi(\cdot, s, \chi_\K)
\end{aligned}
\ee
holds on $I_\xi^\sharp\otimes \frak{M}_{n-1,\K}$. Recall that $\imath_\xi = \imath_\mu \otimes \imath_\nu \in \Hom_{\GL_{n-1}(\K)}(I_{\xi_0}, F_{\xi}^\vee\otimes I_\xi)$. It is clear that 
\[
\imath_\xi(I_{\xi_0}^\sharp)\subset F_\xi^\vee\otimes I_\xi^\sharp.
\]

\begin{prp} \label{prop3.2}
The diagram 
\[
 \begin{CD}
            F_\xi^\vee  \otimes I_\xi^\sharp
                  @> \phi_{\xi,j}\otimes \Lambda_\xi(\cdot, s+j, \chi_\K) >> \chi_{\K, \frac{1}{2}-s}^{(j)}\otimes  \mathfrak{M}^*_{n-1,\K} \\
            @A  \imath_\xi  A A          @ |  \\
          I_{\xi_0}^\sharp @>\Lambda_{\xi_0}(\cdot, s, \chi_\K^{(j)}) >> \chi_{\K, \frac{1}{2}-s}^{(j)}\otimes  \mathfrak{M}^*_{n-1,\K}\\
  \end{CD}
\]
commutes. 
\end{prp}

\begin{proof}
Recall from Proposition \ref{phixij}  that  $\phi_{\xi,j}\in\Hom_{\GL_{n-1}(\K\otimes_\R \C)}(F_\xi^\vee , \otimes_{\iota\in \mathcal E_\K}{\det}^j)$ and 
\[
\phi_{\xi, j}(z^{-1}. v_\xi^\vee)=1,
\]
where $v_\xi^\vee:=v_\mu^\vee\otimes v_\nu^\vee$. For $\phi \in I_{\xi_0}$ and $g\in \GL_{n-1}(\K)\subset \GL_n(\K)\times \GL_{n-1}(\K)$, we have that
\[
\begin{aligned}
 &  \phi_{\xi, j} ( \imath_\xi(\phi) (zg)) \\
 = \ & \phi_{\xi, j}\left(\phi(zg)\cdot (g^{-1}z^{-1}. v_\xi^\vee)\right) \quad (\textrm{see } \eqref{imu}) \\
= \ &  \phi(zg)\cdot \phi_{\xi, j}(g^{-1}z^{-1}. v_\xi^\vee)  \\
= \ &  \phi(zg) \cdot (\otimes_{\iota\in \mathcal E_\K}{\det}^{-j})(g).
\end{aligned}
\]
Assume that $\phi\in I_{\xi_0}^\sharp$. By  \eqref{lambda}, we have that
\[
\begin{aligned}
& \left( \phi_{\xi,j}\otimes \Lambda_\xi(\cdot, s+j, \chi_\K)\right)(\imath_{\xi}(\phi)\otimes m) \\ 
= \ & \int_{  \GL_{n-1}(\K)} \phi_{\xi, j} ( \imath_\xi(\phi) (zg)) \cdot \chi_\K(\det g) \cdot \abs{\det g}^{s+j-\frac{1}{2}}_{\K}\od\! m(g) \\
= \ & \int_{  \GL_{n-1}(\K)}  \phi(zg) \cdot (\otimes_{\iota\in \mathcal E_\K}{\det}^{-j})(g)  \cdot \chi_\K(\det g) \cdot \abs{\det g}^{s+j-\frac{1}{2}}_{\K}\od\! m(g) \\
= \ &  \int_{  \GL_{n-1}(\K)}  \phi(zg) \cdot   \chi_\K^{(j)}(\det g) \cdot \abs{\det g}^{s-\frac{1}{2}}_{\K}\od\! m(g) \\
= \ & \Lambda_{\xi_0}(\phi\otimes m, s, \chi_\K^{(j)}),
\end{aligned}
\]
where $m\in \frak{M}_{n-1,\K}$. This proves the proposition. 
\end{proof}

\begin{corp}
Let the notations be as above. Then
\be \label{const}
\Xi_{\mu, \nu, j}(s) \cdot \frac{\Gamma_{\psi_\K^{(n)}} (s+j, \varrho^\mu, \varrho^\nu, \chi_\K)}{\Gamma_{\psi_\K^{(n)}} (s, \varrho^{0_{n,\K}}, \varrho^{0_{n-1,\K}}, \chi_\K^{(j)})} \cdot 
\frac{\oL(s+j, \pi_\mu\times \pi_\nu)}{\oL(s, \pi_{0_{n,\K}}\times \pi_{0_{n-1,\K}})}=1
\ee
as meromorphic functions of the variable  $s\in \BC$.
\end{corp}

\begin{proof}
This follows from \eqref{cube}, \eqref{LLSS} and Proposition \ref{prop3.2}. 
\end{proof}

\subsection{Archimedean local factors} To finish the proof, it remains to evaluate the function $\Xi_{\mu, \nu, j}(s)^{-1}$ given by \eqref{const} and show that it is equal to the constant $\Omega'_{\mu, \nu, j}$ given by  \eqref{OmegaKj} as in Theorem \ref{thmap}. To this end, we first recall some standard facts about archimedean local L-factors and epsilon factors that we need, following \cite{K}. Let
\[
\Gamma_\K(s)=\left\{ \begin{array}{ll} \pi^{-s/2}\Gamma(s/2), & \textrm{if }\K\cong\BR; \\
2(2\pi)^{-s}\Gamma(s), & \textrm{if }\K\cong\BC, 
\end{array}
\right.
\]
where $\Gamma(s)$ is the standard gamma function. Recall the Legendre duplication formula
\be \label{dup}
\Gamma_\BC(s) = \Gamma_\BR(s) \Gamma_\BR(s+1).
\ee
Recall the additive character   $\psi_\K$ that is defined in \eqref{psiK} by using $\psi_\BR$.

If $\K\cong\BR$, then the followings hold true. 
\begin{itemize}
\item For all $t\in\BC$ and $\delta\in\{0,1\}$,
\[
\oL(s, \abs{\,\cdot\,}_\K^t \sgn_{\K^\times}^\delta)=  \Gamma_\BR(s+t +\delta),
\]
and
\[
\varepsilon(s,  \abs{\,\cdot\,}_\K^t \sgn_{\K^\times}^\delta, \psi_\K^{(n)})=((-1)^n  {\rm i})^\delta;
\]
\item For all  $a, b\in \BC$ with $a-b\in\BZ\setminus\{0\}$, and $t$, $\delta$ as above,
\[
\oL(s, D_{a, b} \times \abs{\,\cdot\,}_\K^t \sgn_{\K^\times}^\delta)=  \Gamma_\BC\left(s+t+\max\{a,b\}\right).
\]
Here and henceforth, for any $a', b'\in \BC$ with $a'-b'\in \BZ$,
\[
\max\{a', b'\}: = \left\{\begin{array}{ll}
a', & \textrm{if }a'-b' \geq 0; \\
b', & \textrm{otherwise}. \end{array} \right.
\]
\item For all $a, b, a', b'\in \BC$ with $a-b, a'-b'\in \BZ\setminus\{0\}$,
\[
\oL(s, D_{a, b} \times D_{a', b'})= \Gamma_\BC\left(s+\max\{a+a', b+b'\}\right) \Gamma_\BC\left(s+ \max\{a+b', b+a'\}\right).
\]
\end{itemize}
If $\K\cong\BC$, then for all $a, b\in \BC$ with $a-b\in \BZ$, 
\[
\oL(s, \iota^a \bar{\iota}^b )= \Gamma_\BC\left(s+\max\{a, b\}\right)
\]
and
\[
\varepsilon(s, \iota^a \bar{\iota}^b, \psi_\K^{(n)}) = ((-1)^n \mathrm{i})^{\abs{a-b}}. 
\]

By the well-known branching rule for $\GL_n(\BC)$, we have that $j\in \BZ$ is a balanced place for $\xi$ if and only if 
\[
-\mu^{\iota}_n \geq \nu^{\iota}_1+ j  \geq - \mu^{\iota}_{n-1} \geq \nu^{\iota}_2+j \geq \cdots \geq \nu^{\iota}_{n-1}+j \geq - \mu^{\iota}_1
\]
for every $\iota\in \CE_\K$.  Equivalently, by  \cite[Corollary 2.35]{Rag2},  $j\in \mathbb{Z}$ is a balanced  place for $\xi$ if and only if 
\[
m_{\mu, \nu}^-\leq j\leq  m_{\mu, \nu}^+, 
\]
where
\[
\begin{aligned}
& m_{\mu, \nu}^- := \max\{-\mu_{n-i}^{\iota}-\nu_i^{\iota}: 1\leq i\leq n-1, \iota \in \CE_\K\}, \\
& m_{\mu, \nu}^+ := \min\{-\mu_{n+1-i}^{\iota}-\nu_i^{\iota}: 1\leq i\leq n-1, \iota \in \CE_\K\}.
\end{aligned}
\]
Recall the half-integers $\tilde{\mu}^\iota_i$ ($i=1,\ldots, n$) given by \eqref{lmu}. The following result can be easily checked by using the above result (see the proof of \cite[Lemma 2.24]{Rag2}).

\begin{lemp} \label{lem:Rag}
Assume that $\xi=(\mu,\nu)$ is balanced. Then
\[
\tilde{\mu}^\iota_i + \tilde{\nu}^\iota_k - \tilde{\mu}^{\bar{\iota}}_{n+1-i}  - \tilde{\nu}^{\bar{\iota}}_{n-k}
\]
is positive if $i+k\leq n$, and is negative otherwise, for every $\iota\in\CE_\K$.
\end{lemp}

Now we  establish  the following result, which thereby finishes the proof of Theorem \ref{thmap2}.

\begin{prp}
The function $\Xi_{\mu, \nu, j}(s)^{-1}$ given by \eqref{const} is equal to the constant $\Omega'_{\mu, \nu, j}$ given by \eqref{OmegaKj}.
 \end{prp}

\begin{proof} We use the notation of Theorem \ref{thmap}. It is clear that
\[
 \frac{\sgn(\varrho^\mu, \varrho^\nu, \chi_\K)}{\sgn(\varrho^{0_{n,\K}}, \varrho^{0_{n-1,\K}}, \chi_\K^{(j)})} 
 = \prod_{i>k, \,  i+k \leq n}(-1)^{\sum_{\iota\in \CE_\K}(\mu_i^\iota+\nu_k^\iota+j)} =(-1)^{ j \frac{n^2(n-1)}{2}[\K:\, \BR]}\cdot \varepsilon_{\mu, \nu}.
\]
To evaluate the contribution from the local gamma and $\oL$-factors, we consider the real and complex cases separately. 

(i) Assume that $\K\cong \BR$. Then $\CE_\K=\{\iota\}$. Using Lemma \ref{lem:Rag}, it is easy to check that
\[
\oL(s, \pi_\mu\times \pi_\nu) =  \prod_{i+k\leq n} \Gamma_\BC(s + \tilde{\mu}^\iota_i + \tilde{\nu}^\iota_k).
\]
For a character $\omega\in \widehat{\K^\times}$, write $\delta(\omega)\in \{0, 1\}$ such that $\omega(-1) = (-1)^{\delta(\omega)}$. Then 
\be\label{R0}
\varepsilon(s, \omega, \psi_\K^{(n)})= ((-1)^n  \mathrm{i})^{\delta(\omega)},
\ee
and it is clear that
\be\label{R1}
j+ \mu^\iota_i + \nu^\iota_k  - \delta(\varrho^\mu_i \varrho^\nu_k \chi_\K)  +  \delta(\varrho^{0_{n,\K}}_i \varrho^{0_{n-1,\K}}_k \chi_\K^{(j)}) \in 2\BZ.
\ee 
We have that
\[
\frac{\oL(1-s, (\varrho^\mu_i \varrho^\nu_k\chi_\K)^{-1})}{\oL(s, \varrho^\mu_i\varrho^\nu_k \chi_\K)} = \frac{\Gamma_\BR(1-s-\tilde{\mu}^{\iota}_i - \tilde{\nu}^\iota_k+\delta(\varrho^\mu_i \varrho^\nu_k \chi_\K))}{\Gamma_\BR(s+\tilde{\mu}^\iota_i + \tilde{\nu}^{\iota}_k+\delta(\varrho^\mu_i \varrho^\nu_k \chi_\K))}.
\]
It follows from \eqref{dup} that
\be \label{R2}
\begin{aligned}
& \left( \prod_{i+k\leq n} \frac{\oL(1-s, (\varrho^\mu_i \varrho^\nu_k \chi_\K )^{-1})}{\oL(s, \varrho^\mu_i\varrho^\nu_k\chi_\K)} \right)  \cdot  \oL(s, \pi_\mu\times \pi_\nu)\\
= \ & \prod_{i+k\leq n}\left(\Gamma_\BR(s+\tilde{\mu}^\iota_i + \tilde{\nu}^\iota_k+1-\delta(\varrho^\mu_i \varrho^\nu_k \chi_\K))
\cdot\Gamma_\BR(1-s- \tilde{\mu}^\iota_i  - \tilde{\nu}^\iota_k+ \delta(\varrho^\mu_i \varrho^\nu_k \chi_\K))\right).
\end{aligned}
\ee
By \eqref{R0}, \eqref{R1}, \eqref{R2} and the formula
\[
\Gamma_\BR(s+  \ell )\cdot \Gamma_\BR(2-s- \ell ) = {\rm i}^\ell \cdot \Gamma_\BR(s)\cdot \Gamma_\BR(2-s),\quad \ell \in 2 \BZ, 
\]
we find that
\[
\begin{aligned}
& \frac{\gamma_{\psi_\K^{(n)}} (s+j, \varrho^\mu, \varrho^\nu, \chi_\K)}{\gamma_{\psi_\K^{(n)}} (s, \varrho^{0_{n,\K}}, \varrho^{0_{n-1,\K}}, \chi_\K^{(j)})} \cdot 
\frac{\oL(s+j, \pi_\mu\times \pi_\nu)}{\oL(s, \pi_{0_{n,\K}}\times \pi_{0_{n-1,\K}})} \\
= \ &  \prod_{i+k\leq n} \left(\frac{ \varepsilon(s+j, \varrho^\mu  \varrho^\nu \chi_\K, \psi_\K^{(n)})}{\varepsilon (s, \varrho^{0_{n,\K}} \varrho^{0_{n-1,\K}} \chi_\K^{(j)}, \psi_\K^{(n)})} \cdot {\rm i}^{j+ \mu_i^\iota+\nu^\iota_k - \delta(\varrho^\mu  \varrho^\nu \chi_\K) + \delta(\varrho^{0_{n,\K}} \varrho^{0_{n-1,\K}} \chi_\K^{(j)})} \right) \\
= \ & \prod_{i+k\leq n}  ((-1)^n   \mathrm{i})^{j+\mu^\iota_i + \nu^\iota_k} \\
= \ &  (-1)^{j \frac{n^2(n-1)}{2}}\cdot   \mathrm{i}^{j\frac{n(n-1)}{2}} \cdot c_\mu' \cdot c_\nu.
\end{aligned}
\]

(ii) Assume that $\K\cong \BC$. Then $\chi_\K^{(j)}$ is trivial, which will be omitted from the notations for convenience. Using Lemma \ref{lem:Rag} again, we find that
\[
\oL(s, \pi_\mu\times \pi_\nu) = \prod_{i+k\leq n, \, \iota\in \CE_\K} \Gamma_\BC(s + \tilde{\mu}^\iota_i +\tilde{\nu}^\iota_k).
\]
We have that
\[
\frac{\oL(1-s, (\varrho^\mu_i \varrho^\nu_k)^{-1})}{\oL(s, \varrho^\mu_i\varrho^\nu_k)} = \frac{\Gamma_\BC(1-s-\min_{\iota \in\CE_\K}\{\tilde{\mu}^{\iota}_i + \tilde{\nu}^\iota_k\})}{\Gamma_\BC(s+\max_{\iota\in \CE_\K}\{\tilde{\mu}^\iota_i + \tilde{\nu}^{\iota}_k\})}.
\]
It follows that
\be \label{C1}
\begin{aligned}
& \left(\prod_{i+k\leq n} \frac{\oL(1-s, (\varrho^\mu_i \varrho^\nu_k)^{-1})}{\oL(s, \varrho^\mu_i\varrho^\nu_k)} \right) \cdot  \oL(s, \pi_\mu\times \pi_\nu)\\
= \ & \prod_{i+k\leq n}\left( \Gamma_\BC(s+\min_{\iota\in\CE_\K}\{\tilde{\mu}^\iota_i + \tilde{\nu}^\iota_k\})
\cdot\Gamma_\BC(1-s-\min_{\iota\in\CE_\K}\{\tilde{\mu}^\iota_i + \tilde{\nu}^\iota_k\})\right).
\end{aligned}
\ee
Using \eqref{C1} and the formula
\[
\Gamma_\BC(s+ \ell)\cdot\Gamma_\BC(1-s-\ell) = (-1)^\ell \cdot \Gamma_\BC(s)\cdot \Gamma_\BC(1-s),\quad \ell \in \BZ,
\]
we find that
\be \label{C2}
\begin{aligned}
& \frac{\gamma_{\psi_\K^{(n)}} (s+j, \varrho^\mu, \varrho^\nu)}{\gamma_{\psi_\K^{(n)}} (s, \varrho^{0_{n,\K}}, \varrho^{0_{n-1,\K}})} \cdot 
\frac{\oL(s+j, \pi_\mu\times \pi_\nu)}{\oL(s, \pi_{0_{n,\K}}\times \pi_{0_{n-1,\K}})} \\
= \ & \prod_{i+k\leq n} \left( \varepsilon(s+j, \varrho^\mu_i \varrho^\nu_k, \psi_\K^{(n)}) \cdot (-1)^{j+\min_{\iota\in\CE_\K} \{ \mu^\iota_i + \nu^\iota_k\}}\right).
\end{aligned}
\ee
We have the local epsilon factor
\[
\varepsilon(s+j, \varrho^\mu_i \varrho^\nu_k, \psi_\K^{(n)}) = ((-1)^n  \mathrm{i})^{\max_{\iota\in\CE_\K}\{\mu^\iota_i + \nu^\iota_k \} - \min_{\iota\in \CE_\K} \{ \mu^{ \iota}_i + \nu^{ \iota}_k\} }.
\]
Hence \eqref{C2} is equal to 
\[
\prod_{i+k \leq n, \, \iota\in \CE_\K}  ((-1)^n \mathrm{i})^{j + \mu^\iota_i + \nu^\iota_k } =  (-1)^{jn^2(n-1)}\cdot  \mathrm{i}^{j n(n-1)}\cdot c'_\mu \cdot c_\nu.
\]
This finishes the proof of the proposition. 
\end{proof}

\section{Non-archimedean period relations} \label{secnap}

In this section, let $\K$ be a non-archimedean local field of characteristic zero. Fix a non-trivial unitary character $\psi_\K: \K\to \BC^\times$, and define the character $\psi_{n,\K}$ of $\mathrm{N}_n(\K)$ as in 
\eqref{psi_n} ($n\geq 1$).

\subsection{Preliminaries}

Let $p$ be the residue characteristic of $\K$, and $\mu_{p^\infty}\subset \BC^\times$ be the subgroup of $p$-th power roots of unity. Recall the cyclotomic character 
\[
{\rm Aut}(\BQ(\mu_{p^\infty})/\BQ)\to \BZ_p^\times, \quad \sigma\mapsto t_{\sigma, p}
\]
defined by requiring that
\be \label{cyclo}
\sigma( \zeta )=\zeta^{t_{\sigma, p}} \quad \textrm{for all } \zeta \in \mu_{p^\infty}.
\ee
Write $\sigma\mapsto t_{\sigma, \K}$ for the composition 
\[
\Aut(\BC/\BQ)  \xrightarrow{\rm restriction} \Aut(\BQ(\mu^\infty_p)/\BQ) \xrightarrow{\sigma\mapsto t_{\sigma, p}}\BZ_p^\times \subset \K^\times. 
\]
Following \cite[page 79--80]{H}
    and \cite[page 594]{Mah}, define 
     \be\label{tnk}
\mathbf{t}_{n, \sigma, \K} := {\rm diag}(t^{-(n-1)}_{\sigma, \K}, \ldots, t_{\sigma, \K}^{-1}, 1) \in \GL_n(\K),
\ee
 and define an action of 
      $\Aut(\BC)$ on   $\Ind^{\GL_n(\K)}_{{\rm N}_n(\K)} \psi_{n,\K}$ (the smooth induction) by
\be \label{sigmaf}
 {}^\sigma \! f(g):=\sigma\left( f(\mathbf{t}_{n, \sigma,\K} \cdot g)  \right), 
\ee
where $\sigma\in \Aut(\BC)$, $f\in \Ind^{\GL_n(\K)}_{{\rm N}_n(\K)}\psi_{n,\K}$, and $g\in \GL_n(\K)$. 

 Let $\Pi_\K$  be a generic irreducible smooth representations of $\GL_n(\K)$, with a fixed Whittaker functional 
\[
\lambda_\K \in \Hom_{\mathrm{N}_n(\K)}(\Pi_\K, \psi_{n, \K})\setminus\{0\}. 
\]
Using $\lambda_\K$, we realize $\Pi_\K$ as a subrepresentation of $\Ind^{\GL_n(\K)}_{{\rm N}_n(\K)} \psi_{n,\K}$ by
\be \label{whit-space}
\Pi_\K \to \Ind^{\GL_n(\K)}_{{\rm N}_n(\K)} \psi_{n,\K},\quad u\mapsto \left(g\mapsto \lambda_\K(g.u)\right).
\ee Put
\[
{}^\sigma\Pi_\K: = \sigma(\Pi_\K) \subset \Ind^{\GL_n(\K)}_{{\rm N}_n(\K)}\psi_{n,\K},
\] 
which is also a generic irreducible smooth representations of $\GL_n(\K)$ with a fixed Whittaker functional (the evaluation map at the identity matrix). 

Let $\chi_\K: \K^\times \to \BC^\times$ be a  character. Let $\frak{c}(\chi_\K)$ and $\frak{c}(\psi_\K)$ be the conductors of $\chi_\K$ and $\psi_\K$, which are respectively ideal and fractional ideal of $\CO_\K$. Here $\CO_\K$ denotes the ring of integers of $\K$.  Fix $y_\K \in \K^\times$ such that
\[
\frak{c}(\psi_\K)= y_\K \cdot \frak{c}(\chi_\K).
\]
The local Gauss sum is defined by
\be \label{gauss}
\CG(\chi_\K):=\CG(\chi_\K, \psi_\K, y_\K): = \int_{\CO_\K^\times} \chi_\K(x)^{-1} \cdot \psi_\K(y_\K x)\od\! x,
\ee
where $\od\! x$ is the normalized Haar measure so that $\CO_\K^\times$ has total volume 1. Note that $\CG(\chi_\K)=1$ when
$\frak{c}(\chi_\K)=\frak{c}(\psi_\K)=\CO_\K$. For every $\sigma\in \Aut(\C)$, it is easily checked  that  $\frak{c}({}^\sigma\chi_\K)=\frak{c}(\chi_\K)$, and 
\be \label{gauss22}
\CG(\chi_\K, \psi_\K, y_\K)={}^\sigma \chi_\K(t_{\sigma, \K})\cdot \CG({}^\sigma\chi_\K, \psi_\K, y_\K).
\ee

\subsection{Non-archimdean period relation}

Suppose that $n\geq 2$, and $\Sigma_\K$ is a generic irreducible smooth representation of $\GL_{n-1}(\K)$ with a fixed   Whittaker functional  
\[
 \lambda'_\K \in \Hom_{\mathrm{N}_{n-1}(\K)}(\Sigma_\K, \psi_{n-1, \K})\setminus \{0\}.
\]
As before, we use $\lambda_\K'$ to realize $\Sigma_\K$ as a subrepresentation of $\Ind^{\GL_{n-1}(\K)}_{{\rm N}_{n-1}(\K)}\psi_{n-1,\K}$, and we have a subrepresentation ${}^\sigma \Sigma_\K \subset \Ind^{\GL_{n-1}(\K)}_{{\rm N}_{n-1}(\K)}\psi_{n-1,\K}$ for every $\sigma\in \Aut(\C)$.

As in the archimedean case, denote by $\frak{M}_{n-1,\K}$ the one-dimensional space of invariant measures on $\GL_{n-1}(\K)$. Fix the Haar measure on $\mathrm{N}_{n-1}(\K)$ to be the product of self-dual Haar measures on $\K$ with respect to $\psi_\K$, as in \eqref{mea1}. Then each $m\in \frak{M}_{n-1,\K}$ induces a quotient measure $\bar{m}$  on $\mathrm{N}_{n-1}(\K)\backslash \GL_{n-1}(\K)$.


Let $\chi_{\Sigma_\K}$ denote the  central character of $\Sigma_\K$. For  every $\sigma\in \Aut(\C)$, it is clear that  ${}^\sigma (\chi_{\Sigma_\K}) = \chi_{{}^\sigma\Sigma_\K}$. Similar to \eqref{gauss22}, we also have that
\be \label{gauss222}
\CG(\chi_{\Sigma_\K}, \psi_\K, y'_\K)={}^\sigma \chi_\K(t_{\sigma, \K})\cdot \CG(\chi_{{}^\sigma\Sigma_\K}, \psi_\K, y'_\K),
\ee
where $y'_\K\in \K^\times$ satisfies that $\frak{c}(\psi_\K)= y'_\K \cdot \frak{c}(\chi_{\Sigma_\K})$.

We call an invariant measure $m$ on $\GL_{n-1}(\K)$ rational if $m(K)\in \BQ$ for  every open compact subgroup $K$ of $\GL_{n-1}(\K)$. All rational measures on $\GL_{n-1}(\K)$ form a rational structure of 
$\frak{M}_{n-1,\K}$. By using this rational structure, we get a $\sigma$-linear isomorphism $\sigma: \frak{M}_{n-1,\K}\rightarrow \frak{M}_{n-1,\K}$. By taking the tensor product of the above  $\sigma$-linear isomorphism with the $\sigma$-linear isomorphisms as defined in \eqref{sigmaf}, we get a $\sigma$-linear isomorphism
\be\label{sigmal00}
\sigma: \Pi_\K \otimes \Sigma_\K \otimes \chi_{\K, s-\frac{1}{2}}\otimes \frak{M}_{n-1,\K}\rightarrow  {}^\sigma \Pi_\K \otimes {}^\sigma\Sigma_\K \otimes {}^\sigma(\chi_{\K, s-\frac{1}{2}})\otimes \frak{M}_{n-1,\K},
\ee
where $s\in \C$.

Similar to \eqref{nrs}, we have the normalized Rankin-Selberg integrals 
\[
\oZ^\circ(\cdot, s, \chi_\K)\in \Hom_{\GL_{n-1}(\K)}(\Pi_\K \otimes \Sigma_\K \otimes \chi_{\K, s-\frac{1}{2}} \otimes \frak{M}_{n-1,\K}, \BC)
\]
and
\[
\oZ^\circ(\cdot, s, {}^\sigma \chi_\K)\in \Hom_{\GL_{n-1}(\K)}( {}^\sigma \Pi_\K \otimes  {}^\sigma \Sigma_\K \otimes  ({}^\sigma \chi_{\K})_{s-\frac{1}{2}} \otimes \frak{M}_{n-1,\K}, \BC),
\]
where $ ({}^\sigma \chi_{\K})_{s-\frac{1}{2}} :={}^\sigma \chi_{\K}\cdot \abs{\,\cdot\,}_\K^{s-\frac{1}{2}}$, which equals ${}^\sigma (\chi_{\K, s-\frac{1}{2}})$ when $s\in \frac{1}{2}+\BZ$.

We reformulate the non-archimedean period relation,  which is essentially due to Harder \cite[Section III]{H} for $n=2$ and Mahnkopf \cite[Section 3.4]{Mah} and Raghuram \cite[Section 3.3]{Rag1} in general, as in  the following proposition. 

\begin{prp} \label{nonarchrel}
For all $s_0\in \frac{1}{2}+\BZ$ and $\sigma\in \Aut(\BC)$, the diagram  
\[
 \begin{CD}
          \Pi_\K \otimes \Sigma_\K \otimes \chi_{\K, s_0-\frac{1}{2}}\otimes \frak{M}_{n-1,\K}
                  @> \CG(\chi_{\Sigma_\K}, \psi_\K, y'_\K)\cdot \CG(\chi_\K, \psi_\K, y_\K)^{\frac{n(n-1)}{2}}\cdot \oZ^\circ(\cdot, s_0, \chi_\K) >> \BC \\
            @V  \sigma  VV            @VV  \sigma  V\\
          {}^\sigma \Pi_\K \otimes {}^\sigma\Sigma_\K \otimes ({}^\sigma\chi_{\K})_{ s_0-\frac{1}{2}}\otimes \frak{M}_{n-1,\K} @> \CG(\chi_{{}^\sigma\Sigma_\K},\psi_\K, y'_\K)\cdot \CG({}^\sigma\chi_\K,\psi_\K, y_\K)^{\frac{n(n-1)}{2}}\cdot \oZ^\circ(\cdot, s_0, {}^\sigma\chi_\K) >> \C \\
  \end{CD}
\]
commutes.   
\end{prp}


\begin{proof}
We give a proof for completeness. Note that 
\[
\oL(s, \Pi_\K \times \Sigma_\K \times \chi_\K) = P(q^{\frac{1}{2}-s})^{-1}
\]
for a polynomial $P(X)\in \BC[X]$. 
 For  $\sigma\in \Aut(\BC)$, denote by ${}^\sigma P(X)\in \BC[X]$ the polynomial obtained by applying 
$\sigma$ to the coefficients of the polynomial $P(X)$. Following the proof of \cite[Lemma 4.6]{Clo}, and by noting that the local Rankin-Selberg $\oL$-function does not depend on $\psi_\K$, it is easy to show that
\[
\oL(s, {}^\sigma \Pi_\K \times {}^\sigma\Sigma_\K \times {}^\sigma\chi_\K)= {}^\sigma P(q^{\frac{1}{2}-s})^{-1}.
\]
Specifying $s$ to  $s_0\in \frac{1}{2}+\BZ$, we obtain that
\be \label{sigmaL}
\oL(s_0, {}^\sigma\Pi_\K \times {}^\sigma\Sigma_\K \times {}^\sigma\chi_\K) = \sigma(\oL(s_0, \Pi_\K\times \Sigma_\K \times \chi_\K)).
\ee 

For $f\in \Pi_\K$, $f'\in  \Sigma_\K$, $m\in \frak{M}_{n-1,\K}$, and $s_0\in \frac{1}{2}+\BZ$ large enough, by \eqref{sigmaf}  and \eqref{sigmaL} we have that
\[
\begin{aligned}
& \oZ^\circ({}^\sigma\!f \otimes {}^\sigma\!f' \otimes {}^\sigma m, s_0, {}^\sigma\chi_\K) )\\
= \ & \frac{1}{\oL(s_0, {}^\sigma\Pi_\K\times {}^\sigma \Sigma_\K \times {}^\sigma\chi_\K)}  \\
  & \cdot \int_{\RN_{n-1}(\K)\backslash\GL_{n-1}(\K)}  {}^\sigma\!f \left( \begin{bmatrix} g & 0 \\ 0 & 1 \end{bmatrix}\right) \cdot
 {}^\sigma\!f' (g)   \cdot  {}^\sigma\chi_\K (\det g) \cdot \abs{\det g}_{\K}^{s_0-\frac{1}{2}} \od\! \overline{{}^\sigma m}(g)\\
=\ &  \frac{1}{\sigma(\oL(s_0, \Pi_\K\times \Sigma_\K \times \chi_\K))}   \int_{\RN_{n-1}(\K)\backslash \GL_{n-1}(\K)}  \sigma\left(f \left( \mathbf{t}_{n, \sigma, \K} \begin{bmatrix} g & 0 \\ 0 & 1\end{bmatrix}\right)\right) \\
& \cdot
\sigma\left( f'  (\mathbf{t}_{n-1, \sigma, \K} \cdot g) \right)   \cdot {}^\sigma\chi_\K (\det g) \cdot \abs{\det g}_{\K}^{s_0-\frac{1}{2}} \od\!\overline{{}^\sigma m}(g)\\
=\ & \sigma\Bigg( \frac{1}{\oL(s_0, \Pi_\K\times \Sigma_\K \times \chi_\K)}  \cdot \int_{\RN_{n-1}(\K)\backslash \GL_{n-1}(\K)}  f \left( \begin{bmatrix} t_{\sigma,\K}^{-1}\mathbf{t}_{n-1, \sigma, \K}\, g & 0 \\ 0 & 1\end{bmatrix}\right) \\
& \cdot
 f'  (\mathbf{t}_{n-1, \sigma, \K} \cdot g)   \chi_\K (\det g) \cdot \abs{\det g}_{\K}^{s_0-\frac{1}{2}} \od\!\overline{ m}(g)\Bigg)\\
 = \ & {}^\sigma \chi_{\Sigma_\K}(t_{\sigma, \K})\cdot {}^\sigma \chi_\K(t_{\sigma, \K})^{\frac{n(n-1)}{2}}\cdot \sigma( \oZ^\circ(f \otimes f' \otimes m, s_0,\chi_\K)).
\end{aligned}
\]
It is well-known from the theory of Rankin-Selberg L-functions that the map $s\mapsto \oZ^\circ(f\otimes f' \otimes m, s, \chi_\K)$ is an element of the ring
$\BC[q^{s-\frac{1}{2}}, q^{\frac{1}{2}-s}]$.  Therefore the above equality holds for all $s_0\in \frac{1}{2}+\BZ$. Hence by \eqref{gauss22} and \eqref{gauss222}, the diagram in the proposition is commutative. 
\end{proof}

\section{Whittaker periods} \label{sec4}

Let $\rk$ be a number field with adele ring $\BA$ as in the Introduction. In this section we define the Whittaker periods for irreducible automorphic representations of $\GL_n(\BA)$ that are {\it tamely isobaric} (see \eqref{uni-iso00}) and regular algebraic. 

\subsection{Canonical generators of the cohomology spaces}  \label{sec4.1}
Put
$
K_{n,\infty} :=\prod_{v|\infty} K_{n, \rk_v}
$ ($n\geq 1$), 
where $K_{n, \rk_v}$ is the standard maximal compact subgroup of $\GL_n(\rk_v)$ as in \eqref{cpt-inf} for an archimedean place $v$ of $\rk$.
Define a one-dimensional real vector space
\[
\omega_{n,\infty}(\BR):=\wedge^{d_{n,\infty}} (\frak{gl}_{n}(\rk_\infty) / \frak{k}_{n,\infty}),
\]
where
\[
d_{n,\infty}:= \sum_{v|\infty} d_{n, \rk_v} = \dim_\BR (\frak{gl}_{n}(\rk_\infty) / \frak{k}_{n,\infty}). 
\]
Put
\[
\omega_{n,\infty}:=\omega_{n,\infty}(\BR)\otimes_\BR \BC.
\] 
Similar to  \eqref{orient} in the archimedean case, denote by $\frak{O}_{n,\infty}$ the complex orientation space of $\omega_{n,\infty}$, and put
 \be\label{wtildefo}
 \widetilde{\frak{O}}_{n,\infty}:= \frak{O}_{n-1,\infty}\otimes \cdots\otimes \frak{O}_{1,\infty}\otimes \frak{O}_{0,\infty}.
 \ee
By convention we  set $ \widetilde{\frak{O}}_{0,\infty}:=\frak{O}_{0,\infty}:=\BC$. 
We   identify $\frak{O}_{n,\infty}\otimes \frak{O}_{n,\infty}$ with $\BC$ in the obvious way. Then we have that
\be \label{tildeO}
 \widetilde{\frak{O}}_{n,\infty} \otimes  \widetilde{\frak{O}}_{n-1,\infty} = \frak{O}_{n-1,\infty}.
\ee

Let $\mu=\{\mu^\iota\}_{\iota\in \CE_\rk} \in (\BZ^n)^{\CE_\rk}$ be a highest weight that is pure as in the Introduction. For every archimedean place $v$ of $\rk$,  view $\CE_{\rk_v}$ as a subset of  $\CE_{\rk}$ in the obvious way, and set 
\be \label{muv}
\mu_v:=\{\mu^{\iota}\}_{\iota\in \CE_{\rk_v}} \in  (\BZ^n)^{\CE_{\rk_v}}.
\ee 
  Put
\[
\Omega(\mu) := \left\{ \widehat \otimes_{v|\infty} \pi_{\mu_v}: \pi_{\mu_v}\in \Omega(\mu_v)\right\}
\]
and
\[
 \CH_{\mu} :=\bigoplus_{\pi_{\mu}\in \Omega(\mu)}  \CH(\pi_\mu),
 \]
where
\[
\CH(\pi_\mu):= \RH_\mathrm{ct}^{b_{n,\infty}}(\BR^\times_+\backslash \GL_n(\rk_\infty)^0; F_\mu^\vee\otimes \pi_{\mu})\otimes \widetilde{\frak{O}}_{n,\infty}.
\]
Here $b_{n,\infty} = \sum_{v|\infty} b_{n, \rk_v}$ is as in the Introduction,
and  $\BR^\times_+$ is identified with a central  subgroup of $\GL_n(\rk_\infty)$ via the diagonal embedding.

Recall from the Introduction that $F_\mu$ is an irreducible algebraic representation of $\GL_n(\rk\otimes_\mathbb Q \C)$ of highest weight $\mu$.  It has a decomposition
\[
  F_\mu=\otimes_{v|\infty} F_{\mu_v}.
\]
For every archimedean place $v$ of $\rk$, we have fixed a generator $v_{\mu_v}\in (F_{\mu_v})^{\mathrm{N}_n(\rk_v\otimes_\BR \BC)}$. This yields a generator  
\[
v_{\mu}:=\otimes v_{\mu_v}\in (F_{\mu})^{\mathrm{N}_n(\rk\otimes_\mathbb Q \BC)}.
\]
We remark that the representation $F_\mu$  is unique up to isomorphism, and the pair $(F_\mu,v_\mu)$ is more rigid in the sense that it is unique up to a unique isomorphism. 
Also recall from \eqref{genw} the Whittaker functional $\lambda_{\mu_v}$ on $\pi_{\mu_v}\in \Omega(\mu_v)$. By tensor product, this induces the Whittaker functional $\lambda_\mu$ on every $\pi_\mu\in \Omega(\mu)$. 

Let $\varepsilon\in \widehat{\pi_0(\rk_\infty^\times)}$. Denote the $\varepsilon$-isotypic component of $\CH_\mu$ by $\CH_\mu[\varepsilon]$ (similar notation will be used without further explanation). Then
$\CH_\mu[\varepsilon]$  is one-dimensional. In what follows, we will define a canonical generator 
$\kappa_{\mu, \varepsilon}$ of $\CH_\mu[\varepsilon]$, which is determined by the pairs $(F_\mu,v_\mu)$ and $(\pi_\mu, \lambda_\mu)$. Here we suppose that $\pi_\mu$ is the unique representation in $\Omega(\mu)$ such that $\mathcal H(\pi_\mu)[\varepsilon]\neq \{0\}$. 

We first consider the case that $\mu=0_{n,\infty}$, the zero weight. For $n=1$, we  naturally  identify $\CH_{0_{1,\infty}}[\varepsilon]$ with $\BC$, and put $\kappa_{0_{1,\infty},\varepsilon}:=1$ under this identification.  

For $n\geq 2$, fix 
\[
\pi_{0_{n,\infty}} = \widehat \otimes_{v|\infty} \pi_{0_{n,\rk_v}}\in \Omega(0_{n,\infty})\quad {\rm and}\quad \pi_{0_{n-1,\infty}}= \widehat \otimes_{v|\infty} \pi_{0_{n-1, \rk_v}}\in \Omega(0_{n-1,\infty}),
\]
with fixed Whittaker functionals   
\[
\lambda_{0_{n, \infty}}\in  \Hom_{ \mathrm{N}_n(\rk_\infty)}(\pi_{0_{n,\infty}}, \psi_{n,\infty})\setminus\{0\} 
\]
and 
\[
\lambda_{0_{n-1, \infty}}\in  \Hom_{ \mathrm{N}_{n-1}(\rk_\infty)}(\pi_{0_{n-1,\infty}}, \psi_{n-1,\infty})\setminus\{0\}.  
\]

Denote by $\frak{M}_{n-1, \infty}$ the one-dimensional space of invariant measures on $\GL_{n-1}(\rk_\infty)$. Similar to  the local  case at each archimedean place, we have an identification  
\[
\frak{M}_{n-1,\infty} = \omega_{n-1, \infty}^* \otimes \frak{O}_{n-1,\infty}
\]
by push-forward of measures. Similar to \eqref{nrs}, we have the normalized Rankin-Selberg integral
\[
\oZ^\circ(\cdot, s )  \in \Hom_{\GL_{n-1}(\rk_\infty)} (\pi_{0_{n,\infty}}\widehat{\otimes}\pi_{0_{n-1,\infty}}\otimes  \abs{\det}_{\rk_\infty}^{s-\frac{1}{2}}\otimes \frak{M}_{n-1,\infty}, \BC). 
\]
In view of \eqref{tildeO}, we define a 
map $\CP_{\infty, 0}$ to be the composition of 
\begin{eqnarray} 
   \nonumber &  \CP_{\infty, 0}: & \CH(\pi_{0_{n,\infty}}) \otimes \CH(\pi_{0_{n-1,\infty}}) \\
 & \rightarrow  & \RH_\mathrm{ct}^{d_{n-1,\infty}}( \GL_{n-1}(\rk_\infty)^0; \pi_{0_{n,\infty}} \widehat\otimes \pi_{0_{n-1,\infty}}) \otimes 
   \frak{O}_{n-1,\infty} \nonumber \\
& \rightarrow  & \RH_\mathrm{ct}^{d_{n-1,\infty}}(\GL_{n-1}(\rk_\infty)^0;  \frak{M}_{n-1,\infty}^*)\otimes \frak{O}_{n-1,\infty}=\BC, \nonumber 
\end{eqnarray}
where  the first arrow is the restriction of cohomology composed with the cup product, and the last arrow is the map induced by the linear functional
\[
\oZ^\circ(\cdot, \frac{1}{2}) : \pi_{0_{n,\infty}} \widehat\otimes \pi_{0_{n-1,\infty}} \rightarrow \frak{M}_{n-1,\infty}^*.
\]
By the nonvanishing hypothesis that is proved in \cite{Sun}, these modular symbols for all $ \pi_{0_{n,\infty}}\in \Omega(0_{n,\infty})$ and $\pi_{0_{n-1,\infty}}\in \Omega(0_{n-1,\infty})$  give the following non-degenerate pairings for all $\varepsilon\in
\widehat{\pi_0(\rk_\infty^\times)}$, still denoted by 
\[
\CP_{\infty, 0}: \CH_{0_{n,\infty}}[\varepsilon]\times \CH_{0_{n-1,\infty}}[\varepsilon] \to \BC.
\]
We inductively define $\kappa_{0_{n,\infty}, \varepsilon}$ by requiring that 
\[
\CP_{\infty, 0}(\kappa_{0_{n,\infty}, \varepsilon}, \kappa_{0_{n-1,\infty}, \varepsilon})=1.
\]
In general, we define 
\[
\kappa_{\mu, \varepsilon} := \jmath_\mu (\kappa_{0_{n,\infty},\varepsilon}),
\]
where 
\[
\jmath_\mu: \CH_{0_{n,\infty}}[\varepsilon] \to \CH_\mu[\varepsilon]
\]
is the isomorphism induced by the local ones in Proposition \ref{prpjmath}.

\subsection{Some actions of $\Aut(\C)$} \label{sec4.3} 

Recall the additive character  $\psi_\BR$ from \eqref{psir}. Denote by $\BA_\BQ$ the adele ring of $\BQ$. Fix a nontrivial additive character of $\rk\backslash \BA$ as the composition of 
\be\label{dpsi}
\psi: \rk\backslash \BA \xrightarrow{{\rm Tr}_{\rk/\BQ}} \BQ\backslash \BA_\BQ  \to  \BQ\backslash \BA_\BQ / \widehat\BZ  =  \BR/\BZ \xrightarrow{\psi_\BR} \BC^\times,
\ee
where ${\rm Tr}_{\rk/\BQ}$ is the trace map, and $\widehat\BZ$ is the profinite completion of $\BZ$.
Write $\psi = \otimes_v\psi_v$, where $\psi_v$ is a character of  $\rk_v$ for each place $v$ of $\rk$.  By using $\psi$, we define the character $\psi_{n}$ of $\mathrm{N}_n(\A)$ as in 
\eqref{psi_n} ($n\geq 1$). 
Then we have a decomposition $\psi_n = \psi_{n, f}\otimes \psi_{n,\infty}$, where $\psi_{n, f}$ and $\psi_{n,\infty}$ are characters of ${\rm N}_n(\BA_f)$ and ${\rm N}_n(\rk_\infty)$ respectively. 
Here  $\BA_f$ denotes the finite adele ring of  $\rk$ so that $\BA=\BA_f\times \rk_\infty$.

For every $\sigma\in \Aut(\C)$, put
      \[
      \mathbf{t}_{n, \sigma}  := (\mathbf t_{n,\sigma, \rk_v})_{v\nmid \infty}\in \GL_n(\BA_f) \qquad (\textrm{see \eqref{tnk}}),
      \]  
 and define an action of 
      $\Aut(\BC)$ on $\Ind^{\GL_n(\BA_f)}_{{\rm N}_n(\BA_f)}\psi_{n,f}$  (the smooth induction)  by 
\be \label{sigmaf0}
{}^\sigma\!f(g):=\sigma\left( f(\mathbf{t}_{ n, \sigma} \cdot g)\right),
\ee
where $f\in \Ind^{\GL_n(\BA_f)}_{{\rm N}_n(\BA_f)} \psi_{n,f}$ and $g\in \GL_n(\A_f)$.

      Let $\Pi_f$ be a generic irreducible smooth representation of $\GL_n(\BA_f)$, with a fixed   Whittaker functional 
      \[
      \lambda_f\in \Hom_{{\rm N}_n(\BA_f)}(\Pi_f, \psi_{n,f})\setminus\{0\}.
      \]
As before, we use $\lambda_f$ to realize $\Pi_f$ as a subrepresentation of $\Ind^{\GL_n(\BA_f)}_{{\rm N}_n(\BA_f)}\psi_{n,f}$. 
       The rationality field of $\Pi_f$, denoted by $\BQ(\Pi_f)$, is the  fixed field of the group of field automorphisms $\sigma\in {\rm Aut}(\BC)$ such that ${}^\sigma \Pi_f :=\sigma(\Pi_f)= \Pi_f$. 

Let $\Pi$ be an irreducible smooth automorphic representation of $\GL_n(\BA)$. For the notion of  smooth automorphic representations and the locally convex topologies on them, see \cite{LS} and \cite{GZ} for example. 
   If $\Pi$ is cuspidal, then the {\it exponent} of $\Pi$ is defined to be the real number ${\rm ex}(\Pi)$ such that $\Pi\otimes \abs{\det}_\BA^{-{\rm ex}(\Pi)}$ is unitarizable, where $\abs{\, \cdot \, }_\BA$ is the normalized absolute value 
on $\BA$.  

We say that 
    $\Pi$ is  {\it tamely isobaric} if 
     \be \label{uni-iso00}
      \Pi \cong \Ind^{\GL_n(\BA)}_{P(\BA)} (\Pi_1\otimes\cdots \otimes \Pi_r)\quad(\textrm{normalized smooth induction}),
      \ee
      for   a standard  parabolic subgroup $P$ of $\GL_n$ with Levi subgroup $M_P= \GL_{n_1}\times\cdots \times \GL_{n_r}$, and irreducible cuspidal smooth automorphic representations $\Pi_i$ of $\GL_{n_i}(\BA)$, $i=1,\ldots, r$, that have the same exponent. 
 When $\Pi$ is tamely isobaric, we view it as a space of smooth automorphic forms by using the Eisenstein series (see \cite[Proposition 2]{L} and \cite[Section 1.4.3]{GL} for more details).

 Suppose that $\Pi_f$ and $\Pi_\infty$ are respectively  the finite and infinite part of 
$\Pi$ so that $\Pi= \Pi_f \otimes \Pi_\infty$.   Now we  assume that $\Pi$ is  tamely isobaric and regular algebraic.    By the the proof of \cite[Lemma 1.2]{G} (see also \cite[Section 1.4.3]{GL}),  for every $\sigma\in \Aut(\BC)$,  ${}^\sigma\Pi_f :=\sigma(\Pi_f)$ given by \eqref{sigmaf0} is the finite part of a unique irreducible smooth automorphic representation ${}^\sigma\Pi$ of $\GL_n(\BA)$. Moreover, ${}^\sigma\Pi$ is  also tamely isobaric and regular algebraic. 
   
 \begin{remarkp}
 More precisely, the above assertion holds when $\Pi$ is cuspidal and regular algebraic by \cite[Theorem 3.13]{Clo}. In general, if $\Pi$ is tamely isobaric as in \eqref{uni-iso00} and is regular algebraic, then
 \[
\Xi_1\widehat \otimes\cdots\widehat \otimes \Xi_r\\
 :=(\Pi_1\widehat\otimes\cdots\widehat \otimes\Pi_r)\otimes\rho_P 
 \]
  is a regular algebraic irreducible cuspidal  smooth automorphic representation of $M_P(\BA)$, where $\rho_P :=\delta_P^{\frac{1}{2}}$ is the square root of the modular character $\delta_P$ of $P(\A)$. Then we have that 
 \[
 {}^\sigma\Pi \cong \Ind^{\GL_n(\BA)}_{P(\BA)} \left({}^\sigma \Xi_1\widehat \otimes\cdots \widehat \otimes {}^\sigma\Xi_r\right)\otimes \rho_P^{-1}.
 \]
 \end{remarkp}

 Recall that the rationality field  of $\Pi$ is defined to be $\BQ(\Pi):=\BQ(\Pi_f)$.  Let $\Aut(\BC/\BQ(\Pi))$  act on $\Pi_f$ by  \eqref{sigmaf0}. It is known that $\BQ(\Pi)$  is a number field and  $(\Pi_f)^{\Aut(\BC/\BQ(\Pi))}$ is a $\BQ(\Pi)$-rational structure of $\Pi_f$ (see
 \cite[Lemma 3.2]{RS2}).

As in the Introduction, suppose that $F_\mu$ is an irreducible algebraic representation of $\GL_n(\rk\otimes_{\mathbb Q}\C)$ whose highest weigh  $\mu=\{\mu^\iota\}_{\iota\in \CE_\rk}\in (\BZ^n)^{\CE_\rk}$ is pure of  weight $w_\mu \in \BZ$ so that
\[
\mu^\iota_1+\mu^{\bar{\iota}}_n = \cdots = \mu^\iota_n + \mu^{\bar{\iota}}_1= w_\mu\quad \textrm{for all }\iota\in \CE_\rk.
\]
 Similar to the local case \eqref{algind1}, we realize $F_\mu$ as the algebraic induction 
 \[
 F_\mu = {}^{\rm alg}\Ind^{\GL_n(\rk \otimes_\BQ \BC)}_{\bar{\RB}_n(\rk\otimes_\BQ\BC)}\chi_\mu,
 \]
 and realize $v_\mu\in F_\mu$ as the $\mathrm N_n(\rk\otimes_{\mathbb Q} \C) $-invariant function that has value $1$ at the identity matrix.  Then the generator $v_\mu^\vee\in (F_\mu^\vee)^{\bar{\mathrm N}_n(\rk\otimes_\mathbb Q \C)}$ is identified with the evaluation map at the identity matrix.  Similarly,  $F_\mu^\vee$ is realized  as the algebraic induction 
 \[
 F_\mu^\vee = {}^{\rm alg}\Ind^{\GL_n(\rk \otimes_\BQ \BC)}_{{\RB}_n(\rk\otimes_\BQ\BC)}\chi_{-\mu}.
 \]

For every $\sigma\in \Aut(\BC)$,  write
     \[
     {}^\sigma\!\mu = \{ \mu^{\sigma^{-1} \circ \iota}\}_{\iota \in \CE_\rk}.
     \]
    As a consequence of the purity lemma \cite[Lemma 4.9]{Clo}, $\mu$ necessarily satisfies the condition that (see \cite[Lemma 1.3]{G})
  \[
  \mu^{\sigma\circ\bar{\iota}} = \mu^{\overline{\sigma\circ\iota}}\quad \textrm{ for all $\sigma\in \Aut(\BC)$ and $\iota\in \CE_\rk$}. 
  \]
Therefore, ${}^\sigma\!\mu$ is also pure of weight  $ w_\mu$. 
 The rationality field $\BQ(F_\mu)$ is defined to be the fixed field of the group of field automorphisms $\sigma\in \Aut(\BC)$ such that
$ {}^\sigma\!\mu  =  \mu$. 

Let  $\Aut(\BC)$ act on the space  of algebraic functions on $\GL_n(\rk \otimes_\BQ\BC)$ by
\be \label{sigmaF}
({}^\sigma\!f)(x):= \sigma( f(\sigma^{-1}x))\quad (x\in \GL_n(\rk\otimes_\BQ\BC)),
\ee
where $\Aut(\BC)$ acts on $\GL_n(\rk \otimes_\BQ\BC)$ through its action on the second factor of $\rk\otimes_\BQ \BC$. 
Then
\[
\sigma(F_\mu)=  F_{ {}^\sigma\!\mu} \quad \textrm{and} \quad \sigma(F_\mu^\vee)=  F_{ {}^\sigma\!\mu}^\vee.
\]



 Define 
  \be \label{CXn}
  \CX_{n}:=(\BR^\times_+\cdot \GL_{n}(\rk))\backslash \GL_{n}(\BA)/ K_{n,\infty}^0.
    \ee
    For every open compact subgroup $K_f$ of $\GL_n(\BA_f)$, the finite-dimensional representation $F_\mu^\vee$ defines a sheaf on $\CX_n/K_f$, which is still denoted by $F_\mu^\vee$. 
    Let
$\RH^{b_{n,\infty}}(\CX_n/ K_f, F_\mu^\vee)$ be the sheaf cohomology group, and define
\be\label{sheafcoh}
\begin{aligned}
\CH^{b_{n,\infty}}(\CX_n, F_\mu^\vee)& := \RH^{b_{n,\infty}}(\CX_n, F_\mu^\vee) \otimes  \widetilde{\frak{O}}_{n,\infty} \\
&:= \varinjlim_{K_f}\RH^{b_{n,\infty}}(\CX_n/ K_f, F_\mu^\vee)\otimes \widetilde{\frak{O}}_{n,\infty},
\end{aligned}
\ee
where $K_f$ runs over the directed system of open compact subgroups of $\GL_n(\BA_f)$. 

Note that $\widetilde{\frak{O}}_{n,\infty}$ has a natural $\BQ$-structure.   For every $\sigma\in \Aut(\BC)$,  the  map \eqref{sigmaF} induces a $\sigma$-linear map
   \be \label{sigF}
   \sigma: \CH^{b_{n,\infty}}(\CX_n, F_\mu^\vee) \to \CH^{b_{n,\infty}}(\CX_n, F_{{}^\sigma\!\mu}^\vee).
   \ee


Put
    \[
    \GL_n(\BA)^\natural:= \GL_n(\BA_f)\times \pi_0(\rk_\infty^\times).
    \]
    Then both the domain and codomain of the map \eqref{sigF} are naturally smooth representations of $\GL_n(\BA)^\natural$, and the map  \eqref{sigF} is  $\GL_n(\BA)^\natural$-equivariant.

 \subsection{Definition of the Whittaker periods}  \label{sec4.4} 

 Fix the Haar measure on $\mathrm{N}_n(\BA)$ to be the product of self-dual Haar measures on $\BA$ with respect to $\psi$, as in \eqref{mea1}.
Then we have a nonzero continuous linear functional 
\be \label{can-whit}
\lambda \in \Hom_{\mathrm{N}_n(\BA)}(\Pi, \psi_n),\quad \varphi\mapsto \int_{\mathrm{N}_n(\rk)\backslash\mathrm{N}_n(\BA)} \varphi(u) \cdot \overline{\psi_n}(u) \od\! u.
\ee
By the uniqueness of 
Whittaker models, we have  a factorization
\be\label{fact}
\lambda  =  \lambda_f \otimes  \lambda_\infty,
\ee
where 
\[
\lambda_f\in \Hom_{\mathrm{N}_n(\BA_f)}(\Pi_f, \psi_{n,f})
\]
as before, and 
\[
\lambda_\infty\in \Hom_{\mathrm{N}_n(\rk_\infty)}(\Pi_\infty, \psi_{n,\infty}).
\]

More generally, for every $\sigma\in \Aut(\C)$, let ${}^\sigma \lambda \in \Hom_{\mathrm{N}_n(\BA)}({}^\sigma\Pi, \psi_n)$ be the Whittaker functional defined by the integrals as in \eqref{can-whit}. Similar to \eqref{fact},  we also have  factorizations 
\[
{}^\sigma \Pi= {}^\sigma\Pi_f \otimes {}^\sigma \Pi_\infty\quad\textrm{ and }\quad  {}^\sigma\lambda  = {}^\sigma \lambda_f \otimes  {}^\sigma \lambda_\infty.
\]
Recall that  ${}^\sigma\Pi_f :=\sigma(\Pi_f)$ is realized as a space of Whittaker functions so that  ${}^\sigma\lambda_f$ is realized as the evaluation map at the identity matrix. 

Suppose that $F_\mu$ is the coefficient system of  $\Pi$ as in the Introduction. Then $F_{{}^\sigma\!\mu}$ is the coefficient system of ${}^\sigma \Pi$  (\cf \cite[Theoerem 3.13]{Clo} and \cite[Corollary 1.4]{G}) so that
 \[
  \CH({}^\sigma\Pi_\infty) :=\RH_\mathrm{ct}^{b_{n,\infty}}(\BR^\times_+\backslash \GL_n(\rk_\infty)^0; F_{ {}^\sigma\!\mu}^\vee\otimes  {}^\sigma\Pi_\infty)\otimes \widetilde{\frak{O}}_{n,\infty}\neq \{0\}.
  \]
  Consequently, $\BQ(F_\mu)\subset \BQ(\Pi)$.
 Put
 \[
  \CH({}^\sigma\Pi) :=\RH_\mathrm{ct}^{b_{n,\infty}}(\BR^\times_+\backslash \GL_n(\rk_\infty)^0; F_{ {}^\sigma\!\mu}^\vee\otimes  {}^\sigma\Pi)\otimes \widetilde{\frak{O}}_{n,\infty}. 
 \]
 Then we have the canonical isomorphism 
 \[
 \iota_{\rm can}: {}^\sigma\Pi_f\otimes  \CH({}^\sigma\Pi_\infty) \to \CH({}^\sigma\Pi). 
 \]

   Following \cite[Lemma 3.15]{Clo} and \cite[Proposition 1.6]{G}, we have a $\GL_n(\BA)^\natural$-equivariant embedding 
  \be \label{iota}
   \iota_{\Pi}: \Pi_f\otimes \CH(\Pi_\infty) = \CH(\Pi)  \hookrightarrow \CH^{b_{n,\infty}}(\CX_n, F_\mu^\vee).
  \ee
Let $\varepsilon \in \widehat{ \pi_0(\rk_\infty^\times)}$ be the character 
$
\varepsilon_{\Pi_\infty}\cdot\sgn_\infty^{\frac{(n-1)(n-2)}{2}}
$
when  $n$  is odd, and be arbitrary when $n$ is even. Then we have a $\GL_n(\BA)^\natural$-equivariant linear embedding 
  \[
  \iota_{\Pi, \varepsilon}:  \Pi_f\otimes \CH(\Pi_\infty)[\varepsilon] = \CH(\Pi)[\varepsilon]  \hookrightarrow \CH^{b_{n,\infty}}(\CX_n, F_\mu^\vee).
  \]

   \begin{prp} \label{prop4.2}
   Let the notations and assumptions be as above. Then 
   \[
   \dim_{\GL_n(\BA)^\natural}(\CH(\Pi)[\varepsilon], \CH^{b_{n,\infty}}(\CX_n, F_\mu^\vee))=1,
   \]
and for every $\sigma\in \Aut(\C)$ the map \eqref{sigF}  induces a commutative diagram 
   \[  
 \begin{CD}
            \CH(\Pi)[\varepsilon]
                  @> \iota_{\Pi,\varepsilon} >>   \CH^{b_{n,\infty}}(\CX_n, F_\mu^\vee) \\
            @ V  \sigma VV          @  V V \sigma V  \\
          \CH({}^\sigma\Pi)[\varepsilon] @> \iota_{{}^\sigma\Pi, \varepsilon} >>  \CH^{b_{n,\infty}}(\CX_n, F_{{}^\sigma\!\mu}^\vee). \\
  \end{CD}
\]
Moreover,  under the action given by the left vertical arrow of the above diagram, $(\CH(\Pi)[\varepsilon])^{\Aut(\BC/\BQ(\Pi))}$ is a $\BQ(\Pi)$-rational structure of $\CH(\Pi)[\varepsilon]$.
   \end{prp}
   
   \begin{proof} The commutative diagram follows from \cite[Propositions 1.19 and 1.21]{GL}, and the fact that the map \eqref{sigF} commutes with the actions of $\pi_0(\rk_\infty^\times)\cong \pi_0(K_{n, \infty})$.
   The last assertion is implied by  Drinfeld-Manin principle (see \cite[Proposition 3.16]{Clo}) and \cite[Lemma 3.2.1]{Clo}. 
\end{proof}

\begin{remarkp}\label{remks54}
It follow from Proposition \ref{prop4.2} that for every $\sigma\in \Aut(\BC)$, the central character of $F_{{}^\sigma\!\mu}^\vee\otimes ({}^\sigma \Pi)_\infty$ equals that of $F_\mu^\vee\otimes \Pi_\infty$. Consequently, $({}^\sigma \Pi)_\infty$ is uniquely determined by $\sigma$ and $\Pi_\infty$. Specifying to the case that $n=1$,  we know that the infinite part of ${}^\sigma\chi$ equals that of $\chi$, for every finite order Hecke character $\chi: \rk^\times \backslash \A^\times\rightarrow \C^\times$.
\end{remarkp}
 

  We equip $\CH(\Pi)[\varepsilon]$ with the action of $\Aut(\BC/\BQ(\Pi))$ given by Proposition \ref{prop4.2}. Write
   \[
   \Pi^\natural:=\Pi_f\otimes \varepsilon \cong \CH(\Pi)[\varepsilon],
   \]
  where  $\varepsilon\in   \widehat{ \pi_0(\rk_\infty^\times)}$ is identified with $\BC$ as a vector space. We equip $\Pi^\natural$ with the action of  $\Aut(\BC/\BQ(\Pi))$ given by its action on $\Pi_f$ as in  \eqref{sigmaf0} and its natural action on $\C$. 
 
 \begin{lemp} \label{inv-gen}
 There exists a generator 
 \[
 \omega_{\Pi^\natural}\in  \Hom_{\GL_n(\BA)^\natural}(\Pi^\natural, \CH(\Pi)[\varepsilon])
 \]
 that is  $\Aut(\BC/\BQ(\Pi))$-equivariant. Moreover, such a generator is unique up to multiplication by  scalar in $\BQ(\Pi)^\times$.
 \end{lemp}
 
 \begin{proof}
 Recall that $(\Pi^\natural)^{\Aut(\BC/\BQ(\Pi))}$ is a $\BQ(\Pi)$-rational structure of $\Pi^\natural$ 
 (\cite[Lemma 3.2]{RS2}), and $(\CH(\Pi)[\varepsilon])^{\Aut(\BC/\BQ(\Pi))}$ is a $\BQ(\Pi)$-rational structure of $\CH(\Pi)[\varepsilon]$ (Proposition  \ref{prop4.2}).
 Let $\overline{\BQ}$ denote the field of algebraic numbers in $\C$.  By  the multiplicity one property of new vectors,  the $\overline{\BQ}$-rational structure of $\Pi_f$
is unique up to homotheties (see the proof of \cite[Theorem 3.13]{Clo}, and \cite[Chapter I]{W}).  It follows that 
$\Pi^\natural$ and $\CH(\Pi)[\varepsilon]$ are isomorphic over 
$\overline{\BQ}$.  Since  $(\Pi^\natural)^{\Aut(\BC/\overline{\BQ})}$ is irreducible (as a smooth representation of $\GL_n(\BA)^\natural$ over $\overline \BQ$), 
$\Aut(\overline{\BQ}/\BQ(\Pi))$ acts continuously on the one-dimensional $\overline{\BQ}$-vector space (with the discrete topology) 
\[
\Hom_{\GL_n(\BA)^\natural}\left((\Pi^\natural)^{\Aut(\BC/\overline{\BQ})}, (\CH(\Pi)[\varepsilon])^{\Aut(\BC/\overline{\BQ})}\right).
\]
This implies the existence of $\omega_{\Pi^\natural}$ by \cite[Proposition 11.1.6]{Spr}. The uniqueness is obvious. 
 \end{proof}

     Fix $\omega_{\Pi^\natural}$ as in Lemma \ref{inv-gen}. For $\sigma\in \Aut(\BC)$, put 
     \[
 {}^\sigma \Pi^{\natural}:={}^\sigma\Pi_f\otimes \varepsilon. 
 \]
 The $\sigma$-linear isomorphisms  $\sigma: \Pi^{\natural}\to {}^\sigma \Pi^{\natural}$ and 
     \[
     \sigma: \CH(\Pi)[\varepsilon] \to \CH({}^\sigma\Pi)[\varepsilon] \qquad (\textrm{see  Proposition \ref{prop4.2}}) 
     \]
      induce a $\sigma$-linear isomorphism 
      \[
      \sigma:  \Hom_{\GL_n(\BA)^\natural}(\Pi^\natural, \CH(\Pi)[\varepsilon]) \to \Hom_{\GL_n(\BA)^\natural}({}^\sigma \Pi^{\natural}, \CH({}^\sigma\Pi)[\varepsilon]).
      \]
      Using this isomorphism, we define
      \[
      \omega_{{}^\sigma\Pi^{\natural}}:=\sigma(\omega_{\Pi^\natural})\in  \Hom_{\GL_n(\BA)^\natural}({}^\sigma\Pi^{\natural}, \CH({}^\sigma\Pi)[\varepsilon]).
      \]
Unraveling definitions, we have a commutative diagram
\be \label{omegasig}
       \begin{CD}
          \Pi^\natural
                  @> \omega_{\Pi^\natural} >> \CH(\Pi)[\varepsilon] \\
            @V  \sigma VV            @VV  \sigma  V\\
          {}^\sigma\Pi^{\natural} @> \omega_{{}^\sigma\Pi^{\natural}} >> \CH({}^\sigma\Pi)[\varepsilon]. \\
  \end{CD}
  \ee
  
Recall form  Section \ref{sec4.1} that the pairs  $(F_\mu,v_\mu)$ and $(\Pi_\infty, \lambda_\infty)$ determine   a generator $\kappa_{\mu, \varepsilon}$ of  $\CH(\Pi_\infty)[\varepsilon]=\CH_\mu[\varepsilon]$. More generally, for every $\sigma\in \Aut(\C)$, the pairs  $(F_{{}^\sigma\!\mu},v_{{}^\sigma\!\mu})$ and $({}^\sigma\Pi_\infty, {}^\sigma\lambda_\infty)$ determine   a generator $\kappa_{{}^\sigma\!\mu, \varepsilon}$ of  $\CH({}^\sigma\Pi_\infty)[\varepsilon]=\CH_{{}^\sigma\!\mu}[\varepsilon]$.

      \begin{dfnp} \label{df-per}
    For every $\sigma\in \Aut(\C)$,   the Whittaker period $\Omega_\varepsilon({}^\sigma\Pi)\in \BC^\times$ is the unique scalar such that the  diagram 
      \be \label{diag-per}
       \begin{CD}
          {}^\sigma\Pi_f\otimes \varepsilon
                  @>  {\rm id}\otimes  \kappa_{{}^\sigma\!\mu, \varepsilon} >>  {}^\sigma\Pi_f\otimes \CH({}^\sigma\Pi_\infty)[\varepsilon] \\
            @V  \Omega_\varepsilon({}^\sigma\Pi) VV            @ VV \iota_{\rm can} V \\
      {}^\sigma\Pi_f\otimes\varepsilon @> \omega_{{}^\sigma\Pi^\natural} >> \CH({}^\sigma\Pi)[\varepsilon] \\
  \end{CD}
      \ee
      commutes. 
      \end{dfnp}
     
     Up to scalar multiplication by $\BQ(\Pi)^\times$, the Whittaker periods defined above are independent of the choice of the generator  $\omega_{\Pi^\natural} $. More precisely, we have the following lemma. 
     \begin{lemp} \label{uniquefamily}
Let $c\in \BQ(\Pi)^\times$ so that $\omega_{\Pi^\natural}' := c\cdot \omega_{\Pi^\natural}\in  \Hom_{\GL_n(\BA)^\natural}(\Pi^\natural, \CH(\Pi)[\varepsilon])$ is another generator  that is  $\Aut(\BC/\BQ(\Pi))$-equivariant, which defines a corresponding family of Whittaker periods
 $\{\Omega'_\varepsilon({}^\sigma\Pi)\}_{\sigma\in\Aut(\BC)}$. Then for all $\sigma\in\Aut(\BC)$,
 \[
 \sigma\left(\frac{\Omega'_\varepsilon(\Pi)}{\Omega_\varepsilon(\Pi)}\right) = \frac{\Omega'_\varepsilon({}^\sigma\Pi)}{\Omega_\varepsilon({}^\sigma\Pi)}= \sigma(c^{-1}).
 \]
\end{lemp}

\begin{proof}
   This is an easy consequence of the commutative diagrams \eqref{omegasig} and \eqref{diag-per}.
\end{proof}

      For every $\sigma\in \Aut(\BC)$, we define a $\sigma$-linear map 
      \be \label{sig-inf}
      \sigma: \CH(\Pi_\infty)[\varepsilon] \to \CH({}^\sigma \Pi_\infty)[\varepsilon]
      \ee
      such that
      \[
      \sigma(\kappa_{\mu,\varepsilon}) = \kappa_{{}^\sigma\!\mu,\varepsilon}.
      \]
      
      \begin{prp}
      For all $\sigma\in \Aut(\BC)$, the  diagram 
       \be \label{per-diag}
      \begin{CD}
          \Pi_f \otimes \CH(\Pi_\infty)[\varepsilon]
                  @>  \Omega_\varepsilon(\Pi)^{-1} \cdot \iota_{\rm can} >> \CH(\Pi)[\varepsilon] \\
            @V  \sigma VV            @VV  \sigma  V\\
          {}^\sigma\Pi_f \otimes \CH({}^\sigma\Pi_\infty)[\varepsilon]@>  \Omega_\varepsilon({}^\sigma \Pi)^{-1} \cdot \iota_{\rm can} >> \CH({}^\sigma\Pi)[\varepsilon] \\  \end{CD}
   \ee
   commutes, where the left vertical arrow is the $\sigma$-linear map induced by the  map $\sigma: \Pi_f\to {}^\sigma\Pi_f$ and the map \eqref{sig-inf}.
      \end{prp}
      
      \begin{proof}
      This follows easily from \eqref{omegasig} and \eqref{diag-per}. 
      \end{proof}
      

\section{Modular symbols and proof of Theorem \ref{thm: global period relation}} \label{sec5}


\subsection{Rankin-Selberg integrals} \label{sec5.1}

In this subsection, let $\Pi$ be an irreducible cuspidal smooth automorphic representation of $\GL_n(\BA)$ ($n\geq 2$), and 
let $\Sigma$ be a tamely isobaric  irreducible smooth automorphic representation of $\GL_{n-1}(\BA)$. As before, we realize $\Pi$ and $\Sigma$ as  spaces 
of smooth automorphic forms on $\GL_n(\BA)$ and $\GL_{n-1}(\BA)$ respectively.

As in \eqref{can-whit}, we have Whittaker functionals
\[
\lambda \in \Hom_{\mathrm{N}_n(\BA)}(\Pi, \psi_n)\quad \textrm{and}\quad \lambda' \in \Hom_{\mathrm{N}_{n-1}(\BA)}(\Sigma, \psi_{n-1})\]
defined by integrals, and as  
in \eqref{fact}, we have decompositions 
\[
\lambda  =  \lambda_f \otimes  \lambda_\infty\quad \textrm{and}\quad 
\lambda'  =  \lambda'_f \otimes  \lambda'_\infty,
\]
with
\[
\lambda_f\in \Hom_{\mathrm{N}_n(\BA_f)}(\Pi_f, \psi_{n,f}), \quad \lambda_\infty \in \Hom_{\mathrm{N}_n(\rk_\infty)}(\Pi_\infty, \psi_{n,\infty}), \quad \Pi=\Pi_f\otimes \Pi_\infty,
\]
and 
\[
\lambda'_f\in \Hom_{\mathrm{N}_{n-1}(\BA_f)}(\Sigma_f, \psi_{n-1,f}), \quad \lambda'_\infty \in \Hom_{\mathrm{N}_{n-1}(\rk_\infty)}(\Sigma_\infty, \psi_{n-1,\infty}), \quad \Sigma=\Sigma_f\otimes \Sigma_\infty.
\]


Let 
$\chi: \rk^\times\backslash \BA^\times\to \BC^\times$ be a Hecke character.  Similar to \eqref{chikt}, for each $t\in \C$ define a character 
\be\label{dchit}
\chi_t:= \chi \cdot \abs{\,\cdot\,}_\BA^t: \A^\times \rightarrow \C^\times. 
\ee
 As usual, write
\[
\chi=\otimes_{v} \chi_v=\chi_f \otimes \chi_\infty \quad \textrm{and}\quad \chi_t = \chi_{f, t} \otimes \chi_{\infty, t}.
\] 

Denote by $\frak{M}_{n-1}$ and $\frak{M}_{n-1,f}$  the one-dimensional spaces of invariant measures on $\GL_{n-1}(\BA)$ and $\GL_{n-1}(\BA_f)$ respectively, so that
\[
\frak{M}_{n-1} = \frak{M}_{n-1, f}\otimes \frak{M}_{n-1,\infty}.
\] 
Similar to \eqref{nrs}, we have the finite part of the normalized Rankin-Selberg integral 
\[
\oZ^\circ(\cdot, s, \chi_f)\in \Hom_{\GL_{n-1}(\BA_f)}(\Pi_f\otimes \Sigma_f \otimes \chi_{f, s-\frac{1}{2}}\otimes \frak{M}_{n-1,f}, \BC), 
\]
and the normalized Rankin-Selberg integral at infinity 
\[
\oZ^\circ(\cdot, s, \chi_\infty)  \in \Hom_{\GL_{n-1}(\rk_\infty)} (\Pi_\infty\widehat{\otimes}\Sigma_\infty\otimes \chi_{\infty, s-\frac{1}{2}}\otimes \frak{M}_{n-1,\infty}, \BC).
\]

 Define the global Rankin-Selberg period integral 
\[
                \oZ(\cdot, s, \chi) :  \, \Pi \widehat{\otimes} \Sigma\otimes \chi_{s-\frac{1}{2}} \otimes \frak{M}_{n-1}  \rightarrow \BC              
                  \]
                  by
                  \[
                  \oZ(\varphi\otimes \varphi'\otimes 1\otimes m):=\int_{ \GL_{n-1}(\rk)\backslash \GL_{n-1}(\BA)} \varphi\left(\begin{bmatrix} g & 0 \\ 0 & 1\end{bmatrix}\right) \cdot \varphi'(g) \cdot \chi(g)\cdot \abs{\det g}_\BA^{s-\frac{1}{2}} \od\! \bar{m}(g), 
                  \]
                  where $\varphi\in \Pi$, $\varphi'\in \Sigma$,  $m\in \frak{M}_{n-1}$ and $\bar{m}$ is the quotient measure of $m$. Here 
                  \[
                  \Pi\widehat\otimes \Sigma:= (\Pi_\infty\widehat{\otimes}\Sigma_\infty)\otimes (\Pi_f\otimes \Sigma_f).
                  \]
                  
The following proposition reformulates the Euler factorization of Rankin-Selberg period integrals established in  \cite[page 796, equation (7)]{JS2}.

\begin{prp} \label{measure}
For  $\Pi$, $\Sigma$, $\chi$ as above, and all $s\in \BC$, the diagram
   \[
                   \begin{CD}
         (\Pi_\infty\widehat{\otimes} \Sigma_\infty\otimes \chi_{\infty, s-\frac{1}{2}}\otimes \frak{M}_{n-1,\infty})\otimes (\Pi_f\otimes\Sigma_f\otimes \chi_{f,s-\frac{1}{2}} \otimes \frak{M}_{n-1, f})  @> \oZ^\circ(\,\cdot\,, s, \chi_\infty)\otimes \oZ^\circ(\,\cdot\,, s, \chi_f)>> \BC \\
            @V= VV           @V \oL(s, \Pi\times \Sigma \times \chi) VV\\
        \Pi\widehat{\otimes} \Sigma \otimes \chi_{s-\frac{1}{2}}\otimes \frak{M}_{n-1} @>\oZ(\cdot, s, \chi) >> \BC\\
                     \end{CD}
           \]
         commutes. 
           \end{prp}

\subsection{Modular symbols and modular symbols at infinity} \label{sec5.2} From now on,  further assume that $\Pi$ and $\Sigma$  are regular algebraic with balanced coefficient systems $F_\mu$ and $F_\nu$ respectively, and   assume that $\chi$ has finite order.

Similar to \eqref{sheafcoh} we have the following space given by  sheaf cohomology with compact support: 
\be\label{sheafcoh2}
\begin{aligned}
\CH_c^{b_{n,\infty}}(\CX_n, F_\mu^\vee)& := \RH_c^{b_{n,\infty}}(\CX_n, F_\mu^\vee) \otimes  \widetilde{\frak{O}}_{n,\infty} \\
&:= \varinjlim_{K_f}\RH_c^{b_{n,\infty}}(\CX_n/ K_f, F_\mu^\vee)\otimes \widetilde{\frak{O}}_{n,\infty},
\end{aligned}
\ee
where $K_f$ runs over the directed system of open compact subgroups of $\GL_n(\BA_f)$. This is also a smooth representation of $\GL_n(\A)^\natural $ and we have a natural $\GL_n(\A)^\natural $-equivariant linear map
\[
 \iota_\mu: \CH_c^{b_{n,\infty}}(\CX_n, F_\mu^\vee)\rightarrow \CH^{b_{n,\infty}}(\CX_n, F_\mu^\vee).
\]
Since $\Pi$ is cuspidal, there is a natural embedding 
\[
 \iota_\Pi':  \CH(\Pi) \hookrightarrow \CH^{b_{n,\infty}}_c(\CX_n, F_\mu^\vee)
\]
such that $\iota_\mu\circ  \iota_\Pi'= \iota_\Pi$ (see \cite[Lemma 3.15]{Clo}).


Put
\[
\widetilde{\CX}_{n-1}:= \GL_{n-1}(\rk)\backslash \GL_{n-1}(\BA)/ K_{n-1,\infty}^0.
\]
The embedding $\imath: \GL_{n-1}(\BA) \hookrightarrow \GL_n(\BA)$ given by \eqref{emh} induces a proper map, still denoted by 
\[
\imath: \widetilde{\CX}_{n-1} \to \CX_n,
\]
which induces a map 
\[
\imath^*: \CH^{b_{n,\infty}}_c(\CX_n, F_\mu^\vee) \to \CH^{b_{n,\infty}}_c(\widetilde{\CX}_{n-1}, F_\mu^\vee).
\]
The natural map $\wp: \widetilde{\CX}_{n-1}\to \CX_{n-1}$ induces a map
\[
\wp^*: \CH^{b_{n-1,\infty}}(\CX_{n-1}, F_\nu^\vee) \to \CH^{b_{n-1,\infty}}(\widetilde{\CX}_{n-1}, F_\nu^\vee).
\]

Since $\xi:=(\mu, \nu)$ is assumed to be balanced,  by \cite[Theorem 2.21]{Rag2} (see also
\cite[Theorem 2.3]{KS} and \cite[Lemma 4.7]{GH}) we have that
\[
\{j\in \BZ \, :\, j\textrm{ is balanced for }\xi\} =\left\{j\in \BZ \, : \, \frac{1}{2}+j \textrm{ is a critical place of }\Pi\times\Sigma\right\}.
\]
Recall that a half-integer $\frac{1}{2}+j$ is a critical place 
of $\Pi \times \Sigma$ if it  is a pole of neither
$\oL(s, \Pi_\infty\times \Sigma_\infty)$ nor $\oL(1-s, \Pi_\infty^\vee\times \Sigma_\infty^\vee)$.

Let $j$ be a balanced place for $\xi$. 
Define an algebraic  character 
\[
\delta_{ j}:=\otimes_{\iota\in \CE_\rk} {\det}^j 
\]
of $\GL_{n-1}(\rk \otimes_\BQ\BC)$. Put
\[
\oH(\chi_j) : = \RH^0_{\rm ct}(\BR^\times_+\backslash \GL_{n-1}(\rk_\infty)^0;   \delta_{j}^\vee\otimes \chi_{ j})\qquad(\textrm{see \eqref{dchit} for the definition of $\chi_j$}).
\]
Then we have a natural  injective map
\[
\iota_j: \oH(\chi_j) \hookrightarrow \RH^0(\CX_{n-1}, \delta_{ j}^\vee). 
\]

With the notation as before, we have the generators
\[
  v^\vee_\mu := \otimes_{\iota\in \CE_\rk} v_{\mu^\iota}^\vee\in (F_\mu^\vee)^{\bar{\mathrm N}_n(\rk\otimes_\BQ \C)}\quad \textrm{and}\quad v_\nu^\vee := \otimes_{\iota\in \CE_\rk} v_{\nu^\iota}^\vee\in (F_\nu^\vee)^{\bar{\mathrm N}_n(\rk\otimes_\BQ \C)}.
\]
Put $F_\xi := F_\mu\otimes F_\nu$ and
$
v_\xi^\vee : =  v_\mu^\vee \otimes v_\nu^\vee.
$
Recall from \eqref{z} the element \[
z=(z_n, z_{n-1})\in \GL_n(\Z)\times \GL_{n-1}(\Z)\subset \GL_n(\rk)\times \GL_{n-1}(\rk)
\]
By Proposition \ref{phixij}, we have a unique element 
\[
\phi_{\xi, j}\in \Hom_{\GL_{n-1}(\rk \otimes_{\BR}\BC)}(F_\xi^\vee \otimes \delta_{ j}^\vee, \mathbb{C})
\]
such that
$
\phi_{\xi, j} \left(  z.v_\xi^\vee \otimes 1 \right)=1.
$
Then $\phi_{\xi, j}$ induces a linear map
\[
\phi_{\xi, j}: \RH^{d_{n-1,\infty}}_c(\widetilde{\CX}_{n-1}, F_\xi^\vee)\otimes \RH^0(\widetilde{\CX}_{n-1}, \delta_{ j}^\vee) \to \RH^{d_{n-1,\infty}}_c(\widetilde{\CX}_{n-1}, \BC).
\]
Put
\[
\frak{M}_{n-1}^\natural: = \frak{M}_{n-1, f}\otimes \frak{O}_{n-1,\infty}.
\]
Note  that $\widetilde{\CX}_{n-1}/K_f$ is an orientable manifold when   $K_f$ is a sufficiently small open compact subgroup of $\GL_{n-1}(\BA_f)$, and pairing with the fundamental class yields a linear map 
\[
\int_{\widetilde{\CX}_{n-1}}:   \RH^{d_{n-1,\infty}}_c (\widetilde{\CX}_{n-1}, \BC) \otimes \frak{M}_{n-1}^\natural  
\to  \BC.
\]
See \cite[Section 5.1]{Mah} for more explanations.

In view of \eqref{tildeO}, 
define the modular symbol $\CP_{j}$ to be the composition of
\begin{eqnarray}
 \label{pj} & \CP_{j}:&  \CH(\Pi)\otimes \CH(\Sigma) \otimes \oH(\chi_j) \otimes \frak{M}_{n-1, f} \\
   \nonumber  &\xrightarrow{\iota'_\Pi\otimes\iota_\Sigma\otimes\iota_j\otimes {\rm id}} &\CH^{b_{n,\infty}}_c(\CX_n, F_\mu^\vee)\otimes  \CH^{b_{n-1,\infty}}(\CX_{n-1}, F_\nu^\vee) \otimes \RH^0(\CX_{n-1}, \delta_{ j}^\vee)  \otimes \frak{M}_{n-1, f} \\
     \nonumber     &\xrightarrow{\imath^*\otimes \wp^* \otimes \wp^*\otimes {\rm id}} & \CH^{b_{n,\infty}}_c(\widetilde{\CX}_{n-1}, F_\mu^\vee)\otimes  \CH^{b_{n-1,\infty}}(\widetilde{\CX}_{n-1}, F_\nu^\vee) \otimes \RH^0(\widetilde{\CX}_{n-1}, \delta_{ j}^\vee) \otimes \frak{M}_{n-1, f} \\
        \nonumber  &\xrightarrow{\cup\otimes {\rm id}} & \RH^{d_{n-1,\infty}}_c(\widetilde{\CX}_{n-1}, F_\xi^\vee \otimes \delta_{ j}^\vee)\otimes \frak{O}_{n-1,\infty}  \otimes \frak{M}_{n-1, f} \\
          \nonumber  &\xrightarrow{\phi_{\xi, j}\otimes {\rm id}} & \RH^{d_{n-1,\infty}}_c (\widetilde{\CX}_{n-1}, \BC)  \otimes \frak{M}_{n-1, f}^\natural \\
       \nonumber   &\xrightarrow{\int_{\widetilde{\CX}_{n-1}}} & \BC.
      \end{eqnarray}

    Recall that we have the normalized Rankin-Selberg integral at infinity 
      \[
      \begin{aligned}
\oZ^\circ(\cdot, s, \chi_\infty)  & \in \Hom_{\GL_{n-1}(\rk_\infty)} (\Pi_\infty\widehat{\otimes}\Sigma_\infty\otimes \chi_{\infty, s-\frac{1}{2}}\otimes \frak{M}_{n-1,\infty}, \BC) \\ 
& =  \Hom_{\GL_{n-1}(\rk_\infty)} ( \Pi_\infty\widehat{\otimes}\Sigma_\infty\otimes \chi_{\infty, s-\frac{1}{2}}, \frak{M}_{n-1,\infty}^*).
\end{aligned}
\]
Put
\[
\begin{aligned}
\oH(\chi_{\infty, j}) : & = \RH^0_{\rm ct}(\BR^\times_+\backslash \GL_{n-1}(\rk_\infty)^0;   \delta_{j}^\vee\otimes \chi_{\infty, j}) \\
& =  \RH^0_{\rm ct}(\BR^\times_+\backslash \GL_{n-1}(\rk_\infty)^0;   \chi_\infty\cdot\sgn_\infty^j ),
\end{aligned}
\]
where  $\sgn_\infty$ is given as in the Introduction.

Analogous to the archimedean modular symbol  defined in Section \ref{secams}, we define the
{\it modular symbol at infinity}, which is denoted by $\CP_{\infty,j}$,  to be the composition of 
\begin{eqnarray} 
   \nonumber &  \CP_{\infty,j}:  & \CH(\Pi_\infty) \otimes \CH(\Sigma_\infty) \otimes \RH(\chi_{\infty, j}) \\
 & \rightarrow  & \RH_\mathrm{ct}^{d_{n-1,\infty}}( \GL_{n-1}(\rk_\infty)^0; (\Pi_\infty\widehat{\otimes} \Sigma_\infty\otimes  \chi_{\infty, j}) \otimes (F_\xi^\vee\otimes \delta_{ j}^\vee)) \otimes 
   \frak{O}_{n-1,\infty} \nonumber \\
& \rightarrow  & \RH_\mathrm{ct}^{d_{n-1,\infty}}(\GL_{n-1}(\rk_\infty)^0;  \frak{M}_{n-1,\infty}^*)\otimes \frak{O}_{n-1,\infty}=\BC, \nonumber 
\end{eqnarray}
where  the first arrow is the restriction of cohomology composed with the cup product, and the last arrow is the map induced by the linear functional
\[
\oZ^\circ(\cdot, \frac{1}{2}+j, \chi_\infty)  \otimes \phi_{\xi,j}:(\Pi_\infty\widehat\otimes \Sigma_\infty\otimes \chi_{\infty, j})\otimes (F_\xi^\vee \otimes \delta_{ j}^\vee) \rightarrow \frak{M}_{n-1,\infty}^*.
\]

 \begin{prp}\label{prop: main diagram} (\cf \cite[(5.3)]{Mah} and \cite[Section 4.6]{Jan5})
         Let the notations and  assumptions be as above.  Then the diagram
         $$\begin{CD}
              \Pi_f\otimes \Sigma_f\otimes \chi_{f, j}  \otimes  \frak{M}_{n-1,f}  \otimes \CH(\Pi_\infty)\otimes \CH(\Sigma_\infty)\otimes \oH(\chi_{\infty, j}) @> \oZ^\circ(\,\cdot\,, \frac{1}{2}+j, \chi_f) \otimes \CP_{\infty,j}>> \BC \\
           @V \iota_{\rm can}  VV @V \oL(\frac{1}{2}+j, \Pi\times \Sigma \times \chi) V V  \\
           \CH(\Pi)\otimes \CH(\Sigma)\otimes \oH(\chi_j) \otimes \frak{M}_{n-1, f} @>\CP_j>> \C
         \end{CD}$$
commutes, where the left vertical arrow $\iota_{\rm can}$ is the natural isomorphism. 
      \end{prp}
      
     \begin{proof}
     Define $\frak{q}_{n,\infty} := \left( \frak{gl}_{n}(\rk_\infty)/ (\BR \oplus \frak{k}_{n,\infty})\right)\otimes_\BR \BC$. We have a map
     \[
    ( \wedge^{b_{n,\infty}}\frak{q}_{n,\infty})^* \otimes (\wedge^{b_{n-1,\infty}}\frak{q}_{n-1,\infty})^* \to \omega_{n-1,\infty}^* = \wedge^{d_{n-1,\infty}}\left((\frak{gl}_{n-1}(\rk_\infty)/\frak{k}_{n-1,\infty})\otimes_\BR\BC\right)^*
     \]
     induced by restriction. By the identification of  continuous cohomology and relative Lie algebra cohomology (\cite[Theorem 6.1]{HM}), as well as the explicit determination of the  relative Lie algebra cohomology (\cite[Proposition 9.4.3]{Wa1}), we have that 
     \[
     \CH(\Pi_\infty) =\big((\wedge^{b_{n,\infty}}\frak{q}_{n,\infty})^* \otimes \Pi_\infty \otimes F_\mu^\vee\big)^{K_{n,\infty}^0} \otimes \widetilde{\frak{O}}_{n,\infty}
     \]
     and 
     \[
       \CH(\Sigma_\infty) =\big((\wedge^{b_{n-1,\infty}}\frak{q}_{n-1,\infty})^*  \otimes \Sigma_\infty \otimes F_\nu^\vee\big)^{K_{n-1,\infty}^0}\otimes \widetilde{\frak{O}}_{n-1,\infty}.
     \]
     By definition of $\CP_{\infty, j}$, the top horizontal arrow of the diagram is identified with the composition of
     \begin{eqnarray*}
&& \Pi_f\otimes \Sigma_f\otimes \chi_{f, j} \otimes \frak{M}_{n-1,f} \otimes \big((\wedge^{b_{n,\infty}}\frak{q}_{n,\infty})^* \otimes \Pi_\infty \otimes F_\mu^\vee\big)^{K_{n,\infty}^0}\otimes  \widetilde{\frak{O}}_{n,\infty}  \\ 
&& \otimes \big((\wedge^{b_{n-1,\infty}}\frak{q}_{n-1,\infty})^*  \otimes \Sigma_\infty \otimes F_\nu^\vee\big)^{K_{n-1,\infty}^0} \otimes \widetilde{\frak{O}}_{n-1,\infty} \otimes \delta_{j}^\vee \otimes \chi_{\infty, j} \\
           &\xrightarrow{\textrm{restriction}} &   \Pi_f\otimes\Sigma_f \otimes  \chi_{f, j}\otimes \frak{M}_{n-1, f} \otimes \omega_{n-1,\infty}^*\otimes \Pi_\infty \widehat{\otimes} \Sigma_\infty \otimes F_\xi^\vee  \otimes \frak{O}_{n-1,\infty}  \otimes\delta_{ j}^\vee \otimes  \chi_{\infty, j} 
            \\
         &=& (\Pi_f\otimes\Sigma_f \otimes  \chi_{f, j}\otimes \frak{M}_{n-1,f})  \otimes ( \Pi_\infty \widehat{\otimes} \Sigma_\infty\otimes \chi_{\infty, j} \otimes \frak{M}_{n-1,\infty})\otimes (F_\xi^\vee\otimes \delta_{j}^\vee)   \\
         &\rightarrow & \BC,
                  \end{eqnarray*}
                  where the last map is given by
                  \[
                   \oZ^\circ(\cdot, \frac{1}{2}+j, \chi_f)\otimes \oZ^\circ(\cdot, \frac{1}{2}+j, \chi_\infty)  \otimes \phi_{\xi,j}.
                  \] 
                  Using fast decreasing differential forms as in \cite[Section 5.6]{Bo}, the bottom arrow of the  diagram is identified with the composition of
                  \begin{eqnarray*}
&&\big((\wedge^{b_{n,\infty}}\frak{q}_{n,\infty} )^* \otimes \Pi \otimes F_\mu^\vee\big)^{K_{n,\infty}^0} \otimes \widetilde{\frak{O}}_{n,\infty}  \\
&&  \otimes \big((\wedge^{b_{n-1,\infty}}\frak{q}_{n-1,\infty} )^*  \otimes \Sigma \otimes F_\nu^\vee\big)^{K_{n-1,\infty}^0} \otimes \widetilde{\frak{O}}_{n-1,\infty}  \otimes \delta_{ j}^\vee \otimes \chi_j \otimes \frak{M}_{n-1, f}\\
           &\xrightarrow{\textrm{restriction}} & \omega_{n-1,\infty}^*\otimes \Pi \widehat{\otimes} \Sigma \otimes F_\xi^\vee   \otimes \frak{O}_{n-1,\infty}\otimes \delta_j^\vee\otimes   \chi_j  \otimes \frak{M}_{n-1,f } \\
       &=&  (  \Pi \widehat{\otimes} \Sigma \otimes \chi_j \otimes \frak{M}_{n-1}) \otimes  (F_\xi^\vee \otimes \delta_{ j}^\vee) \\
         &\xrightarrow{\oZ(\cdot, \frac{1}{2}+j, \chi)\otimes \phi_{\xi, j}} & \BC.
                  \end{eqnarray*}
                  The proposition then follows from Proposition \ref{measure}.
     \end{proof}

\subsection{ Two commutative diagrams} 

For every $\sigma\in \Aut(\C)$, we note that the infinite part of $({}^\sigma \chi)_j$ coincides with $\chi_{\infty, j}$. 
Denote the corresponding modular symbol at infinity by 
\[
{}^\sigma\CP_{\infty, j}: \CH({}^\sigma\Pi_\infty) \otimes \CH({}^\sigma\Sigma_\infty) \otimes \RH(\chi_{\infty, j})  \to \BC,
\]
and introduce the normalized modular symbol at infinity
\be\label{norminf}
{}^\sigma \CP^\circ_{\infty, j} :=  \Omega'_{\mu, \nu, j}\cdot {}^\sigma\CP_{\infty, j}, \quad\textrm{where }\ 
\Omega'_{\mu, \nu, j} := \prod_{v|\infty} \Omega'_{\mu_v, \nu_v, j}.
\ee
In particular, we have the 
 normalized modular symbol at infinity
\[
\CP^\circ_{\infty, j} :=  \Omega'_{\mu, \nu, j}\cdot \CP_{\infty, j}.
\]

 As in   \eqref{sig-inf}, we have  a $\sigma$-linear isomorphisms 
     \[
      \sigma: \CH(\Pi_\infty) \to \CH({}^\sigma \Pi_\infty)
    \]
      such that
      \[
      \sigma(\kappa_{\mu,\varepsilon}) = \kappa_{{}^\sigma\!\mu,\varepsilon}\quad \textrm{for all $\varepsilon\in \widehat{\pi_0(\rk_\infty^\times)}$ that occurs in $\CH(\Pi_\infty)$}.
      \]
We have 
a similar $\sigma$-linear isomorphism 
     \[
      \sigma: \CH(\Sigma_\infty) \to \CH({}^\sigma \Sigma_\infty), 
    \]
    as well as a $\sigma$-linear isomorphism 
    \[
     \sigma: \RH(\chi_{\infty, j}) \rightarrow \RH(\chi_{\infty, j}) 
      \]
      such that $\sigma(1)=1$.
      By tensor product, we get a $\sigma$-linear isomorphism
      \[
      \sigma:  \CH(\Pi_\infty) \otimes \CH(\Sigma_\infty) \otimes \RH(\chi_{\infty, j}) \rightarrow  \CH({}^\sigma\Pi_\infty) \otimes \CH({}^\sigma\Sigma_\infty) \otimes \RH(\chi_{\infty, j}).
      \]
\begin{prp}

For all $\sigma\in \Aut(\C)$, the diagram
\be \label{inf-sig}
 \begin{CD}
            \CH(\Pi_\infty) \otimes \CH(\Sigma_\infty) \otimes \RH(\chi_{\infty, j}) 
                  @> \CP_{\infty, j}^\circ >>   \BC  \\
            @ V \sigma VV         @ VV \sigma V \\
            \CH({}^\sigma\Pi_\infty) \otimes \CH({}^\sigma\Sigma_\infty) \otimes \RH(\chi_{\infty, j})  @> {}^\sigma\CP_{\infty, j}^{ \circ} >> \BC
  \end{CD}
\ee
commutes. 
    
\end{prp}
\begin{proof}
      Let $\Pi_{0_{n,\infty}}:=\widehat\otimes_{v|\infty}\Pi_{0_{n, \rk_v}}$ and $\Sigma_{0_{n-1, \infty}}:=\widehat \otimes_{v|\infty}\Sigma_{0_{n-1, \rk_v}}$  be the cohomological representations (as in Section \ref{sec21}) of $\GL_n(\rk_\infty)$ and $\GL_{n-1}(\rk_\infty)$  that have trivial coefficient systems and respectively have the same central characters 
  as that of $F_{\mu}^\vee\otimes \Pi_\infty$ and $F_{\nu}^\vee \otimes \Sigma_\infty$.  

Applying Theorem \ref{thmap} for all $v|\infty$, we obtain a commutative diagram
 \be \label{inf-red}
 \begin{CD}
            \CH(\Pi_\infty) \otimes \CH(\Sigma_\infty) \otimes \RH(\chi_{\infty, j}) 
                  @>\Omega'_{\mu,\nu, j} \cdot \CP_{\infty, j}>>   \BC  \\
            @A  \jmath_\mu \otimes \jmath_\nu \otimes{\rm id}   A  A         @ | \\
           \CH(\Pi_{0_{n, \infty}}) \otimes \CH(\Sigma_{0_{n-1, \infty}}) \otimes \RH(\chi_{\infty, j}) @> \CP_{\infty, 0} >> \BC.  \\
  \end{CD}
\ee 
This easily implies the proposition. 
\end{proof}

 Pick an element $y=(y_v)_{v\nmid \infty}\in \A_f^\times$ such that  
\[
\frak{c}(\psi_v)= y_v \cdot \frak{c}(\chi_v)\quad \textrm{for all $v\nmid \infty$.}
\]
Define the Gauss sum  \be\label{gauss-sum000}
    \CG(\chi):=\CG(\chi,\psi,y):=\prod_{v\nmid \infty} \CG(\chi_v, \psi_v, y_v),
     \ee
     where $\CG(\chi_v, \psi_v, y_v)$ is the local Gauss sum given by \eqref{gauss}. 
    Similarly, pick an element $y'=(y'_v)_{v\nmid \infty}\in \A_f^\times$ such that  
\[
\frak{c}(\psi_v)= y'_v \cdot \frak{c}(\chi_{\Sigma_v})\quad \textrm{for all $v\nmid \infty$},
\]
and define the Gauss sum  \be\label{gauss-sum002}
    \CG(\chi_\Sigma):= \CG(\chi_{\Sigma},\psi,y'):=\prod_{v\nmid \infty} \CG(\chi_{\Sigma_v}, \psi_v, y'_v).
     \ee
   Here we write $\Sigma_f=\otimes_{v\nmid \infty}' \Sigma_v$ as usual, and  $\chi_\Sigma$ and  $\chi_{\Sigma_v}$ respectively denote the central characters of $\Sigma$ and $\Sigma_v$. More generally, we have the Gauss sums
   \[
    \CG({}^\sigma\chi):=\CG({}^\sigma\chi,\psi,y)\quad\textrm{and}\quad \CG(\chi_{{}^\sigma\Sigma}):=\CG(\chi_{{}^\sigma\Sigma},\psi,y'),
   \]
where $\chi_{{}^\sigma \Sigma}$ denotes the central character of ${}^\sigma \Sigma$.

Similar to \eqref{sigmal00},  for all $s\in \C$ we have a $\sigma$-linear isomorphism 
\[
   \sigma: \Pi_f \otimes \Sigma_f\otimes \chi_{f, s-\frac{1}{2}} \otimes \frak{M}_{n-1,f}\rightarrow  {}^\sigma  \Pi_f \otimes {}^\sigma\Sigma_f \otimes {}^\sigma (\chi_{f, s-\frac{1}{2}}) \otimes \frak{M}_{n-1, f}.
\]
Note that ${}^\sigma (\chi_{f, s-\frac{1}{2}})=({}^\sigma \chi)_{f, s-\frac{1}{2}}$ when $s\in \frac{1}{2}+\Z$. 

     \begin{prp} \label{Af-rel}
    For all $s_0\in \frac{1}{2}+\BZ$ and $\sigma\in \Aut(\BC)$, the  diagram  \[
 \begin{CD}
          \Pi_f \otimes \Sigma_f\otimes \chi_{f, s_0-\frac{1}{2}} \otimes \frak{M}_{n-1,f}
                  @> \CG(\chi_{\Sigma}) \cdot \CG(\chi)^{\frac{n(n-1)}{2}}\cdot \oZ^\circ(\cdot, s_0, \chi_f) >> \BC \\
            @V   \sigma VV            @VV  \sigma  V\\
        {}^\sigma  \Pi_f \otimes {}^\sigma\Sigma_f \otimes {}^\sigma (\chi_{f, s_0-\frac{1}{2}}) \otimes \frak{M}_{n-1, f}@> \CG(\chi_{{}^\sigma\Sigma})\cdot \CG({}^\sigma\chi)^{\frac{n(n-1)}{2}}\cdot \oZ^\circ(\cdot, s_0, {}^\sigma\chi_f) >> \C \\
  \end{CD}
\]
commutes.
\end{prp}

     \begin{proof}
    Write 
     $
     \Pi_f= \otimes'_{v\nmid\infty}\Pi_v$ as usual. 
     By the uniqueness of Whittaker functionals, we  write 
     $\lambda_f=\otimes_{v\nmid\infty} \lambda_v$, $\lambda_f' = \otimes_{v\nmid\infty} \lambda'_v$ and assume that 
     \[
     \lambda_v(e_v) = \lambda'_v (e_v') =1
     \]
     for all but finitely many  $v\nmid\infty$ such that $\Pi_v$ and $\Sigma_v$ are unramified, where $e_v\in \Pi_v$ and $e_v'\in \Sigma_v$ are the spherical vectors used in the definition of the restricted tensor products 
     $\Pi_f$ and $\Sigma_f$. For places $v$ as above, if moreover $\chi_v$ is unramified and $\psi_v$ has conductor $\CO_{\rk_v}$, then it is known that  (see \cite[Proposition 2.4]{JS1})
     \[
     \oZ^\circ(e_v\otimes e_v' \otimes m_{n-1,\rk_v}^\circ, s, \chi_v)=1,
     \]
where $m_{n-1, \rk_v}^\circ\in \frak{M}_{n-1, \rk_v}$ is the Haar measure on $\GL_{n-1}(\rk_v)$ such that a maximal open compact subgroup has total volume 1.
 The proposition then follows from Proposition \ref{nonarchrel}.  
     \end{proof}

In analogy to \eqref{norminf}, for the finite part we introduce 
\[
{}^\sigma\CP^{\circ}_{f, j} := \CG(\chi_{{}^\sigma\Sigma})\cdot \CG({}^\sigma\chi)^{\frac{n(n-1)}{2}}\cdot \oZ^\circ(\cdot, \frac{1}{2}+j, {}^\sigma\chi_f).
\]
Specifically, we have 
\[
\CP^{\circ}_{f, j} := \CG(\chi_\Sigma)\cdot \CG(\chi)^{\frac{n(n-1)}{2}}\cdot \oZ^\circ(\cdot, \frac{1}{2}+j, \chi_f).
\]
Then Proposition \ref{Af-rel} can be rephrased as the following commutative diagram 
\be \label{f-sig}
 \begin{CD}
           \Pi_f\otimes \Sigma_f \otimes \chi_{f, j} \otimes \frak{M}_{n-1, f}
                  @> \CP_{f, j}^\circ >>   \BC  \\
            @ V \sigma VV         @ VV \sigma V \\
            {}^\sigma\Pi_f \otimes {}^\sigma\Sigma_f \otimes {}^\sigma( \chi_{f, j} )\otimes \frak{M}_{n-1, f}@> {}^\sigma \CP_{f, j}^{ \circ} >> \BC.
  \end{CD}
\ee

      
      


     \subsection{Proof of Theorem \ref{thm: global period relation}}


As in \eqref{pj}, we have the modular symbol map
\[
{}^\sigma \CP_j:  \CH({}^\sigma\Pi)\otimes \CH({}^\sigma\Sigma )\otimes \RH({}^\sigma\chi_j) \otimes \frak{M}_{n-1,f } \to \BC.
\]
Put
\[
{}^\sigma \oL^*_j :=  \frac{\oL(\frac{1}{2}+j, {}^\sigma\Pi \times {}^\sigma\Sigma \times {}^\sigma\chi )} { \Omega'_{\mu,\nu, j} \cdot \CG(\chi_{{}^\sigma\Sigma})\cdot \CG({}^\sigma\chi)^{\frac{n(n-1)}{2}}}. 
\]
Then by Proposition \ref{prop: main diagram} the diagram
\[
\begin{CD}
               {}^\sigma\Pi_f \otimes {}^\sigma\Sigma_f\otimes {}^\sigma(\chi_{f, j}) \otimes \frak{M}_{n-1,f} \otimes \CH({}^\sigma\Pi_\infty)\otimes \CH({}^\sigma \Sigma_\infty )\otimes \oH(\chi_{\infty, j}) @>{}^\sigma\CP^{ \circ}_{f,j}\otimes  {}^\sigma\CP^{\circ}_{\infty, j} >> \BC \\
           @V \iota_{\rm can}  VV @ VV {}^\sigma\oL^*_j  V  \\
           \CH({}^\sigma\Pi)\otimes \CH({}^\sigma\Sigma )\otimes \RH({}^\sigma\chi_j) \otimes \frak{M}_{n-1,f }  @> {}^\sigma \CP_j  >> \C
         \end{CD}
\]
commutes. 

We are now ready to prove Theorem \ref{thm: global period relation}. It is clear that \eqref{eq:mainthm0} is a consequence of 
\eqref{eq:galois}, and we will prove the latter.  To save space,  denote the subspaces of  $\pi_0(\rk_\infty^\times)$-fixed vectors in   the two spaces in the left vertical arrow of  the last diagram by $ \CH({}^\sigma\Pi, {}^\sigma\Sigma, {}^\sigma\chi, j)_{\rm loc}$ and $ \CH({}^\sigma\Pi, {}^\sigma\Sigma, {}^\sigma\chi, j)_{\rm glob}$ respectively, so that the last diagram reads 
\be \label{diag1}
\CD
           \CH({}^\sigma\Pi, {}^\sigma\Sigma, {}^\sigma\chi, j)_{\rm loc} @> {}^\sigma \CP^{ \circ}_{f,j}\otimes   {}^\sigma \CP^{ \circ}_{\infty, j}  >> \BC \\
           @V \iota_{\rm can} V V @ VV {}^\sigma\oL^*_j V \\
            \CH({}^\sigma\Pi, {}^\sigma\Sigma, {}^\sigma\chi, j)_{\rm glob} @>  {}^\sigma \CP_j  >> \BC.
         \endCD
        \ee
        
                By \eqref{per-diag}, we have a commutative diagram 
         \be \label{diag2}
         \CD
           \CH(\Pi, \Sigma, \chi, j)_{\rm loc}  @> \sigma >> \CH({}^\sigma\Pi, {}^\sigma\Sigma, {}^\sigma\chi, j)_{\rm loc}  \\
           @V  \Omega_{(j)}^{-1} \cdot \iota_{\rm can} VV @VV  {}^\sigma\Omega_{(j)}^{-1}\cdot\iota_{\rm can} V \\
           \CH(\Pi, \Sigma, \chi, j)_{\rm glob} @> \sigma >>  \CH({}^\sigma\Pi, {}^\sigma\Sigma, {}^\sigma\chi, j)_{\rm glob}, 
         \endCD
         \ee
         where the top horizontal arrow is the tensor product of the left vertical arrows in \eqref{inf-sig} and \eqref{f-sig}, and we write for short
    \[
         \Omega_{(j)}:= \Omega_{\varepsilon_n}(\Pi)\cdot\Omega_{\varepsilon_{n-1}}(\Sigma)\quad\textrm{ and }\quad  {}^\sigma \Omega_{(j)}:= \Omega_{\varepsilon_n}({}^\sigma\Pi)\cdot\Omega_{\varepsilon_{n-1}}({}^\sigma\Sigma).
         \]

       It is well-known that the global modular symbol is $\Aut(\BC)$-equivarient (see \cite[Proposition 3.14]{Rag1}), that is, the following diagram
       \be \label{diag3}
          \CD
            \CH(\Pi, \Sigma, \chi, j)_{\rm glob} @> \CP_j >> \BC \\
           @V \sigma VV @VV \sigma V  \\
           \CH({}^\sigma\Pi, {}^\sigma\Sigma, {}^\sigma\chi, j)_{\rm glob}  @>{}^\sigma\CP_j>> \BC
         \endCD
       \ee 
        commutes. 
        
Since $\Omega_{\mu,\nu,j}$ and $\Omega'_{\mu,\nu,j}$ only differ by a sign, \eqref{eq:galois} is equivalent to the equation 
\[
\sigma\left(\frac{\oL^*_j}{\Omega_{(j)}}\right) = \frac{{}^\sigma\oL^*_j}{{}^\sigma\Omega_{(j)}},
\]
which amounts to the commutativity of the diagram
\be \label{diag4}
          \CD
            \BC @> \sigma >> \BC \\
           @V \frac{\oL^*_j}{\Omega_{(j)}}\ VV @VV  \frac{{}^\sigma\oL^*_j}{{}^\sigma\Omega_{(j)}} V  \\
           \BC  @> \sigma >> \BC.
         \endCD
\ee
Here 
\[
\oL^*_j :=  \frac{\oL(\frac{1}{2}+j, \Pi \times \Sigma \times \chi )} { \Omega'_{\mu,\nu, j} \cdot \CG(\chi_{\Sigma})\cdot \CG(\chi)^{\frac{n(n-1)}{2}}}.
\]

 The commutative diagrams \eqref{diag1}, \eqref{diag2} and \eqref{diag3}, together with \eqref{inf-sig} and \eqref{f-sig}, give us the following diagram
         \[
         \xymatrix{
          \CH(\Pi, \Sigma, \chi, j)_{\rm loc} \ar[rd]^{\sigma} \ar[ddd]_{\Omega_{(j)}^{-1}\cdot\iota_{\rm can}} \ar[rrr]^{  \CP^{ \circ}_{f,j}\otimes  \CP^{ \circ}_{\infty, j} } 
  &&& \BC \ar^\sigma[rrd] \ar^{ \sigma}[rrd] \ar[ddd]^{\frac{\oL^*_j}{\Omega_{(j)}}}|!{[dll];[dr]}\hole & \\
    &\CH({}^\sigma\Pi, {}^\sigma\Sigma, {}^\sigma\chi, j)_{\rm loc}  \ar[ddd]_{{}^\sigma\Omega_{(j)}^{-1}\cdot\iota_{\rm can}}  \ar[rrrr]^(0.4){ {}^\sigma\CP^{\circ}_{f,j}\otimes  {}^\sigma\CP^{ \circ}_{\infty, j} } 
    &&& & \BC  \ar[ddd]^{ \frac{{}^\sigma\oL^*_j}{{}^\sigma\Omega_{(j)}}}\\
    &&& & \\
     \CH(\Pi, \Sigma, \chi, j)_{\rm glob} \ar[rd]_{\sigma} \ar[rrr]^(0.6){\CP_j} |!{[ur];[dr]}\hole
    && & \BC \ar[drr]^{\sigma}  & \\
    &  \CH({}^\sigma\Pi,{}^\sigma\Sigma, {}^\sigma\chi, j)_{\rm glob} \ar[rrrr]^{{}^\sigma\CP_j} 
    &&& & \BC,
}
\]
where all squares are commutative except \eqref{diag4}. This forces \eqref{diag4} to be commutative as well. 
 This proves \eqref{eq:galois}, hence  finishes the proof of Theorem \ref{thm: global period relation}.

\section*{Acknowledgements}

D. Liu was supported in part by the Natural Science Foundation of Zhejiang Province (No. LZ22A010006) and the National Natural Science Foundation of 
China (No. 12171421). B. Sun was supported in part by National Key $\textrm{R}\,\&\,\textrm{D}$ Program of China (No. 2020YFA0712600).

\end{document}